\documentclass[aos,MSNbibl,nameyear,dvips]{arximspdf}
\usepackage{graphicx}
\doi{10.1214/13-AOS1128} %kopijuoti is PTS
\volume{41}
\issue{1}
\pubyear{2013}
\firstpage{1816}
\lastpage{1864}

\makeatletter
\newcommand{\rrVert}{\Vert}
\newcommand{\llVert}{\Vert}
\def\cal{\mathcal}
\newtheorem{theorem}{Theorem}
\newtheorem{lemma}[theorem]{Lemma}
\newtheorem{proposition}[theorem]{Proposition}
\newproclaim{remark}{Remark}

\renewcommand{\hat}{\widehat}

\newcommand{\boeta}{\bolds{\eta}}

\newcommand{\cF}{{\cal F}}

\newcommand{\bA}{\mathbf{A}}
\newcommand{\bB}{\mathbf{B}}
\newcommand{\bG}{\mathbf{G}}
\newcommand{\bI}{\mathbf{I}}
\newcommand{\bQ}{\mathbf{Q}}
\newcommand{\bU}{\mathbf{U}}
\newcommand{\bV}{\mathbf{V}}
\newcommand{\bX}{\mathbf{X}}
\newcommand{\bY}{\mathbf{Y}}

\newcommand{\be}{\mathbf{e}}
\newcommand{\bmu}{\bolds{\mu}}
\newcommand{\btau}{\bolds{\tau}}
\newcommand{\bsigma}{\bolds{\sigma}}
\newcommand{\bUpsilon}{\bolds{\Upsilon}}
\newcommand{\bvarepsilon}{\bolds{\varepsilon}}
\newcommand{\bGamma} {\bolds{\Gamma}}
\newcommand{\bD}{\mathbf{D}}
\newcommand{\beps}{\bvarepsilon}
\newcommand{\bPhi}{\bolds{\Phi}}

\makeatother

\begin{document}
\begin{frontmatter}

\title{Optimal sparse volatility matrix estimation for high-dimensional
It\^{o} processes with measurement~errors}
\runtitle{Optimal large sparse volatility matrix estimation}

\begin{aug}
\author[A]{\fnms{Minjing} \snm{Tao}\ead[label=e1]{minjing@stat.wisc.edu}},
\author[A]{\fnms{Yazhen} \snm{Wang}\corref{}\thanksref{t1}\ead[label=e2]{yzwang@stat.wisc.edu}}
\and
\author[B]{\fnms{Harrison H.} \snm{Zhou}\thanksref{t2}\ead[label=e3]{huibin.zhou@yale.edu}}
\thankstext{t1}{Supported in part by NSF Grants DMS-10-5635 and DMS-12-65203.}
\thankstext{t2}{Supported in part by NSF Career Award DMS-0645676 and NSF FRG Grant DMS-08-54975.}
\runauthor{M. Tao, Y. Wang and H. H. Zhou}
\affiliation{University of Wisconsin-Madison, University of
Wisconsin-Madison and~Yale~University}
\address[A]{M. Tao\\
Y. Wang\\
Department of Statistics\\
University of Wisconsin-Madison\\
Madison, Wisconsin 53706\\
USA\\
\printead{e1}\\
\phantom{E-mail:\ }\printead*{e2}} %adresu isvedimo komanda gale!
\address[B]{H. H. Zhou\\
Department of Statistics\\
Yale University\\
New Haven, Connecticut 06511\\
USA\\
\printead{e3}}
\end{aug}

% HISTORY:
\received{\smonth{11} \syear{2012}}
\revised{\smonth{2} \syear{2013}}

% ABSTRACT
%
\begin{abstract}
Stochastic processes are often used to model complex scientific
problems in
fields ranging from biology and finance to engineering and physical science.
This paper investigates rate-optimal estimation of the volatility
matrix of a high-dimensional It\^{o} process
observed with measurement errors at discrete time points.
The minimax rate of convergence is established for estimating sparse
volatility matrices. By combining the multi-scale
and threshold approaches we construct a volatility matrix estimator to
achieve the optimal convergence rate.
The minimax lower bound is derived by considering a subclass of It\^{o}
processes for which the minimax
lower bound is obtained through a novel equivalent model of covariance
matrix estimation for independent but
nonidentically distributed observations and through a delicate
construction of the least favorable
parameters. In addition, a simulation study was conducted to test the
finite sample performance of the
optimal estimator, and the simulation results were found to support the
established asymptotic theory. %our theoretical findings.

%a $p$-dimensional It\^{o} process is observed with measurement errors
%at $n$ distinct time points, and our goal is to
%estimate the volatility matrix of the process. We establish the
%minimax theory for estimating large sparse volatility matrices under
%matrix
%norm as both $n$ and $p$ go to infinity. The theory shows that the
%optimal convergence rate depends on $n$ and $p$
%through $n^{-1/4} \sqrt{\log p}$ and a volatility matrix estimator is
%explicitly constructed to achieve the optimal convergence rate.
%diffusion processes play an important role in many modern scientific
%studies.
\end{abstract}

% KEYWORDS
% Pirmas kwd is didziosios raides
%

\end{frontmatter}

%s1 #&#
\section{Introduction}\label{sec1}

Modern scientific studies in fields ranging from biology and finance to
engineering and physical science often
need to model complex dynamic systems where it is essential to
incorporate internally or externally originating
random fluctuations in the systems [\citet{AtSMykZha05},
\citet{MueAnd06} and \citet{Whi95}].
%As covariance matrix estimation is widely studied in statistical
%problems with i.i.d. observations,
Continuous-time diffusion processes, or more generally, It\^o
processes, are frequently employed to model such complex dynamic systems.
Data collected in the studies are treated as the processes observed at
discrete time points with possible
noise contamination. For example, the prices of financial assets are
usually modeled by It\^o processes, and
the price data observed at high-frequencies are contaminated by market
microstructure noise. In this paper we
investigate estimation of the volatilities of
the It\^o processes based on noisy data.
%In financial applications such as risk management, portfolio
%allocation, and option pricing, we need to estimate the volatilities
%of the processes based on the noisy data.

Several volatility estimation methods have been developed in the past
several years. For estimating a univariate integrated
volatility, popular estimators include two-scale realized volatility
[\citet{ZhaMykAtS05}], multi-scale realized volatility
[\citet{Zha06} and \citet{FanWan07}], realized kernel volatility [\citet{Baretal08}] and pre-averaging based realized
volatility [\citet{Jacetal09}]. For estimating a bivariate integrated
co-volatility, common methods are the previous-tick approach
[\citet{Zha11}], the refresh-time scheme and realized kernel volatility
[\citet{Baretal11}],
the generalized synchronization scheme [\citet{AtSFanXiu10}]
and the pre-averaging
approach [\citet{ChrKinPod10}]. %, Kinnebrock, and Podolskij
%(2010))
Optimal volatility and co-volatility estimation has been investigated
in the parametric or nonparametric setting
[\citet{AtSMykZha05}, \citet{Bib11}, Gloter
and Jacod (\citeyear{GloJac01N1,GloJac01N2}), \citet{Rei11} and \citet{Xiu10}].
%For the univariate case, Gloter and Jacod (\citeyear{GloJac01N1}, b) investigated
%optimal estimation of a scalar volatility under a diffusion model with
%noise contamination.
These works are for estimating scalar volatilities or volatility
matrices of small size.
\citet{WanZou10} and \citet{Taoetal11} studied the problem of
estimating a large sparse
volatility matrix based on noisy high-frequency financial data.
\citet{FanLiYu12} employed a large volatility matrix estimator based on
high-frequency data for portfolio allocation.
The large volatility matrix estimation is a high-dimensional extension
of the univariate case. It can be also considered
as a generalization of large covariance matrix estimation for i.i.d.
data to volatility matrix estimation for dependent
data with measurement errors. Despite recent progress on volatility
matrix estimation, there has been remarkably little
fundamental theoretical study on optimal estimation of large volatility
matrices.
Consistent estimation of large matrices based on high-dimensional data
usually requires some sparsity, and the sparsity
may naturally result from appropriate formulation of some
low-dimensional structures in the high-dimensional data.
%For large matrix estimation with high dimensional data, there are
%usually some low-dimensional structures in
%the high dimensional data, and with appropriate formalization the
%low-dimensional structures will give
%rise to sparsity.
For example, in large volatility matrix estimation with high-frequency
financial data sparsity means that a relatively
small number of market factors play a dominate role in driving
volatility movements and capturing the market risk.
In this paper we establish the optimal rate of convergence for large
volatility matrix estimation under various
matrix norms over a wide range of classes of sparse volatility
matrices. We expect that our work will stimulate further
theoretical and methodological research as well as more application
orientated study on large volatility matrix estimation.

%More
Specifically we consider the problem of estimating the sparse
integrated volatility matrix for a $p$-dimensional
It\^o process observed with additive noises at $n$ equally spaced
discrete time points. The minimax upper bound is obtained by
constructing a new procedure through a combination of the multi-scale
and threshold approaches and by studying its
risk properties. We first construct a multi-scale volatility matrix
estimator and show that its elements obey
subGaussian tails with a convergence rate $n^{-1/4}$. Then we threshold
the constructed estimator to obtain a
threshold volatility matrix estimator and derive its convergence rate.
The upper bound depends on $n$ and $p$ through $n^{-1/4} \sqrt{\log p}$.

A key step in obtaining the optimal rate of convergence is the
derivation of the minimax lower bound for the high-dimensional
It\^o process with measurement errors. We succeed in establishing the
risk lower bound in three steps. First we select a
particular subclass of It\^o processes with a zero drift and
a constant volatility matrix so that the volatility matrix estimation
problem becomes a covariance matrix estimation problem
where the observed data are dependent and have measurement errors;
second, take a special transformation of the observations
to convert the problem into a new covariance matrix estimation problem
where the observed data have no measurement errors and
are independent but not identically distributed, with covariance
matrices equal to the constant volatility matrix plus an
identity matrix multiplying by a shrinking factor depending on the
sample size $n$; third, adopt the minimax lower bound
technique developed in \citet{CaiZho} for sparse covariance matrix
estimation based on i.i.d. data to establish a minimax
lower bound for independent but nonidentically distributed
observations. The minimax lower bound matches the upper bound obtained
by the new procedure up to a constant factor, and thus the upper bound
is rate-optimal.

The volatility matrix estimation is closely related to large covariance
matrix estimation which received lots of attentions recently in the literature.
% as well as in other areas such as compressive sensing and financial
%econometrics. There is extensive and vibrant research on large
%covariance matrix estimation in past several years.
While the covariance matrix %estimation and principal component analysis
plays a key role in statistical analysis, its classic estimation
procedures, like the sample
covariance matrix estimator, may behave very poorly when the matrix
size is comparable to or exceeds the sample
size. To overcome the curse of dimensionality, various regularization
techniques have been developed for estimation of large covariance
matrices in recent years.
\citet{WuPou03} explored nonparametric estimation of large
covariance matrices by local stationarity.
\citet{LedWol04} proposed to boost diagonal elements and downgrade
off-diagonal elements of the sample covariance matrix estimator.
\citet{Huaetal06} used a penalized likelihood method to estimate
large covariance matrices.
\citet{YuaLin07} considered large covariance matrix estimation in a
Gaussian graph model.
Bickel and Levina (\citeyear{BicLev08N1,BicLev08N2}) developed regularization methods by
banding or thresholding the sample covariance matrix estimator when the
matrix size is comparable to the sample size.
\citet{ElK08} employed a graph model approach to characterize sparsity
and investigated consistent estimation of large covariance matrices.
\citet{FanFanLv08} utilized factor models for estimating large
covariance matrices.
\citet{JohLu09} studied consistent estimation of leading
principal components in principal component analysis.
\citet{LamFan09} established sparsistency and convergence rates for
large covariance matrix estimation.
\citet{CaiZhaZho10} and \citet{CaiZho} studied minimax
estimation of covariance matrices when both sample size and matrix size
are allowed to go to
infinity and derived optimal convergence rates for estimating decaying or
sparse covariance matrices.

The rest of the paper proceeds as follows. Section~\ref{sec2} presents the model
and the data and
constructs volatility matrix estimators. Section~\ref{secasymp} establishes the
asymptotic theory under sparsity for
the constructed matrix estimators as both sample size and matrix size
go to infinity.
Section~\ref{optimalsec} derives the minimax lower bound for estimating a large sparse
volatility matrix
and shows that the threshold volatility matrix estimator asymptotically
achieves the minimax lower
bound. Thus combining results in Sections~\ref{secasymp} and~\ref{optimalsec} together, we
establish the optimality for
large sparse volatility matrix estimation. Section~\ref{sec5} features a
simulation study to illustrate the finite sample
performances of the volatility matrix estimators. To facilitate the
reading we relegate all proofs to
Section~\ref{sec6} and two \hyperref[app]{Appendix} sections, where we first provide the main
proofs of the theorems in Section~\ref{sec6} and
then collect additional technical proofs in the two appendices.

%s2 #&#
\section{Volatility matrix estimation}\label{sec2}
%s2.1 #&#
\subsection{The model set-up}
\label{secmodel}
Suppose that $\bX(t)=(X_1(t),\ldots,X_p(t))^T$ is an It\^{o} process
following the model %a continuous-time diffusion model
%
%
%e1 #&#
\begin{equation}
\label{modeldiffusion} d \bX(t) = \bmu_t \,dt +
\bsigma_t^T \,d \bB_t,\qquad t\in[0,1],
\end{equation}
where stochastic processes $\bX(t)$, $\bB_t$, $\bmu_t$ and $\bsigma_t$
are defined on a filtered probability space
$(\Omega, \cF, \{\cF_t, t \in[0,1]\}, P)$ with filtration $\cF_t$
satisfying the usual conditions,
$\bB_t$ is a $p$-dimensional standard Brownian motion with respect to
$\cF_t$,
$\bmu_t$ is a $p$-dimensional drift vector, $\bsigma_t$ is a $p$ by $p$
matrix, and
$\bmu_t$ and $\bsigma_t$ are assumed to be predictable processes with
respect to $\cF_t$.

We assume that the continuous-time process $\bX(t)$ is observed with
measurement errors only at equally spaced discrete
time points; that is, the observed discrete
data %$\bY(t_\ell)=(Y_1(t_\ell, \ldots, Y_p(t_\ell))^T$
$Y_i(t_\ell)$ obey
%
%
%e2 #&#
\begin{equation}
\label{modelnoisy} Y_i(t_{\ell}) =
X_i(t_{\ell})+\varepsilon_i(t_{\ell}),\qquad
i=1,\ldots, p, t_\ell= \ell/n, \ell=1,\ldots,n,
\end{equation}
where $\varepsilon_i(t_{\ell})$ are noises with mean zero.\vadjust{\goodbreak}
%and variance $\eta_i$, and $\varepsilon_i(\cdot)$ and $X_i(\cdot)$ are
%independent.

Let $\bolds{\gamma}(t) = \bsigma_t^T \bsigma_t$ be the
volatility matrix of
$\bX(t)$.
We are interested in estimating the following integrated volatility matrix
of $\bX(t)$,
\[
\bGamma= (\Gamma_{ij})_{1\leq i,j \leq p} = \int_0^1
\bolds{\gamma}(t) \,dt = \int_0^1
\bsigma_t^T \bsigma_t \,dt
\]
based on noisy discrete data $Y_i(t_{\ell})$, $i=1,\ldots, p$,
$\ell=1,\ldots, n$.

%s2.2 #&#
\subsection{Estimator}
\label{secestimator}
Let $K$ be an integer and $\lfloor n/K \rfloor$ be the largest integer
$\leq n/K$.
We divide $n$ time points $t_1, \ldots, t_n$ into $K$ nonoverlap
groups $\btau^k = \{ t_\ell, \ell= k, K+k, 2 K+k, \ldots\}$,
$k=1,\ldots,K$.
Denote by $|\btau^k|$ the number of time points in $\btau^k$.
Obviously, the
value of $|\btau^k|$ is either $\lfloor n/K \rfloor$ or $\lfloor n/K
\rfloor+ 1$.
For $k=1,\ldots,K$,
we write the $r$th time point in $\btau^k$ as $\tau_r^k = t_{(r-1) K+k}$,
$r=1,\ldots,|\btau^k|$.
With each $\btau^k$, we define %co-volatility between components $i$
%and $j$
the volatility matrix estimator
%
%
%e3 #&#
\begin{eqnarray}
\label{RVM} \tilde\Gamma_{ij} \bigl(\btau^k \bigr) &=& \sum
_{r=2}^{|\btau^k|} \bigl[Y_i \bigl(
\tau_{r}^k \bigr)- Y_i \bigl(
\tau_{r-1}^k \bigr) \bigr] \bigl[Y_j \bigl(
\tau_{r}^k \bigr)-Y_j \bigl(
\tau_{r-1}^k \bigr) \bigr],
\nonumber
\\[-8pt]
\\[-8pt]
\nonumber
\tilde{ \bGamma} \bigl(
\btau^k \bigr)&=& \bigl(\tilde\Gamma_{ij} \bigl(
\btau^k \bigr) \bigr)_{1\leq i,j \leq p}.
\end{eqnarray}
Here in (\ref{RVM}), to account for noises in data $Y_i(t_\ell)$, we use
$\btau^k$ to subsample the
data and define $\tilde{\bGamma}(\btau^k)$. To reduce the noise effect
we average $K$ volatility matrix estimators $\tilde{\bGamma}(\btau
^k)$ to
define one-scale volatility matrix estimator
%
%
%e4 #&#
\begin{equation}
\label{RVM1} \tilde\Gamma_{ij}^K = \frac1K \sum
_{k=1}^K \tilde\Gamma_{ij} \bigl(
\btau^k \bigr),\qquad \tilde{\bGamma}{}^K = \bigl(\tilde
\Gamma_{ij}^K \bigr) = \frac1K \sum
_{k=1}^K \tilde{\bGamma} \bigl(\btau^k
\bigr).
\end{equation}
Let $N= [c n^{1/2}]$ for some positive constant $c$,
%Let $N= [n^{1/2}]$
and $K_m = m + N$, $m=1,\ldots, N$. We use each $K_m$ to define a one-scale
volatility matrix estimator $\tilde{\bGamma}{}^{K_m}$
and then combine them together to form a multi-scale
volatility matrix estimator
%
%
%e5 #&#
\begin{equation}
\label{MSRVM} \tilde{\bGamma} = \sum_{m=1}^N
a_m \tilde{\bGamma}{}^{K_m} + \zeta\bigl(\tilde{
\bGamma}{}^{K_1} - \tilde{\bGamma}{}^{K_N} \bigr),
\end{equation}
where
%
%
%e6 #&#
\begin{equation}
\label{coeff} \zeta= \frac{K_1 K_N}{n(N-1)},\qquad a_m = \frac{12
K_m(m-N/2-1/2)}{N(N^2-1)},
\end{equation}
which satisfy
%from simple algebra that
%
\[
\sum_{m=1}^N a_m = 1, \qquad\sum
_{m=1}^N \frac{a_m}{K_m} = 0,\qquad \sum
_{m=1}^N |a_m| = 9/2 + o(1).
% = 3/2 + 3N/N + o(1) = 9/2 + o(1)
\]
The one-scale matrix estimator in (\ref{RVM1}) was studied in \citet{WanZou10}, and the multi-scale scheme
(\ref{MSRVM})--(\ref{coeff}) in the univariate case was investigated in
\citet{Zha06}.

We threshold $\tilde{\bGamma}$ to obtain our final volatility matrix
estimator
%
%
%e7 #&#
\begin{equation}
\label{ThresholdMSRVM} \hat\bGamma= \bigl( \tilde
\Gamma_{ij} 1 \bigl(|\tilde\Gamma_{ij}| \geq\varpi\bigr) \bigr),
\end{equation}
where $\varpi$ is a threshold value to be specified in Theorem~\ref
{Tthreshold}.

In the estimation construction we use only time scales corresponding to
$K_m$ of order $\sqrt{n}$ to form increments and averages. In Section~\ref{secasymp}
we will
demonstrate that the data at these scales contain essential information
for estimating $\bGamma$ and show that $\hat{\bGamma}$
%in (\ref{ThresholdMSRVM})
is asymptotically an optimal estimator of $\bGamma$.

%s3 #&#
\section{Asymptotic theory% for the constructed volatility matrix
%estimators}
for volatility matrix estimators}
\label{secasymp}
First we fix notation for our asymptotic analysis.
%For two sequences $u_{n,p}$
%and $v_{n,p}$ we write $u_{n,p} \asymp v_{n,p}$ if there exist positive
%constants $c_1$ and $c_2$ free of $n$ and $p$ such that
%$c_1 \leq u_{n,p}/v_{n,p} \leq c_2$.
Let ${\mathbf x}=(x_1, \ldots, x_p)^T$ be a $p$-dimensional vector
and $\bA=(A_{ij})$ be a $p$ by $p$ matrix, and
define their $\ell_d$ norms
\[
\| {\mathbf x}\|_d = \Biggl( \sum_{i=1}^p
|x_i|^d \Biggr)^{1/d},\qquad \| \bA\|_d
= \sup\bigl\{ \|\bA{\mathbf x}\|_d, \|{\mathbf x}\|_d=1 \bigr\},\qquad 1 \leq d \leq
\infty.
\]
For the case of matrix, the $\ell_2$ norm is called the matrix spectral
norm. $\| \bA\|_2$ is equal to the square root of the largest
eigenvalue of $\bA\bA^T$,
%
%
%e8 #&#
\begin{equation}
\label{ell-1infty-norm} \| \bA\|_1 = \max_{1 \leq j \leq p}
\sum_{i=1}^p |A_{ij}|, \qquad\| \bA
\|_\infty= \max_{1 \leq i \leq p} \sum
_{j=1}^p |A_{ij}|
\end{equation}
and
%
%
%e9 #&#
\begin{equation}
\label{ell-12infty-norm} \|\bA\|_2^2 \leq\| \bA
\|_1 \| \bA\|_\infty.
\end{equation}
For symmetric $\bA$, (\ref{ell-1infty-norm})--(\ref{ell-12infty-norm})
imply that
$\| \bA\|_2 \leq\| \bA\|_1=\| \bA\|_\infty$,
and $\|\bA\|_2$ is equal to the largest absolute eigenvalue of $\bA$.

Second we state some technical conditions for the asymptotic analysis.
\begin{longlist}[A3.]
\item[A1.] Assume $n^{\beta/2} \leq p \leq\exp( \beta_0 \sqrt
{n})$ for some
constants $\beta>1$ and $\beta_0>0$, and that $\varepsilon_i(t_{\ell
})$ and
$\bX(t)$ in models (\ref{modeldiffusion})--(\ref{modelnoisy})
are independent. %$\varepsilon_i(t_\ell)$, $i=1,\ldots,p$, $\ell=1,
Suppose that $(\varepsilon_1(t_\ell), \ldots, \varepsilon_p(t_\ell))$,
$\ell=1,\ldots, n$, is a strictly stationary
$M$-dependent multivariate time series with mean zero and
$\operatorname{Var}[\varepsilon_i(t_\ell)] = \eta_i \leq\kappa^2$, where $M$ is a
fixed integer, and $\kappa$ is a finite positive constant.
%and its auto-covariance and cross-covariance satisfying
% \sum_{r=0}^{n-1} |Cov[\varepsilon_i(t_\ell), \varepsilon_i(t_{\ell+
%r})]| \leq C_1,
% \sum_{r=0}^{n-1} |Cov[\varepsilon_i^2(t_\ell), \varepsilon_i^2(t_{
% |Cov[\varepsilon_i(t_\ell), \varepsilon_j(t_{\ell+ r})]| \leq
%C_3/(r^2+1), i, j = 1, \ldots, p,
%where $C_1, C_2, C_3$ are positive constants.
Assume further that
%$\varepsilon_i(t_{\ell})$ are linear processes with representation
%$\varepsilon_i(t_{\ell}) = \sum_{j=-\infty}^\infty\phi_{ij}
%$\omega_i(t_\ell)$ are zero mean and unit variance i.i.d. random
%variables obeying subGaussianity in the sense that there is %$\tau>0$
%such that for any $s>0$ and infinite dimensional vector $v=(v_1, v_2,
% P( |v^T \{\omega_i(t_\ell)\}_{\ell= -\infty}^\infty|>s) \leq C e
%^{-s^2/(2 \tau)}, i=1, \ldots, p.
$\varepsilon_i(t_\ell)$ %$(\varepsilon_1(t_1), \ldots,
are subGaussian in the sense that there exist constants $\tau_0 >0$
and $c_0>0$ such that for all $x>0$ and ${\mathbf u}=(u_1, \ldots,
u_n)^T$ with
$\|{\mathbf u}\|_2=1$,
%
%
%e10 #&#
\begin{equation}
\label{subGaussian-2} P \bigl(\bigl| \bigl(\varepsilon_i(t_1),
\ldots, \varepsilon_i(t_n) \bigr) {\mathbf u}\bigr|> x \bigr) \leq
c_0 e ^{-x^2/(2 \tau_0)},\qquad i=1,\ldots, p.
\end{equation}

\item[A2.] Assume that %each component of drift $\bmu(t)$ has bounded
%variation, and
there exist positive constants $c_1$ and $c_2$ such that
\[
\max_{1 \leq i \leq p} \max_{0 \leq t \leq1} \bigl|
\mu_i(t)\bigr| \leq c_1,\qquad \max_{1 \leq i \leq p} \max
_{0 \leq t \leq1} \gamma_{ii}(t) \leq c_2.
\]
%
%where $c_1$ and $c_2$ are positive constants.
Further we assume with probability one for $t \in[0,1]$,
%
%e11 #&#
\begin{eqnarray}
\gamma_{ii}(t) > 0, i=1, \ldots, p,\qquad \gamma_{ii}(t) +
\gamma_{jj}(t) \pm2 \gamma_{ij}(t)>0,\nonumber\\
\eqntext{ i \neq j, i,j=1,
\ldots, p.}
\end{eqnarray}

% \max_{1 \leq i \leq p} \max_{0 \leq t \leq1} |\mu_i(t)| \leq c_3,
% \max_{1 \leq i \leq p} \max_{0 \leq t \leq1} E
% \left[ \exp( \beta\mu^2_{i}(t)) \right] < \infty,
% \max_{1 \leq i \leq p} \max_{0 \leq t \leq1}
% E \left[ \exp(\beta\gamma_{ii}(t)) \right] < \infty,
% \max_{1 \leq i \leq p} E \left[ \exp( \beta
% \varepsilon^2_i(t_{i\ell}))
% \right] < \infty.

\item[A3.] Assume that $\bGamma$ is sparse in the sense that
%impose sparsity on $\bGamma$,
%
%
%e12 #&#
\begin{equation}
\label{Csparsity} \sum_{j=1}^p
|\Gamma_{ij}|^q \leq\Psi\pi_n(p), \qquad i=1,
\ldots, p,
\end{equation}
where $\Psi$ is a positive random variable with finite second moment,
$0 \leq q <1$, and $\pi_n(p)$ is a deterministic function with slow growth
%growing very slowly
in $p$ such as $\log p$.
%For $q=0$, matrix $\gamma(t)$ has at most $\Psi\pi(p)$ nonzero
%elements on each column/row.
\end{longlist}

%Let $\rho_i(r)= Cov(\varepsilon_i(t_\ell), \varepsilon_i(t_{\ell+
%r})$, $\rho_i^{Sq}(r)= Cov(\varepsilon_i^2(t_\ell), \varepsilon_i^2(t_{

Condition A1 allows noises to have cross sectional correlations as well
as cross temporal correlations.
%Specifically, noises $(\varepsilon_1(t_{\ell}), \ldots,
%are serially independent if lags are more than $M$ lags apart,
In particular we may have any contemporaneous correlations between
$\varepsilon_i(t_{\ell})$ and $\varepsilon_j(t_{\ell})$
as well as lagged serial auto-correlations for individual noise
$\varepsilon_i(\cdot)$ and lagged serial cross-correlations between
$\varepsilon_i(\cdot)$ and $\varepsilon_j(\cdot)$ with lags up to $M$.
%Subgaussianity (\ref{subGaussian-2}) and the strict stationarity imply
%that $(\varepsilon_1(t_\ell), \ldots, \varepsilon_p(t_\ell))$, $
As in covariance matrix estimation, the subGaussianity (\ref
{subGaussian-2}) is essentially required to obtain an optimal
convergence rate depending on $p$ through $\sqrt{\log p}$. It is
obvious that independent normal noises satisfy these assumptions.
%the normality and independence assumptions on noise $\varepsilon_i(t_{
%for easy management of technical proofs. These assumptions can be
%relaxed to nonnormal noises with
%cross sectional correlation as well as cross temporal correlation(see
%Supplemental Materials in Appendix).
%%contemporaneous correlations as well as serial correlations
%%the first part is a typical assumption in the literature of
%measurement error models;
%%and volatility analysis for high-frequency financial data.
The constraint $p \geq n^{\beta/2}$ is needed to obtain a
high-dimensional minimax lower bound;
otherwise the problem will be similar to usual asymptotics with large
$n$ but fixed $p$;
$p \leq\exp(\beta_0\sqrt{n})$ is to ensure the existence of a
consistent estimator of $\bGamma$.
%otherwise %$\bGamma$ cannot be consistently estimated.
%the minimax lower bound will stay away from zero as $n$ and $p$ go to
%infinity.
Condition A2 is to impose proper assumptions on the drift and
volatility of
the It\^o process so that we can obtain subGaussian tails for the
quadratic forms of $X_i(t_\ell)$, which together with the
subGaussianity (\ref{subGaussian-2}) are used to derive subGaussian
tails for
the elements of the volatility matrix estimator $\tilde{\bGamma}$.
Condition A3 is a common sparsity
assumption required for consistently estimating large matrices [\citet{BicLev08N2}, \citet{CaiZho}, and \citet{JohLu09}].

The following two theorems establish asymptotic theory for the
estimators $\tilde{\bGamma}$ and $\hat\bGamma$ defined by (\ref{MSRVM})
and (\ref{ThresholdMSRVM}), respectively.

%
%th1 #&#
\begin{theorem}\label{UnivariateConv}
Under models (\ref{modeldiffusion})--(\ref{modelnoisy}) and
conditions \textup{A1--A2},
the estimator $\tilde{\bGamma}$ in (\ref{MSRVM}) satisfies that for
$1\leq i,j \leq p$
and positive $x$ in a neighbor of~$0$,
%
%
%e13 #&#
\begin{equation}
\label{convelement} \mathbb{P} \bigl(\vert\tilde
\Gamma_{ij}-\Gamma_{ij}\vert\geq x \bigr) \leq
\varsigma_1 \exp\bigl\{ \log n % \log(n/x)
- \sqrt{n}
x^2/\varsigma_0 \bigr\}, % \mbox{ for positive $x$ in a neighbor of
%$0$},
\end{equation}
where $\varsigma_0$ and $\varsigma_1$ are positive constants free of
$n$ and $p$.
\end{theorem}

%
%re1 #&#
\begin{remark}\label{rem1}
Theorem~\ref{UnivariateConv} establishes subGaussian
tails for
the elements of the matrix estimator $\tilde{\bGamma}$. It is known
that, when univariate
or bivariate continuous It\^o processes are observed with measurement
errors at $n$ discrete time
points, the optimal convergence rates for estimating a univariate
integrated volatility or
a bivariate integrated co-volatility are $n^{-1/4}$ [Gloter and Jacod
(\citeyear{GloJac01N1,GloJac01N2}), \citet{Rei11}, and \citet{Xiu10}].
The $\sqrt{n} x^2$ factor in the exponent of the tail probability
bound on the right-hand side of (\ref{convelement}) %given by Theorem~\ref{UnivariateConv}
indicates a $n^{-1/4}$ convergence rate for
$\tilde\Gamma_{ij}-\Gamma_{ij}$, which matches the optimal convergence
rate for the
univariate integrated volatility estimation. This is in contrast to
%%comparison with
sub-optimal convergence rate results %on limiting distributions and
%tail probabilities
in the literature where a $n^{-1/6}$ convergence\vadjust{\goodbreak} rate was obtained;
see, for example,
\citet{FanLiYu12}, \citet{WanZou10}, \citet{ZhaMykAtS05}, and \citet{ZheLi11}.
\end{remark}
%
%th2 #&#
\begin{theorem} \label{Tthreshold}
For the threshold estimator $\hat\bGamma$ in (\ref{ThresholdMSRVM})
we choose
threshold $\varpi= \hbar n^{-1/4}\sqrt{\log(np)}$ with any fixed constant
$\hbar\geq5 \sqrt{\varsigma_0}$, where $\varsigma_0$ is the
constant in
the exponent of the tail probability bound on the right-hand side of
(\ref{convelement}).
Denote by $\mathcal{P}_{q}(\pi_{n}(p))$ the set of distributions of
$Y_i(t_\ell)$, $i=1,\ldots, p$, $\ell=1, \ldots, n$, from models
(\ref{modeldiffusion})--(\ref{modelnoisy}) satisfying conditions \textup{A1--A3}.
Then as $n, p \rightarrow\infty$,
%
%
%e14 #&#
\begin{eqnarray}
\label{finalrate} \sup_{\mathcal{P}_{q}(\pi_{n}(p))} \mathbb{E}
\llVert
\hat\bGamma- \bGamma\rrVert_2^2 &\leq&\sup
_{\mathcal{P}_{q}(\pi_{n}(p))} \mathbb{E} \llVert\hat\bGamma-
\bGamma\rrVert
_1^2
\nonumber
\\[-8pt]
\\[-8pt]
\nonumber
& \leq& C^* \bigl[ \pi_n(p) \bigl(
n^{-1/4} \sqrt{\log p} \bigr) ^{1-q} \bigr]^2,
\end{eqnarray}
where $C^*$ is a constant free of $n$ and $p$.
\end{theorem}

%
%re2 #&#
\begin{remark}\label{rem2} For sparse covariance
matrix estimation, \citet{CaiZho} has shown that the threshold
estimator in
\citet{BicLev08N2} is rate-optimal, and the optimal convergence
rate depends
on $n$ and $p$ through $n^{-1/2} \times\sqrt{\log p}$.
The convergence rate obtained in Theorem~\ref{Tthreshold}
depends on the sample size $n$ and the matrix size $p$ through
$n^{-1/4}  \sqrt{\log p}$.
Note that $n^{-1/4}$ is the optimal convergence rate for estimating a
univariate integrated volatility or a bivariate integrated
co-volatility based on noisy data.
%(Gloter and Jacod (\citeyear{GloJac01N1},b) and \citet{MunSch10})
Since our estimation problem is a generalization of covariance matrix
estimation for i.i.d. data
to volatility matrix estimation for an It\^o process with measurement
errors on one hand and a high-dimensional extension of
univariate volatility estimation on the other hand, it is interesting
to see that
the convergence rate in Theorem~\ref{Tthreshold} is a natural blend of
convergence rates in the two cases.
Also as Theorem~\ref{Tthreshold} implies that the maximum of the
eigenvalue differences between $\hat\Gamma$
and $\Gamma$ is bounded
by $\sqrt{C^*} \pi_n(p) ( n^{-1/4} \sqrt{\log p} ) ^{1-q}$.
Thus if the eigenvalues of $\bGamma$ all exceed
$\sqrt{C^*} \pi_n(p) ( n^{-1/4} \sqrt{\log p} ) ^{1-q}$,
asymptotically the eigenvalues of $\hat\Gamma$ are positive, and
$\hat
\Gamma$ is a positive definite matrix.
In particular, if $\pi_n(p) ( n^{-1/4} \sqrt{\log p} )
^{1-q}$ goes to zero as $n$ and $p$ go to infinity,
and $\bGamma$ is positive definite and well conditioned, then $\hat
\Gamma$ is asymptotically positive definite and
well conditioned.
In Section~\ref{optimalsec} we will establish the minimax lower bound
for estimating $\bGamma$ and show that the convergence rate
in Theorem~\ref{Tthreshold} is optimal.
\end{remark}
%In comparison of Theorem~\ref{Tthreshold}
%with the optimal convergence rate for covariance matrix estimation,
%the convergence rate in Theorem~\ref{Tthreshold} has a similar form
%but
%depends on $n$ in terms of $n^{-1/4}$ instead of $n^{-1/2}$ for
%covariance
%matrix estimation. The slower convergence rate here is intrinsically
%due to
%the noise contamination in the observed data under our set-up. The
%$n^{-1/4}$
%convergence rate conforms with the optimal convergence rate for
%estimating
%univariate integrated volatility, and will be heuristically explained
%in
%Remark~5 after the minimax lower bound result in Section~3.

%For threshold $\varpi= M n^{-1/4}\sqrt{\log(np)}$ we may take
%any $M \geq5 \sqrt{C_0}$, where $C_0$ is the same constant as in
%the exponent of tail probability bound in (\ref{convelement}).
%Also the assumption that $p \leq\exp( b_2 \sqrt{n})$
%in Theorem~3 is necessary for the consistency of $\hat\bGamma$.
%Otherwise,
%its risk $ \pi^2_n(p) \left( n^{-1/4} \sqrt{\log p} \right) ^{2-2q}$
%will be bounded below by zero.

%s4 #&#
\section{Optimal convergence rate}
\label{optimalsec}

This section establishes the minimax lower bound for estimating
$\bGamma$
under models (\ref{modeldiffusion})--(\ref{modelnoisy}) and shows
that asymptotically
$\hat{\bGamma}$ achieves the lower bound and thus is optimal.
We state the minimax lower bound for estimating $\bGamma$
with $\mathcal{P}_{q}(\pi_{n}(p))$ under the matrix spectral norm as follows.

%
%th3 #&#
\begin{theorem} \label{lowerbound}
For models (\ref{modeldiffusion})--(\ref{modelnoisy}) satisfying conditions
\textup{A1--A3}, if %$\pi_{n}(p)\leq M n^{(1-q)/4}/ (\log p)^{(3-q)/2}$
for some constant $\aleph>0$,
%
%
%e15 #&#
\begin{equation}
\pi_{n} ( p ) \leq\aleph n^{ ( 1-q ) /4}/ ( \log p )
^{ ( 3-q ) /2}, \label{condc}
\end{equation}
the minimax risk for estimating $\bGamma$ with
$\mathcal{P}_{q}(\pi_{n}(p))$ satisfies that as $n, p \rightarrow
\infty$,
%
%
%e16 #&#
\begin{equation}
\label{lowerbd1} \inf_{\check{\bGamma}}\sup_{\mathcal{P}_{q}(\pi_{n}(p))}
\mathbb{E}\llVert\check{\bGamma}-\bGamma\rrVert_2^{2}\geq
C_* \bigl[ \pi_{n}(p) \bigl( n^{-1/4} \sqrt{\log p}
\bigr)^{1-q} \bigr]^2,
\end{equation}
where $C_*$ is a positive constant free of $n$ and $p$, and the infimum
is taken over all estimators
$\check{\bGamma}$ based on the data $Y_i(t_\ell)$, $i=1,\ldots, p$,
$\ell=1,\ldots, n$, from models~(\ref{modeldiffusion})--(\ref
{modelnoisy}).
\end{theorem}

%re3 #&#
\begin{remark}\label{rem3} Note that the lower bound convergence rate in
Theorem~\ref{lowerbound}
matches the convergence rate of the estimator $\hat\bGamma$ obtained in
Theorem~\ref{Tthreshold}.
%Since for a symmetric matrix $\bA$, the Riesz-Thorin interpolation
%theorem
%implies that for all $d \geq1$, $\| \bA\|_d \leq\| \bA\|_1 = \|
%%\| \cdot\|_1 = \| \cdot\|_\infty$.
Combining Theorems~\ref{Tthreshold} and~\ref{lowerbound} together we conclude
that the optimal convergence rate is $\pi_{n}(p) ( n^{-1/4}
\sqrt{\log p} )^{1-q}$, and the estimator $\hat\bGamma$ in
(\ref{ThresholdMSRVM}) achieves the optimal convergence rate. Moreover,
such optimal estimation results hold for any matrix $\ell_d$ norm with
$1 \leq d \leq\infty$. Indeed, it can be shown that
under the conditions of Theorems~\ref{Tthreshold} and \ref
{lowerbound}, we have that as $n$ and $p$
go to infinity,
%
%
%e17 #&#
\begin{eqnarray}
\label{lowerbd-d-norm} &&\frac{C_*}{4} \bigl[ \pi_{n}(p) \bigl(
n^{-1/4} \sqrt{\log p} \bigr)^{1-q} \bigr]^2 \nonumber\\
&&\qquad\leq
\inf_{\check{\bGamma}}\sup_{\mathcal{P}_{q}(\pi_{n}(p))} \mathbb
{E}\llVert
\check{\bGamma}-\bGamma\rrVert_d^{2} \leq\sup
_{\mathcal{P}_{q}(\pi_{n}(p))} \mathbb{E} \llVert\hat\bGamma-
\bGamma\rrVert
_d^2
\\
& &\qquad\leq C^* \bigl[ \pi_n(p) \bigl( n^{-1/4} \sqrt{\log p}
\bigr) ^{1-q} \bigr]^2,\nonumber
\end{eqnarray}
where $C^*$ and $C_*$ are constants in Theorems~\ref{Tthreshold} and
\ref{lowerbound}, respectively, %free of $n$ and $p$,
$\hat\bGamma$ is the threshold estimator given by
(\ref{ThresholdMSRVM}) with the threshold value specified in Theorem~\ref{Tthreshold} and
the infimum is taken over all estimators $\check{\bGamma}$ based on the
data $Y_i(t_\ell)$,
$i=1,\ldots, p$, $\ell=1,\ldots, n$, from models (\ref
{modeldiffusion})--(\ref{modelnoisy}).
\end{remark}

%re4 #&#
\begin{remark}\label{rem4} Condition (\ref{condc}) is a technical condition
that we need to establish the minimax lower
bound. It is compatible with conditions A1 and~A3 regarding the
constraint on $n$ and $p$ as well as
the slow growth of $\pi_n(p)$ in the sparsity condition (\ref{Csparsity}).

%The constraint $\pi_{n}(p)\leq M \left( n^{1/4}/\sqrt{\log p}
%is necessary for the existence of a consistent estimator. In fact, for
%any
%$M_{1}>0$, the right-hand side of Equation (\ref{lowerbd1}) satisfies
%C_* \pi^2_{n}(p) \left( n^{-1/4} \sqrt{\log p}\right)^{2-2q} \geq
%M_{1} \text{if
%and only if }\pi_{n}(p)\geq\sqrt{M_1/C_*} %\left( M_{1}/C_*
%If the assumption $\pi_{n}(p)\leq M \left( n^{1/4}/\sqrt{\log p}
%bounded away from zero and thus no consistent estimator exists.

Models (\ref{modeldiffusion})--(\ref{modelnoisy}) are complicated
nonparametric models, and the observations from the models are
dependent and have
subGaussian measurement errors. To derive the minimax lower bound for
models (\ref{modeldiffusion})--(\ref{modelnoisy}),
we find a special subclass of the models to attain the minimax lower
bound of the models.
Such an approach is often referred to as the method of hardest
subproblem. Since generally a minimax problem
has lower bound no larger than any of its subproblems, the mentioned
special subclass corresponds to the hardest subproblem
and is referred to as the least favorable submodel. We will show in
Sections~\ref{transformationsec}
and~\ref{lowerbdsec} that the least favorable submodel for models
(\ref
{modeldiffusion})--(\ref{modelnoisy})
can be taken as i.i.d. Gaussian measurement errors $\varepsilon
_i(t_\ell
)$ and process $\bX(t)$ with zero drift and constant volatilities.
To establish the minimax lower bound for the least favorable submodel,
luckily we are able to find a
nice trick in Section~\ref{transformationsec} that transforms the
minimax lower bound problem for the least favorable
submodel into a new covariance matrix estimation problem with
independent but nonidentically distributed observations.
%We usually rely on Le Cam's method or Assouad's lemma to establish a
%minimax lower bound.
\citet{CaiZho} have developed an approach combining both Le Cam's
method and Assouad's lemma, which are two popular methods to establish minimax
lower bounds, to derive the minimax lower bound for estimating a
large sparse covariance matrix based on i.i.d. observations.
We adopt the approach in \citet{CaiZho} to derive the minimax lower
bound for the new covariance matrix estimation problem with independent but
nonidentically distributed observations, which is stated in Theorem~\ref{Operlowerbdthm} of
Section~\ref{lowerbdsec}. The derived minimax lower bound in Theorem~\ref{Operlowerbdthm} corresponds
to the least favorable submodel and thus is the minimax lower bound for
models (\ref{modeldiffusion})--(\ref{modelnoisy}).
Therefore, we prove Theorem~\ref{lowerbound}.
\end{remark}

%s4.1 #&#
\subsection{Model transformation}
\label{transformationsec}
We take a subclass of models (\ref{modeldiffusion})--(\ref
{modelnoisy}) as follows.
For the It\^o processes $\bX(t)$ we let $\bmu_t=0$ and $\bsigma_t$ be
a constant matrix~$\bsigma$; for the noises we let $\varepsilon
_i(t_\ell
)$, $i=1,\ldots, p$, $\ell=1,\ldots, n$, be
i.i.d. random variables with $N(0, \kappa^2)$ distribution, where
$\kappa>0$ is specified in condition A1.
Then $\bGamma= (\Gamma_{ij})=\bsigma^T \bsigma$,
and the sparsity condition (\ref{Csparsity}) becomes
%
%
%e18 #&#
\begin{equation}
\sum_{j=1}^p |\Gamma_{ij}|^q
\leq c_3 \pi_n(p), \label{sparseparaspace}
\end{equation}
where $c_3=E(\Psi)$ and $\Psi$ is given by (\ref{Csparsity}).
%Denote by $\mathcal{G}_{q}(\pi_{n}(p))$ the class of constant matrices
%$\bGamma$ satisfying (\ref{sparseparaspace}) and
%$\| \bGamma\|_2 \leq\tau$, where $\tau$ is a positive constant.

Let $\bY_l = (Y_1(t_l), \ldots, Y_{p}(t_l))^T$, and
$\beps_l = (\varepsilon_1(t_l), \ldots, \varepsilon_p(t_l))^T$.
Then models (\ref{modeldiffusion})--(\ref{modelnoisy}) become
%
%
%e19 #&#
\begin{equation}
\label{low1} \bY_l = \bsigma\bB_{t_l} +
\beps_l, \qquad l = 1, \ldots, n, t_l = l/n
\end{equation}
and $\beps_l \sim N(0, \kappa^2 I_p)$. % where $\kappa>0$ is specified
%in Condition A1.
As $\bY_l$ are dependent, we take differences in (\ref{low1}) and obtain
%
%
%e20 #&#
\begin{equation}
\label{low2} \bY_l - \bY_{l-1} = \bsigma(
\bB_{t_l} - \bB_{t_{l-1}}) + \beps_l -
\beps_{l-1},\qquad l=1,\ldots, n,
\end{equation}
here $\bY_0=\beps_0 \sim N(0, \kappa^2 \bI_p)$. For matrix
\( (\beps_l - \beps_{l-1}, 1 \leq l \leq n) = (\varepsilon_i(t_l) -
\varepsilon_i(t_{l-1}),\allowbreak 1 \leq i \leq p, 1 \leq l \leq n) \),
its elements are independent at different rows but correlated
at the same rows. At the $i$th row, elements $\varepsilon_{i}(t_l) -
\varepsilon_{i}(t_{l-1})$, $l=1, \ldots, n$, have covariance matrix
$\kappa^2 \bUpsilon$, where $\bUpsilon$ is a $n \times n$ tridiagonal
matrix with $2$ along diagonal entries, $-1$ next to diagonal entries
%$=(a_{l\ell})$ is a n by n matrix with $a_{ll}=2$, $a_{l,l\pm1}=-1$
and $0$ elsewhere. $\bUpsilon$ is a Toeplitz matrix [\citet{Wil88}]
that can be diagonalized as follows:
%
%
%e21 #&#
\begin{equation}
\label{low21} \bUpsilon= \bQ\bPhi\bQ^T, \qquad\bPhi= \operatorname{diag}(
\varphi_1, \ldots, \varphi_n),
\end{equation}
where $\varphi_l$ are eigenvalues with expressions
%
%
%e22 #&#
\begin{equation}
\label{low22} \varphi_l= 4 \sin^2 \biggl[
\frac{\pi l}{2(n+1)} \biggr],\qquad l = 1, \ldots, n,
\end{equation}
and $\bQ$ is an orthogonal matrix formed by the eigenvectors of
$\bUpsilon$.
Using (\ref{low21}) we transform the $i$th row of the matrix $(\beps
_l -
\beps_{l-1}, 1 \leq l \leq n)$ by $\bQ$, and obtain
\[
\operatorname{Var} \bigl[ \bigl(\varepsilon_{i}(t_1) -
\varepsilon_{i}(t_{0}), \ldots, \varepsilon
_{i}(t_n) - \varepsilon_{i}(t_{n-1})
\bigr) \bQ\bigr] = \kappa^2 \bQ^T \bUpsilon\bQ=
\kappa^2 \bPhi.
\]
For $i=1,\ldots, p$, let
\begin{eqnarray*}
(e_{i 1},\ldots, e_{i n}) &= &\bigl(\sqrt{n} \bigl[
\varepsilon_{i}(t_1) - \varepsilon_{i}(t_{0})
\bigr], \ldots, \sqrt{n} \bigl[\varepsilon_{i}(t_n) -
\varepsilon_{i}(t_{n-1}) \bigr] \bigr) \bQ,
\\
(u_{i 1}, \ldots, u_{i n})& = &\bigl(\sqrt{n}
\bigl[Y_{i}(t_1) - Y_{i}(t_{0})
\bigr], \ldots, \sqrt{n} \bigl[Y_{i}(t_n) -
Y_{i}(t_{n-1}) \bigr] \bigr) \bQ,
\\
(v_{i 1}, \ldots, v_{i n}) &=& \bigl(\sqrt{n} \bigl[
B_{i}(t_1) - B_{i}(t_{0}) \bigr],
\ldots, \sqrt{n} \bigl[ B_{i}(t_n) -
B_{i}(t_{n-1}) \bigr] \bigr) \bQ.
\end{eqnarray*}
%
% \varepsilon_{i}(t_{l-1})], 1 \leq l \leq n) \bQ, \]
% 1 \leq l \leq n) \bQ, \]
% 1 \leq l \leq n) \bQ. \]
Then as $\bQ$ diagonalizes $\bUpsilon$, $e_{il}$ are independent, with
$e_{il} \sim
N(0, n \kappa^2 \varphi_l)$; %and $\be_i \sim N(0, n \kappa^2 \bPhi)$;
because $B_{i}(t_{l})- B_{i}(t_{l-1})$ are i.i.d. normal random
variables with mean zero
and variance $1/n$, and $\bQ$ is orthogonal, $v_{il}$ are i.i.d.
standard normal random variables.

Put (\ref{low2}) in a matrix form and
%[B_{i,t_l} - B_{i, t_{l-1}} ], 1 \leq i \leq p, 1 \leq l \leq n) \]
%l \leq n) \]
right multiply by $\sqrt{n} \bQ$ on both sides to
%n) \bQ\]
%l \leq n) \bQ\]
obtain
\[
(u_{il}) = \bsigma(v_{il}) + (e_{il}).
\]
Denote by $\bU_l$, $\bV_l$ and $\be_l$ the column vectors of the
matrices $(u_{il})$, $(v_{il})$ and $(e_{il})$, respectively. Then
the above matrix equation is equivalent to
%and write above matrix equation in terms of their columns,
%
%
%e23 #&#
\begin{equation}
\label{low3} \bU_l = \bsigma\bV_l +
\be_l,\qquad l =1, \ldots, n,
\end{equation}
where $\be_l \sim N(0, \kappa^2 n \varphi_l \bI_p)$ and $\bV_l
\sim
N(0, \bI_p)$.

From (\ref{low3}) we have that the data transformed random vectors
$\bU_{1}, \ldots, \bU_n$ are independent with
$\bU_{l}\sim N ( 0,\bGamma+ ( a_{l}-1 ) \bI_{p} )$,
where $a_{l}= 1 + \kappa^2 n \varphi_l$
%1+4\kappa^{2}n\sin^{2}\left( \frac{\pi l}{2\left( n+1\right) }
with \mbox{$0<\kappa<\infty$}.

%s4.2 #&#
\subsection{Lower bound}
\label{lowerbdsec}

We convert the minimax lower bound problem stated in Theorem~\ref{lowerbound}
%for estimating matrix $\bGamma$ under models
%(\ref{modeldiffusion})-(\ref{modelnoisy})
into a much simpler problem of estimating $\bGamma$ based on the
observations $\bU_1,\ldots, \bU_n$ from model (\ref{low3}), where
$\bGamma$ are
constant matrices satisfying
(\ref{sparseparaspace}) and $\| \bGamma\|_2 \leq\tau$ for some constant
$\tau>0$. We denote the new minimax estimation problem by
$\mathcal{Q}_{q}(\pi_{n}(p))$, and
%Denote by $\mathcal{G}_{q}(\pi_{n}(p))$ the class of constant matrices
%$\bGamma$ satisfying (\ref{sparseparaspace}) and
%$\| \bGamma\|_2 \leq\tau$, where $\tau$ is a positive constant,
%and $\mathcal{Q}_{q}(\pi_{n}(p))$ the problem of estimating
%$\bGamma\in\mathcal{G}_{q}(\pi_{n}(p))$ based on observations
%$\bU_1,\ldots, \bU_n$ from (\ref{low3}).
the theorem below derives its minimax lower bound.
%
%
%th4 #&#
\begin{theorem}
\label{Operlowerbdthm}
%Let $\bU_{l}$ be independent with $\bU_{l}\sim
%N\left( 0,\bGamma+\left( a_{l}-1\right) \bI_{p}\right)$, $l=1,\ldots,
%n$.
Assume $p \geq n^{\beta/2}$ for some $\beta>1$.
If $\pi_{n}(p)$ obeys (\ref{condc}), %$\pi_{n}(p)\leq M n^{(1-q)/4}/ (
the minimax risk for estimating matrix $\bGamma$
with $\mathcal{Q}_{q}(\pi_{n}(p))$
%over the parameter space $\mathcal{G}_{q}(\pi_{n}(p))$
satisfies that as $n, p \rightarrow\infty$,
%
%
%e24 #&#
\begin{equation}
\inf_{\check{\bGamma}}\sup_{\mathcal{Q}_{q}(\pi_{n}(p))} \mathbb
{E}\llVert
\check{\bGamma}-\bGamma\rrVert_2^{2}\geq C_* \bigl[
\pi_{n}(p) \bigl( n^{-1/4} \sqrt{\log p} \bigr)^{1-q}
\bigr]^2, \label{lowerbd}
\end{equation}
where $C_*$ is a positive constant free of $n$ and $p$, and the infimum
is taken over all estimators
$\check{\bGamma}$ based on the observations $\bU_1, \ldots, \bU_n$ from
model (\ref{low3}).
%where $\Vert\cdot\Vert$ denotes the matrix spectral norm$.
\end{theorem}

%re5 #&#
\begin{remark}\label{rem5}As we discussed in Remarks~\ref{rem1} and~\ref{rem2}
in Section~\ref{secasymp},
due to noise contamination, the optimal convergence rate depends on sample
size through $n^{-1/4}$, instead of $n^{-1/2}$ for covariance matrix
estimation. %Through the transformation in Section~\ref{transformationsec}
For the univariate case, discrete sine transform was used to
construct a realized volatility estimator [\citet{AtSMykZha05} and \citet{CurCor12}] and reveal
some intrinsic insight into how the $n^{-1/4}$ convergence rate is
obtained [\citet{MunSch10}].
%Actually Section~\ref{transformationsec}
The similar insight for the high-dimensional case can be seen from the
transformation
in Section~\ref{transformationsec}, which converts model
(\ref{low2}) with noisy data into model (\ref{low3}) where the
independent random
vector $\bU_l$ follows a multivariate normal distribution with mean
zero and
covariance matrix $\bGamma+ \kappa^2 n \varphi_l \bI_p$,
$l=1,\ldots, n$.
The transformation via orthogonal matrix $\bQ$, which diagonalizes
Toeplitz matrix $\bUpsilon$ and is equal to
$ ( \sin(\ell r \pi/(n+1)), 1 \leq\ell, r \leq n
) $ normalized by $\sqrt{2/(n+1)}$ [see \citet{Kho06}],
%Thus the transformation from model (\ref{low2}) to model (\ref{low3})
corresponds to a discrete sine transform, with (\ref{low3}) in frequency
domain and $\bU_l \sim N(0, \bGamma+ \kappa^2 n \varphi_l \bI_p)$
corresponding to the discrete sine transform of the data at
frequency $l \pi/(n+1)$.
%Note that (\ref{low22}) indicates that $n \varphi_l$
%behave like $l^2/n$. Since for $l$ with much higher order than $
%$\bGamma+ \kappa^2 n \varphi_l \bI_p$ are dominated by $\kappa^2
%n \varphi_l \bI_p$, the corresponding $\bU_l$ essentially behave like
%noise $\be_l$ and thus are not informative for estimating $\bGamma$.
%On the other hand, for $l$ up to the order of $\sqrt{n}$,
%$n \varphi_l$ are relatively small, and $\bGamma+ \kappa^2 n \varphi_l
%statistically
%similar to $\bV_l$. Hence, there are only $\sqrt{n}$ number of
%frequencies
%at which the transformed data $\bU_l$ are informative %contain
%essential information
%for estimating $\bGamma$, and we use %the $\sqrt{n}$ number of these $
%estimate covariance matrix $\bGamma$ and obtain $(
By comparing the order of $n \varphi_l$, we derive that only at those
frequencies with $l$ up to $\sqrt{n}$,
the transformed data $\bU_l$ are informative for estimating $\bGamma$,
and we use these $[\sqrt{n}]$ number of $\bU_l$ to
estimate %covariance matrix
$\bGamma$ and obtain $(\sqrt{n})^{-1/2}=n^{-1/4}$ convergence rate.
In fact, we have seen the phenomenon in Section~\ref{secmodel} where the $N$
scales used in the construction
of $\tilde\Gamma$ in (\ref{MSRVM}) correspond to $K_m$, with both $N$
and $K_m$ of order $\sqrt{n}$.
\end{remark}

%s5 #&#
\section{A simulation study}\label{sec5}
A simulation study was conducted to compare the finite sample
performances of
the MSRVM estimator in (\ref{MSRVM}) and the threshold MSRVM estimator
in (\ref{ThresholdMSRVM}) with those of the ARVM estimator and the
threshold ARVM estimator introduced in
\citet{WanZou10}.
%To simplify the simulation computation,
We generated $\bX(t)=(X_1(t), \ldots, X_p(t))^T$ at discrete time
points $t_{\ell}=\ell/n$, $\ell=1,\ldots,n$, from model~(\ref
{modeldiffusion}) with $\bmu_t=0$ by the Euler scheme, where
univariate standard Brownian motions were stimulated by the normalized
partial sums of independent standard normal random variables, $\bsigma
_{t_\ell}$ was taken to be a Cholesky decomposition of
\[
\bolds{\gamma}(t_\ell) = \bigl(\gamma_{ij}(t_\ell)
\bigr),\qquad \gamma_{ij}(t_\ell) = \sqrt{\gamma_{ii}(t_\ell)
\gamma_{jj}(t_\ell) } \varrho^{|i-j|},
\]
$\varrho$ was independently generated from a uniform distribution on
$[0.47, 0.53]$,
$(\gamma_{ii}(t_1), \ldots, \gamma_{ii}(t_n))$, $i = 1, \ldots, p$,
were independently drawn from %the stationary distribution of
a geometric Ornstein--Uhlenbeck process satisfying
$d \log\gamma_{ii}(t) = 6 [0.5-\log\gamma_{ii}(t)] \,dt + d W_i(t)$ and
$W_i(t)$ are independent one-dimensional standard Brownian motions that
are independent of $\bB_t$ in model (\ref{modeldiffusion}).
%whose stationary distribution has mean $1/2$ and variance $1/12$,
We computed $\bGamma$ by the average of $\bolds{\gamma}(t_1),
\ldots, \bolds{\gamma}(t_n)$.
We simulated noises $\varepsilon_i(t_\ell)$ independently from a normal
distribution with mean 0 and standard deviation $\theta\sqrt{\Gamma_{ii}}$,
$i=1, \ldots, p$, where $\theta$ is the relative noise level ranging
from $0$ to $0.7$. Finally data $Y_i(t_\ell)$ were obtained
by adding the simulated $\varepsilon_i(t_\ell)$ to the generated
$X_i(t_\ell)$ according to model (\ref{modelnoisy}).
%With the usage of the relative noise level, the selection of the mean
%and variance for $\gamma_{ii}$ will actually not affect the simulation
%results too much.
%the value which represents the decay speed of the volatility matrix,
%The chosen of the interval is inspired by the simulation in Wang and
%Zou (2010) when their $\kappa(0) = 0.5$.
%With the usage of the relative noise level, the selection of the mean
%and variance for $\gamma_{ii}$ will actually not affect the simulation
%results too much.
Using the simulated data $Y_i(t_\ell)$ we computed the MSRVM estimator
and the threshold MSRVM estimator as well as the ARVM
estimator and the threshold ARVM estimator.
In the simulation study we took $n=200$ and $p=100$. We repeated the
whole simulation procedure $200$ times.
For a given matrix estimator $\check{\bGamma}$, a relative matrix
spectral norm error $\Vert\check{\bGamma} - \bGamma\Vert_2/\Vert
\bGamma\Vert_2$ was
used to measure its performance. We evaluated the mean relative matrix
spectral norm error (MRE) by the
average of the relative matrix spectral norm errors over the $200$
repetitions. % for comparing estimators.
As in \citet{WanZou10} we selected tuning parameters like threshold
of the estimators by
minimizing the respective MREs.

%
%f1 #&#
\begin{figure}

\includegraphics{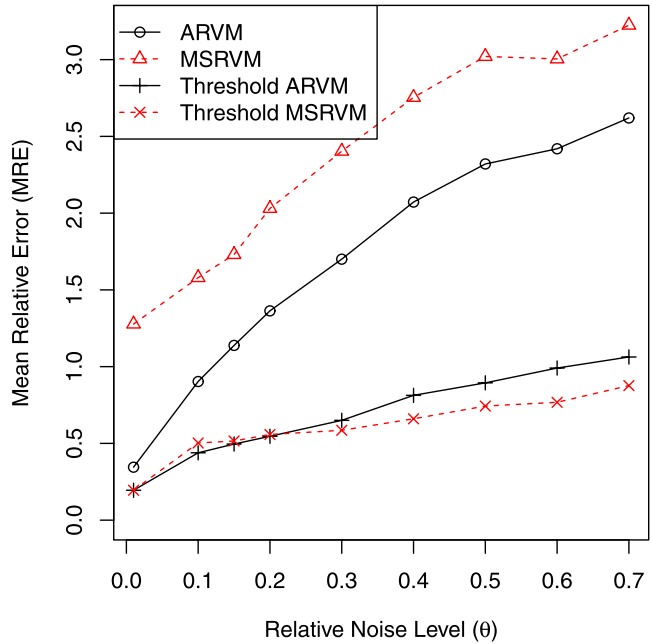}

\caption{The MRE plots of the four estimators for $n=200$ and $p=100$.}
\label{Flargep}
\end{figure}

Figure~\ref{Flargep} is the plots of MRE versus relative noise level
$\theta$ for the MSRVM, ARVM, threshold MSRVM and threshold ARVM
estimators. The basic findings are that while the MREs of the threshold
MSRVM and threshold ARVM estimators are comparable at
low relative noise levels, the threshold MSRVM estimator has smaller
MRE than the threshold ARVM estimator at high relative noise levels;
regardless of relative noise levels, the threshold MSRVM and threshold
ARVM estimators have significantly smaller MREs than the MSRVM and ARVM
estimators. The simulation results support the theoretical conclusions
that the threshold procedure is needed for constructing consistent
estimators of $\bGamma$, %large sparse volatility matrices,
and the threshold MSRVM estimator is asymptotically optimal, while the
threshold ARVM estimator is suboptimal.

We point out that it is important to have a data-driven choice of
tuning parameters for volatility matrix
estimator defined in (\ref{ThresholdMSRVM}). This is largely an open
issue. We briefly describe an approach
for developing a data-dependent selection of the tuning parameters as
follows. For data $\{Y_i(t_\ell), i=1, \ldots, p,\break \ell=1, \ldots,
n\}$ observed from models (\ref{modeldiffusion})--(\ref{modelnoisy}), we
may divide the whole data time interval
into $L$ subintervals $I_1, \ldots, I_L$, and partition data
$Y_i(t_\ell
)$ into $L$ subsamples
$\{Y_i(t_\ell), i=1, \ldots, p, t_\ell\in I_k\}$, $k=1, \ldots, L$,
over the $L$ corresponding time periods.
To estimate integrated volatility $\int_{I_k} \bolds{\gamma}(t)
\,dt/|I_k|$ over
the $k$th period, according to the procedure described
in Section~\ref{secestimator}, we use the $k$th subsample to construct volatility
matrix estimator, which is denoted by $\hat{\Gamma}_k(N, \varpi)$ to
emphasize its dependence on $N$ and $\varpi$, where $|I_k|$ denotes the
length of $I_k$, $\varpi$ is a threshold value
and $N$ is an integer that specifies scales used in the volatility
matrix estimator given by (\ref{ThresholdMSRVM}).
We predict one period ahead volatility matrix estimator $\hat{\Gamma
}_{k+1}(N,\varpi)$
by current period volatility matrix estimator $\hat{\Gamma
}_{k}(N,\varpi
)$ and compute the predication error.
We minimize the sum of the spectral norms of the predication errors to
select $N$ and $\varpi$.
For example, we often have high-frequency
financial data over many days, and it is natural to use data in each
day to estimate the integrated volatility matrix over
the corresponding day. We predict one day ahead daily volatility matrix
estimator by current daily volatility matrix estimator
and compute the predication error. The tuning parameters are then
selected by minimizing the sum of the spectral norms of the
prediction errors.

%s6 #&#
\section{Proofs}\label{sec6}
Denote by $C$'s generic constants whose values are free of $n$ and~$p$ and
may change from appearance to appearance. Let $u \vee v$ and $u \wedge
v$ be the
maximum and minimum of $u$ and $v$, respectively. For two sequences $u_{n,p}$
and $v_{n,p}$ we write $u_{n,p} \asymp v_{n,p}$ if there exist positive
constants $C_1$ and $C_2$ free of $n$ and $p$ such that
$C_1 \leq u_{n,p}/v_{n,p} \leq C_2$.
Without loss of generality we take $N=[n^{1/2}]$ in the construction of
$\tilde{\bGamma}$ given by (\ref{MSRVM}) in Section~\ref{secestimator}.

%s6.1 #&#
\subsection{\texorpdfstring{Proofs of Theorems \protect\ref{UnivariateConv} and \protect\ref{Tthreshold}}
{Proofs of Theorems 1 and 2}}
Let
\begin{eqnarray*}
\bY_r^{k_m} &=& \bigl(Y_1 \bigl(
\tau_r^{k_m} \bigr), \ldots, Y_p \bigl(
\tau_r^{k_m} \bigr) \bigr)^T,\\
\bX_r^{k_m} &=& \bigl(X_1 \bigl(
\tau_r^{k_m} \bigr), \ldots, X_p \bigl(
\tau_r^{k_m} \bigr) \bigr)^T,
\\
\bvarepsilon_r^{k_m} &=& \bigl(\varepsilon_1
\bigl(\tau_r^{k_m} \bigr), \ldots, \varepsilon
_p \bigl(\tau_r^{k_m} \bigr)
\bigr)^T,
\end{eqnarray*}
which are random vectors corresponding to the data, the It\^o process
and the noises at the time
point $\tau_r^{k_m}$, $r=1,\ldots, |\btau^{k_m}|$, $k_m=1,\ldots, K_m$,
and $m=1,\ldots, N$. Note that\vadjust{\goodbreak} we choose index $k_m$ to specify that
the analyses are associated with the study of $\bGamma^{K_m}$ here and
below. We decompose $\tilde{ \bGamma}{}^{K_m}$ defined in (\ref{RVM1})
as follows:
%
%
%e25 #&#
\begin{eqnarray}
\label{VG1} \tilde{\bGamma}{}^{K_m} & =& \frac1 {K_m} \sum
_{k_m=1}^{K_m} \sum
_{r=2}^{|\btau^{k_m}|} \bigl( \bY_r^{k_m}
- \bY_{r-1}^{k_m} \bigr) \bigl( \bY_r^{k_m}
- \bY_{r-1}^{k_m} \bigr)^T
\nonumber\\
& = &\frac1 {K_m} \sum_{k_m=1}^{K_m}
\sum_{r=2}^{|\btau^{k_m}|} \bigl(
\bX_r^{k_m} - \bX_{r-1}^{k_m} +
\bvarepsilon_r^{k_m} - \bvarepsilon_{r-1}^{k_m}
\bigr)\nonumber\\
&&\hspace*{59pt}{}\times \bigl( \bX_r^{k_m} - \bX_{r-1}^{k_m}
+ \bvarepsilon_r^{k_m} - \bvarepsilon_{r-1}^{k_m}
\bigr)^T
\nonumber\\
& =& \frac1 {K_m} \sum_{k_m=1}^{K_m}
\sum_{r=2}^{|\btau^{k_m}|} \bigl\{ \bigl(
\bX_r^{k_m} - \bX_{r-1}^{k_m} \bigr)
\bigl( \bX_r^{k_m} - \bX_{r-1}^{k_m}
\bigr) ^T
\nonumber
\\[-8pt]
\\[-8pt]
\nonumber
&&\hspace*{66pt}{} + \bigl(\bvarepsilon_r^{k_m} -
\bvarepsilon_{r-1}^{k_m} \bigr) \bigl(\bvarepsilon_r^{k_m}
- \bvarepsilon_{r-1}^{k_m} \bigr)^T
\\
&&\hspace*{66pt}{} + \bigl( \bX_r^{k_m} - \bX_{r-1}^{k_m}
\bigr) \bigl(\bvarepsilon_r^{k_m} - \bvarepsilon_{r-1}^{k_m}
\bigr)^T \nonumber\\
&&\hspace*{66pt}{}+ \bigl(\bvarepsilon_r^{k_m} -
\bvarepsilon_{r-1}^{k_m} \bigr) \bigl( \bX_r^{k_m}
- \bX_{r-1}^{k_m} \bigr) ^T \bigr\}
\nonumber\\
& \equiv&\bV^{K_m} + \bG^{K_m}(1) + \bG^{K_m}(2) +
\bG^{K_m}(3),\nonumber
\end{eqnarray}
and thus from (\ref{MSRVM}) we obtain the corresponding decomposition for
$\tilde{\bGamma}$,
%
%
%e26 #&#
\begin{eqnarray}
\label{VG} \tilde{\bGamma} %&=& \sum_{m=1}^N a_m \tilde\bGamma^{K_m}
%+ \zeta(
&=& \sum
_{m=1}^N a_m \bV^{K_m} +
\zeta\bigl(\bV^{K_1} - \bV^{K_N} \bigr)\nonumber \\
&&{}+ \sum
_{r=1}^3 \Biggl[ \sum
_{m=1}^N a_m \bG^{K_m}(r) +
\zeta\bigl(\bG^{K_1}(r) - \bG^{K_N}(r) \bigr) \Biggr]
\\
&\equiv& \bV+ \bG(1) + \bG(2) + \bG(3),\nonumber
\end{eqnarray}
where the $\bV^{k_m}$ and $\bV$ terms are associated with the process
$\bX(t)$ only, the $\bG^{K_m}(1)$ and $\bG(1)$ terms are related to the
noises $\varepsilon_i(t_\ell)$ only and the terms denoted by $\bG
^{K_m}(2)$, $\bG^{K_m}(3)$,
$\bG(2)$ and $\bG(3)$ depend on both $\bX(t)$ and~$\varepsilon
_i(t_\ell)$.

Now we may heuristically explain the basic ideas for proving Theorems~\ref{UnivariateConv}
and~\ref{Tthreshold} as follows. With the expression (\ref{VG}) we
prove the tail probability result for
$\tilde{\bGamma}$ in Theorem~\ref{UnivariateConv} by
establishing tail probabilities for these $\bV$ and $\bG$ terms in
the following
three propositions whose proofs will be given in %Section~\ref{appendixI} of
Appendix~\ref{appendixI}. %the subsequent subsections.

%
%pr5 #&#
\begin{proposition} \label{p_v}
Under the assumptions of Theorem~\ref{UnivariateConv}, we have for
$1\leq i,j \leq p$ and positive $d$ in a neighbor of $0$,
\[
P \bigl( \vert V_{ij}- \Gamma_{ij} \vert\geq d \bigr)
\leq C_1 n \exp\bigl\{-\sqrt{n} d^2/C_2
\bigr\}.
\]
\end{proposition}

%
%pr6 #&#
\begin{proposition}\label{pG2}
Under the assumptions of Theorem~\ref{UnivariateConv}, we have for
$1\leq i,j \leq p$ and positive $d$ in a neighbor of $0$,
\begin{eqnarray*}
P \bigl( \bigl\vert G_{ij}(2) \bigr\vert\geq d \bigr) &\leq&
C_1 n \exp\bigl\{-\sqrt{n} d^2/C_2 \bigr\},\\
P \bigl( \bigl\vert G_{ij}(3) \bigr\vert\geq d \bigr) &\leq&
C_1 n \exp\bigl\{-\sqrt{n} d^2 /C_2 \bigr\}.
\end{eqnarray*}
\end{proposition}

%
%pr7 #&#
\begin{proposition}\label{p_g1}
Under the assumptions of Theorem~\ref{UnivariateConv}, we have for
$1\leq i,j \leq p$ and positive $d$ in a neighbor of $0$,
\[
P \bigl( \bigl\vert G_{ij}(1) \bigr\vert\geq d \bigr) \leq
C_1 \sqrt{n} \exp\bigl\{- \sqrt{n} d^2/ C_2
\bigr\}.
\]
\end{proposition}
Because $V_{ij}$ are quadratic forms in the process $\bX(t_\ell)$ only,
we derive their tail probability in Proposition~\ref{p_v}
from the boundedness of the drift and volatility in condition A2; as
$G_{ij}(1)$ are quadratic forms in the noises
$\varepsilon_i(t_\ell)$ only, we establish the tail probability of
$G_{ij}(1)$ in Proposition~\ref{p_g1} from the subGaussianity
of $\varepsilon_i(t_\ell)$ imposed by condition A1; $G_{ij}(2)$ and
$G_{ij}(3)$ are bilinear forms in
$\bX(t_\ell)$ and $\varepsilon_i(t_\ell)$, thus we obtain the tail
probabilities for $G_{ij}(2)$ and $G_{ij}(3)$ in Proposition~\ref{pG2}
from the subGaussian tails of $\varepsilon_i(t_\ell)$ and $V_{ij}$ as
well as the independence between
$\varepsilon_i(t_\ell)$ and $\bX(t)$ given by condition~A1.
%The normality assumption can be relaxed to a subGaussian tail
%assumption, which is essentially required to obtain the convergence
%rate results in Theorems~\ref{UnivariateConv} and~\ref{Tthreshold}
%(see Supplemental Materials).
Since $\hat{\bGamma}$ is the matrix estimator obtained by thresholding
$\tilde{\bGamma}$, we use the tail probability
result in Theorem~\ref{UnivariateConv} %to analyze $\hat{\bGamma}$ and
%the properties of sparse matrices to handle $\bGamma$
and the sparsity of $\bGamma$ to analyze $\hat{\bGamma} - \bGamma$ and
control its matrix norm for proving Theorem~\ref{Tthreshold}.

\begin{pf*}{Proof of Theorem~\ref{UnivariateConv}}
From (\ref{VG}) we have
\[
P \bigl( \vert\tilde\Gamma_{ij} - \Gamma_{ij} \vert
\geq x \bigr) \leq P \bigl( \vert V_{ij} - \Gamma_{ij}
\vert\geq x/4 \bigr) + \sum_{r=1}^3 P
\bigl( \bigl\vert G_{ij}(r) \bigr\vert\geq x/4 \bigr), %P\left(
%P\left( \left| G_{ij}(2) \right| \geq x/4 \right) +
%P\left( \left| G_{ij}(3) \right| \geq x/4 \right),
\]
and thus the theorem is a consequence of Propositions~\ref{p_v}--\ref{p_g1}.
\end{pf*}

\begin{pf*}{Proof of Theorem~\ref{Tthreshold}}
Define
\begin{eqnarray*}
A_{ij} &=& \bigl\{ |\hat\Gamma_{ij} - \Gamma_{ij}|
\leq2\min\bigl\{|\Gamma_{ij}|,\varpi\bigr\} \bigr\},\qquad D_{ij} = (\hat
\Gamma_{ij} - \Gamma_{ij}) 1 \bigl(A_{ij}^c
\bigr), \\
\bD&=& (D_{ij})_{1 \leq i,j \leq p}.
\end{eqnarray*}
%
%A_{ij} & = & \left\{ |\hat\Gamma_{ij} - \Gamma_{ij}| \leq2\min\{|
%D_{ij} & = & (\hat\Gamma_{ij} - \Gamma_{ij}) 1(A_{ij}^c), \bD=
%(D_{ij})_{1 \leq i,j \leq p}.
As the matrix norm of a symmetric matrix is bounded by its $\ell
_1$-norm, then
%
%
%e27 #&#
\begin{equation}
\label{Dterm} E\Vert\hat\bGamma- \bGamma\Vert_2^2 \leq
E\Vert\hat\bGamma- \bGamma\Vert_1^2 \leq2 E\Vert\hat
\bGamma- \bGamma- \bD\Vert_1^2 + 2 E\Vert\bD\Vert
_1^2.
\end{equation}
We can bound $ E\Vert \hat\bGamma- \bGamma- \bD\Vert^2_1$ as follows:
\begin{eqnarray*}
 &&E\Vert\hat\bGamma- \bGamma- \bD\Vert_1^2\\
 &&\qquad = E \Biggl[
\max_{1\leq j \leq p} \sum_{i=1}^p
| \hat\Gamma_{ij} - \Gamma_{ij}| 1 \bigl(|\hat
\Gamma_{ij} - \Gamma_{ij}| \leq2\min\bigl\{|\Gamma_{ij}|,
\varpi\bigr\} \bigr) \Biggr]^2
\\
&&\qquad\leq E \Biggl[ \max_{1\leq j \leq p} \sum
_{i=1}^p 2 |\Gamma_{ij}| 1\bigl(|
\Gamma_{ij}| < \varpi\bigr) \Biggr]^2 \\
&&\qquad\quad{}+ E \Biggl[ \max
_{1\leq j \leq p} \sum_{i=1}^p 2
\varpi1\bigl(|\Gamma_{ij}| \geq\varpi\bigr) \Biggr]^2
\\
&&\qquad\leq %E\left( 4 \Psi\pi_n(p) \varpi^{1-q} \right) ^2 =
8 E \bigl[\Psi^2 \bigr] \pi_n^2(p)
\varpi^{2(1-q)} \leq C \pi_n^2(p) \bigl(
n^{-1/4} \sqrt{\log p} \bigr)^{2-2 q},
\end{eqnarray*}
where the second inequality is due to the fact that the sparsity of
$\bGamma$
implies
\begin{eqnarray*}
\max_{1 \leq j \leq p} \sum_{i=1}^p
1\bigl(|\Gamma_{ij}| \geq\varpi\bigr) &\leq&\Psi\pi_n(p)
\varpi^{-q},\\ \max_{1 \leq j \leq p} \sum
_{i=1}^p |\Gamma_{ij}| 1\bigl(|
\Gamma_{ij}| < \varpi\bigr) &\leq&\Psi\pi_n(p)
\varpi^{1-q},
\end{eqnarray*}
which are the respective bounds on the number of those entries on each
row with absolute values larger than or equal to $\varpi$
and the sum of those absolute entries on each row with magnitudes less
than $\varpi$; see Lemma~1 in \citet{WanZou10}.
The rest of the proof is to show that $E\Vert\bD\Vert_1^2 =
O(n^{-2})$, a
negligible term. Indeed,
the threshold rule indicates that $\hat{\Gamma}_{ij}=0$ if $|\tilde
{\Gamma}_{ij}| < \varpi$ and
$\hat{\Gamma}_{ij}=\tilde{\Gamma}_{ij}$ if $|\tilde{\Gamma}_{ij}|
\geq
\varpi$, thus
\begin{eqnarray*}
 E\Vert\bD\Vert_1^2 %= E \left[ \max_{1\leq j \leq p} \sum_{i=1}^p |
% 1(|\hat\Gamma_{ij} - \Gamma_{ij}| > 2\min\{|\Gamma_{ij}|,\varpi\})
%&
&\leq& p \sum_{i,j=1}^p
E \bigl[ |\Gamma_{ij}|^2 1 \bigl( |\Gamma_{ij}|
> 2 \min\bigl\{|\Gamma_{ij}|, \varpi\bigr\} \bigr) 1( \hat{
\Gamma}_{ij} =0 ) \bigr]
\\
&&{} + p \sum_{i,j}^p E \bigl[ |\tilde
\Gamma_{ij} - \Gamma_{ij}|^2 1 \bigl(|\tilde
\Gamma_{ij} - \Gamma_{ij}| > 2\min\bigl\{|\Gamma_{ij}|,
\varpi\bigr\} \bigr) 1(\hat{\Gamma}_{ij}=\tilde\Gamma_{ij})
\bigr] %\\ &
\\
&\equiv& I_1 + I_2.
\end{eqnarray*}

For term $I_1$, we have
\begin{eqnarray*}
I_1  %= p \sum_{i,j=1}^p E\left[ |\Gamma_{ij}|^2 1( |\Gamma_{ij}| > 2
% 1( |\tilde{\Gamma}_{ij}| < \varpi) \right] \\ &
&\leq& p \sum
_{i,j=1}^p E \bigl[ |\Gamma_{ij}|^2
1\bigl( |\tilde\Gamma_{ij} - \Gamma_{ij}| > \varpi\bigr) \bigr]
\leq C p \sum_{i,j=1}^p P\bigl(|
\tilde\Gamma_{ij} - \Gamma_{ij}| > \varpi\bigr)
\\
& \leq &C p^3 \exp\bigl\{ \log n %(n/\varpi)
- \sqrt{n}
\varpi^2/\varsigma_0 \bigr\} % \leq C p^3 \exp\left\{ \log n -
\leq C n^{-2},
\end{eqnarray*}
where the third inequality is from Theorem~\ref{UnivariateConv}, and
the last
inequality is due to $\varpi= \hbar n^{-1/4} \sqrt{ \log(np)}$ with
$\hbar^2/\varsigma_0 > 4$.

On the other hand, we can bound term $I_2$ as follows:
\begin{eqnarray*}
I_2  %= p \sum_{i,j=1}^p E \left[ |\tilde\Gamma_{ij} - \Gamma_{ij}|^2
% 1(|\tilde\Gamma_{ij} - \Gamma_{ij}| > 2\min\{|\Gamma_{ij}|,\varpi\})
% 1(|\tilde\Gamma_{ij}| \geq\varpi)\right]\\
% & = p\sum_{i,j=1}^p E \left[ |\tilde\Gamma_{ij} - \Gamma_{ij}|^2
% 1(|\tilde\Gamma_{ij} - \Gamma_{ij}| > 2\min\{|\Gamma_{ij}|,\varpi
% 1(|\tilde\Gamma_{ij}| \geq\varpi)\right]\\
% & + p \sum_{i,j=1}^p E \left[ |\tilde\Gamma_{ij} - \Gamma_{ij}|^2
% 1(|\tilde\Gamma_{ij} - \Gamma_{ij}| > 2\min\{|\Gamma_{ij}|,\varpi
% 1(|\tilde\Gamma_{ij}| \geq\varpi)\right]\\
%&
&\leq& p\sum_{i,j=1}^p E \bigl[ |
\tilde\Gamma_{ij} - \Gamma_{ij}|^2 1\bigl(|\tilde
\Gamma_{ij} - \Gamma_{ij}| > \varpi\bigr) \bigr] \\
&&{}+ p\sum
_{i,j=1}^p E \bigl[ |\tilde\Gamma_{ij} -
\Gamma_{ij}|^2 1\bigl(|\Gamma_{ij}| < \varpi/2, |
\tilde\Gamma_{ij}| \geq\varpi\bigr) \bigr]
\\
& \leq& 2p \sum_{i,j=1}^p E \bigl[ |\tilde
\Gamma_{ij} - \Gamma_{ij}|^2 1\bigl(|\tilde
\Gamma_{ij} - \Gamma_{ij}| > \varpi/2 \bigr) \bigr]
\\
& \leq& 2p \sum_{i,j=1}^p \bigl\{ E \bigl[
|\tilde\Gamma_{ij} - \Gamma_{ij}|^4 \bigr] P
\bigl(|\tilde\Gamma_{ij} - \Gamma_{ij}| > \varpi/2 \bigr) \bigr\}
^{1/2}
\\
& \leq& C p^3 \exp\bigl\{ \log n/2 %(n/\varpi)/2
- \sqrt{n}
\varpi^2 / (8 \varsigma_0) \bigr\} % \leq C p^3 \exp\left\{ -\hbar^2
\leq C
n^{-2},
\end{eqnarray*}
where the third inequality is due to H\"older's inequality,
the fourth inequality is from Theorem~\ref{UnivariateConv} and
%(\ref{tildeG4}) below, Lemma~\ref{tildeG4} below,
%Under the assumptions of Theorem~\ref{UnivariateConv}, we have for
%$1\leq i,j \leq p$,
%
%
%e28 #&#
\begin{equation}
\label{tildeG4} \max_{1 \leq i, j \leq p} E \bigl[ |\tilde{
\Gamma}_{ij} - \Gamma_{ij}|^4 \bigr] \leq C
\end{equation}
and the last inequality is due to the fact that
$\varpi=\hbar n^{-1/4} \sqrt{ \log(np)}$ with $\hbar^2/(8 \varsigma_0)
> 3$.

To complete the proof we need to show (\ref{tildeG4}).
%To prove (\ref{tildeG4}), since $\Gamma_{ij}$ is bounded, we need to
%show
%only that $ E \left[ |\tilde{\Gamma}_{ij}|^4 \right] $ is bounded.
As in \citet{ZhaMykAtS05}, we adjust $\tilde{\bGamma}{}^{K_m}$ to account
for the noise variances. Let
%
%
%e29 #&#
\begin{equation}
\label{eta0} \tilde{\boeta}= %(\tilde{\eta}_{ij}) =
\operatorname{diag}(\tilde\eta_1,\ldots,
\tilde\eta_p), \qquad \tilde\eta_i=\frac{1}{2n}\sum
_{\ell=2}^{n} \bigl[Y_i(t_{\ell})-Y_i(t_{\ell-1})
\bigr]^2,
\end{equation}
and define
%
%
%e30 #&#
\begin{equation}
\label{Gammatilde*} \tilde{\bGamma}{}^{*K_m}=\tilde{\bGamma}{}^{K_m}-2
\frac{n-K_m+1}{K_m}\tilde{\boeta},
\end{equation}
which are the average realized volatility matrix (ARVM) estimators
where the convergence rates
for any finite moments of $\tilde\Gamma{}^{*K_m}_{ij}-\Gamma_{ij}$ are
derived in
Wang and Zou [(\citeyear{WanZou10}), Theorem~1].
Applying Theorem~1 of \citet{WanZou10} to the fourth moment of
$\tilde\Gamma{}^{*K_m}_{ij}-\Gamma_{ij}$, %$\tilde{\bGamma}^{*K_m}$
we have for $1 \leq i,j \leq p$ and $1 \leq m \leq N$,
%
%
%e31 #&#
\begin{eqnarray}
\label{Gammatilde*1}&& E \bigl(\bigl|\tilde\Gamma{}^{*K_m}_{ij}-
\Gamma_{ij}\bigr|^4 \bigr)
\nonumber
\\[-8pt]
\\[-8pt]
\nonumber
&&\qquad\leq C \bigl[ \bigl(K_m
n^{
-1/2} \bigr)^{-4}+ K_m^{-2}+(n/K_m)^{-2}+K_m^{-4}+n^{-2}
\bigr] \leq C.
\end{eqnarray}
%
%Again with bounded $\Gamma_{ij}$ we have for $1 \leq i,j \leq p$ and
%$1 \leq m \leq N$,
% E \left[ |\tilde{\Gamma}^{*K_m}_{ij}|^4 \right] \leq C.
From (\ref{MSRVM}), (\ref{coeff}) and (\ref{Gammatilde*}) together with
simple algebraic manipulations we can express $\tilde{\bGamma}$ by
$\tilde{\bGamma}{}^{*K_m}$ as follows:
\[
\tilde{\bGamma} =\sum_{m=1}^N
a_m\tilde{\bGamma}{}^{*K_m}+ \zeta\bigl(\tilde{
\bGamma}{}^{*K_1}- \tilde{\bGamma}{}^{*K_N} \bigr),
\]
and thus
%
%
%e32 #&#
\begin{equation}
\label{Gammahattilde*} \tilde{\bGamma} -\bGamma=\sum
_{m=1}^N a_m \bigl( \tilde{
\bGamma}{}^{*K_m} - \bGamma\bigr) + \zeta\bigl[ \bigl(\tilde{
\bGamma}{}^{*K_1} - \bGamma\bigr) - \bigl(\tilde{\bGamma}{}^{*K_N} -
\bGamma\bigr) \bigr].
\end{equation}

Combining (\ref{Gammatilde*1}) and (\ref{Gammahattilde*}) and using
(\ref{coeff}) we conclude for $1 \leq i,j \leq p$,
\begin{eqnarray*}
&&E \bigl[ |\tilde{\Gamma}_{ij} -\Gamma_{ij}|^4
\bigr] \\
&&\qquad\leq(N+2)^3 \Biggl[ \sum_{m=1}^N
a_m^4 E \bigl(\bigl|\tilde{\Gamma}{}^{*K_m}_{ij}
- \Gamma_{ij}\bigr|^4 \bigr) \\
&&\hspace*{75pt}\qquad{}+ \zeta^4 E \bigl(\bigl|
\tilde{\bGamma}{}^{*K_1} - \Gamma_{ij}\bigr|^4 + \bigl|\tilde
{\bGamma}{}^{*K_N} - \Gamma_{ij}\bigr|^4 \bigr) \Biggr]\\
&&\qquad
\leq C.
\end{eqnarray*}
\upqed\end{pf*}

%s6.2 #&#
\subsection{\texorpdfstring{Proofs of Theorems \protect\ref{lowerbound} and \protect\ref{Operlowerbdthm}}
{Proofs of Theorems 3 and 4}}\label{lowerboundproofssc}
Section~\ref{transformationsec} shows that Theorem~\ref{lowerbound} is
a consequence of
Theorem~\ref{Operlowerbdthm}. The proof of Theorem~\ref
{Operlowerbdthm} is similar to
but much more involved than the proof of Theorem~2 in \citet{CaiZho} which considered
only i.i.d. observations. It contains four major steps. In the first
step we construct
in detail a finite subset $\mathcal{F}_{\ast}$ of the parameter space
$\mathcal{G}_{q}(\pi_{n}(p))$ in the minimax problem $\mathcal
{Q}_{q}(\pi_{n}(p))$
such that the difficulty of estimation over $\mathcal{F}_{\ast}$ is
essentially the
same as that of estimation over $\mathcal{G}_{q}(\pi_{n}(p))$,
where $\mathcal{G}_{q}(\pi_{n}(p))$ is the class of constant matrices
$\bGamma$ satisfying (\ref{sparseparaspace}) and $\| \bGamma\|_2
\leq
\tau$
for constant $\tau>0$. The second step applies the lower bound
argument in
Cai and Zhou [(\citeyear{CaiZho}), Lemma~3]
%the application of Lemma~3 in Cai and Zhou (2012)
to the carefully constructed parameter set $\mathcal{F}_{\ast}$.
In the third step we calculate the factor $\alpha$ defined in
(\ref{alpha1}) below and the total variation affinity between two
average of
products of $n$ independent but nonidentically distributed
multivariate normals.
%Bounding the affinity is technically involved. The detailed proofs are
%then deferred to subsequent sections.
%Section~\ref{lemmaaffsec}.
%Since only matrix norm is used in the rest of the paper we will write
%$\Vert\cdot\Vert_2$ simply as $\Vert\cdot\Vert$.
%For probability $\mathbb P$, we denote by $\Vert\mathbb P \Vert$ its
%total variation.
The final step combines together the results in steps 2 and 3 to obtain
the minimax
lower bound.

\textit{Step \textup{1:} Construct parameter set ${\cal
F}_*$.} Set $
r=\lceil p/2\rceil$, where $ \lceil x \rceil$ denotes the
smallest integer greater than or equal to $x$, and let $B$ be the
collection of all row
vectors $b= ( v_{j} ) _{1\leq j\leq p}$ such that $v_{j}=0$
for $
1\leq j\leq p-r$ and $v_{j}=0$ or $1$ for $p-r+1\leq j\leq p$ under the
constraint $\llVert b\rrVert_{0}=k$ (to be specified later).
Each element $\lambda=(b_{1},\ldots,b_{r})\in B^{r}$ is treated as an $
r\times p$ matrix with the $i$th row of $\lambda$ equal to $b_{i}$.
Let $
\Delta= \{ 0,1 \} ^{r}$. Define $\Lambda\subset B^{r}$ to be the
set of all elements in $B^{r}$ such that each column sum is less than or
equal to $2k$. For each $b\in B$ and each $1\leq m\leq r$, define a
$p\times
p$ symmetric matrix $A_{m}(b)$ by making the $m$th row of $A_{m}(b)$ equal
to $b$, $m$th column equal to $b^{T}$ and the rest of the entries $0$. Then
each component $\lambda_{i}$ of $\lambda=(\lambda_{1},\ldots,\lambda
_{r})\in\Lambda$ can be uniquely associated with a $p\times p$ matrix $
A_{i}(\lambda_{i})$. Define $\Theta=\Delta\otimes\Lambda$, and let $
\epsilon_{n,p}\in\mathbb{R}$ be fixed (the exact value of $\epsilon_{n,p}
$ will be chosen later). For each $\theta=(\gamma,\lambda)\in\Theta$
with\vadjust{\goodbreak} $\gamma=(\gamma_{1},\ldots,\gamma_{r})\in\Delta$ and $\lambda
=(\lambda_{1},\ldots,\lambda_{r})\in\Lambda$, we associate $\theta
= (
\gamma_{1},\ldots,\gamma_{r},\lambda_{1},\ldots,\lambda_{r} ) $ with a
volatility matrix $\bGamma(\theta)$ by
%
%
%e33 #&#
\begin{equation}
\bGamma(\theta)=\bI_{p}+\epsilon_{n,p}\sum
_{m=1}^{r}\gamma_{m}A_{m}(
\lambda_{m}). \label{lowbdspace}
\end{equation}
For simplicity %Without loss of generality
we assume that $\tau>1$ in the definition of the parameter space
$\mathcal{G}_{q}(\pi_{n}(p))$ for the minimax problem $\mathcal
{Q}_{q}(\pi_{n}(p))$;
%the assumption (\ref{sparseparaspace}),
otherwise we replace $\bI_{p}$ in (\ref{lowbdspace}) by
$C \bI_{p}$ with a small constant $C>0$. Finally we define $\mathcal
{F}_{\ast}$
to be a collection of covariance matrices as
%
%
%e34 #&#
\begin{equation}
\mathcal{F}_{\ast}= \Biggl\{ \bGamma(\theta)\dvtx\bGamma(\theta
)=\bI
_{p}+\epsilon_{n,p}\sum_{m=1}^{r}
\gamma_{m} A_m(\lambda_{m}), %\lambda_{m}(\lambda_{m}),
\theta=(\gamma,\lambda)\in\Theta\Biggr\}. \label{F*}
\end{equation}
Note that each matrix $\bGamma\in\mathcal{F}_{\ast}$ has value $1$
along the main
diagonal and contains an $r\times r$ submatrix, say, $A$, at the upper
right corner, $A^{T}$ at the lower left corner and $0$ elsewhere; each row
of the submatrix $A$ is either identically $0$ (if the corresponding
$\gamma
$ value is $0$) or has exactly $k$ nonzero elements with value
$\epsilon
_{n,p}$.

Now we specify the values of $\epsilon_{n,p}$ and $k$:
%
%
%e35 #&#
\begin{equation}
\label{k-epsilon} \epsilon_{n,p}=\upsilon\biggl( \frac{\log
p}{\sqrt{n}}
\biggr) ^{1/2},\qquad  k= \biggl\lceil\frac{1}{2} \pi_{n}(p)
\epsilon_{n,p}^{-q} \biggr\rceil-1,
\end{equation}
where $\upsilon$ is a fixed small constant that
%Let $\epsilon_{n,p}=\upsilon\left( \frac{\log p}{\sqrt{n}}\right)
%^{1/2}$ for a fixed
%small constant $\upsilon$, and $k=\left\lceil\frac{1}{2}\pi_{n}(p)
%which implies $\max_{j\leq p}\sum_{i\neq j}|\Gamma_{ij} %\sigma
%_{ij}|^{q}\leq2k\epsilon_{n,p}^{q}\leq\pi_{n}(p)$.
we require
%0<\upsilon<\left[ \min\left\{ \frac{1}{3},\tau-1\right\} \frac{1}{M}
% \upsilon^{2}<\frac{\kappa(\beta-1)}{\beta}
%% \upsilon^{2}<\frac{2\kappa(\beta-1)}{\pi\beta}.
%%see Equation (\ref{v1}) for why $\upsilon$ needs to satisfy the
%second inequlaity in (\ref{v}).
%
%
%e36 #&#
\begin{equation}
0<\upsilon< \biggl[ \min\biggl\{ \frac{1}{3},\tau-1 \biggr\}
\frac
{1}{\aleph} \biggr] ^{{1}/{(1-q)}} \label{v}
\end{equation}
and
%
%
%e37 #&#
\begin{equation}
0<\upsilon^{2}<\frac{\beta-1}{27 c_{\kappa}\beta}, \label{v1}
\end{equation}
where $c_{\kappa} = (2 \kappa)^{-1}$ satisfies %is a constant such
%that
%
%
%e38 #&#
\begin{equation}
\label{ck} \sum_{l=1}^n
a_l^{-2} \leq c_{\kappa}\sqrt{n},
\end{equation}
since
\begin{eqnarray*}
\sum_{l=1}^n a_l^{-2}
&\leq&\int_0^n \biggl[1 + 4
\kappa^2 n \sin^2 \biggl( \frac{\pi x}{2(n+1)} \biggr)
\biggr]^{-2} \,dx \leq\frac{n+1}{\pi\kappa\sqrt{n}} \int
_0^\infty
\bigl[1 + v^2 \bigr]^{-2} \,dv
\\
&=&
\frac{\sqrt{n}+1/\sqrt{n}}{ 4\kappa}.
\end{eqnarray*}
Note that $\epsilon_{n,p}$ and $k$ satisfy
$\max_{j\leq p}\sum_{i\neq j}|\Gamma_{ij} |^{q}\leq2k\epsilon
_{n,p}^{q}\leq\pi_{n}(p)$,
%
%
%e39 #&#
\begin{equation}
2k\epsilon_{n,p}\leq\pi_{n}(p)\epsilon_{n,p}^{1-q}
\leq\aleph\upsilon^{1-q}<\min\bigl\{ \tfrac{1}{3},\tau-1 \bigr
\}, \label{1normbound}
\end{equation}
and consequently every $\bGamma(\theta)$ is diagonally dominant and
positive definite, and $\Vert\bGamma(\theta)\Vert_2 \leq\Vert
\bGamma
(\theta)\Vert_{1}\leq2k\epsilon_{n,p}+1<\tau$. Thus we have
$\mathcal{F}%
_{\ast}\subset\mathcal{G}_{q}(\pi_{n}(p))$.\vadjust{\goodbreak}

\textit{Step \textup{2:} Apply the general lower bound argument}.
Let $\bU_{l}$ be independent with
\[
\bU_{l}\sim N \bigl( 0,\bGamma( \theta) + ( a_{l}-1 )
\bI_{p} \bigr),
\]
where $l=1, \ldots, n$, $\theta\in\Theta$, and we denote the joint
distribution by
$P_{\theta}$.
Applying Lemma~3 in \citet{CaiZho} to the parameter space $\Theta$,
we have
%
%
%e40 #&#
\begin{equation}
\inf_{\check{\bGamma}}\max_{\theta\in\Theta} \mathbb{E}_{\theta
}
\bigl\llVert\check{ \bGamma}-\bGamma(\theta) \bigr\rrVert_2^{2}
\geq\alpha\cdot\frac
{r}{8}\cdot\min_{1\leq i\leq r}\llVert
\bar{\mathbb{P}}_{i,0}\wedge\bar{\mathbb{P}}%
_{i,1}\rrVert, \label{lowerbound*}
\end{equation}
where we use $\Vert\mathbb P \Vert$ to denote the total variation of
$\mathbb P$,
%
%
%e41 #&#
\begin{eqnarray}\label{alpha1}
\alpha&\equiv&\min_{ \{ (\theta,\theta^{\prime})\dvtx H(\gamma
(\theta
),\gamma(\theta^{\prime}))\geq1 \} }\frac{\llVert\bGamma
(\theta
)-\bGamma(\theta^{\prime})\rrVert_2^{2}}{H(\gamma(\theta
),\gamma
(\theta^{\prime}))},
\nonumber
\\[-8pt]
\\[-8pt]
\nonumber
  H \bigl(\gamma(
\theta),\gamma\bigl(\theta^{\prime} \bigr) \bigr) &= &\sum
_{i=1}^r\bigl |\gamma_i(\theta) -
\gamma_i \bigl(\theta^{\prime} \bigr)\bigr|
\end{eqnarray}
and
%
%
%e42 #&#
\begin{equation}
\bar{\mathbb{P}}_{i,a}=\frac{1}{2^{r-1}D_{\Lambda}}\sum
_{\theta\in
\Theta}%
\mathbb{P}_{\theta}\cdot\bigl\{
\theta\dvtx\gamma_{i}(\theta)=a \bigr\}, \label{avepi}
\end{equation}
where $a\in\{ 0,1 \} $ and $D_{\Lambda}=\mathrm{Card} \{
\Lambda\} $.

\textit{Step \textup{3:} Bound the affinity and per comparison
loss}. We need to bound the two factors $\alpha$ and $\min_{i}\llVert
\bar{\mathbb{P}}_{i,0}\wedge\bar{\mathbb{P}}_{i,1}\rrVert$ in
(\ref
{lowerbound*}). A lower bound for $\alpha$ is given by the following
proposition
whose proof is the same as that of Lemma~5 in \citet{CaiZho}.\vspace*{-2pt}
%straightforward and thus omitted.

%
%pr8 #&#
\begin{proposition}
\label{dffbd} For $\alpha$ defined in equation (\ref{alpha1}) we have
\[
\alpha\geq\frac{(k\epsilon_{n,p})^{2}}{p}.
\]
\end{proposition}

A lower bound for $\min_{i}\llVert\bar{\mathbb{P}}
_{i,0}\wedge\bar{\mathbb{P}}_{i,1}\rrVert$ is provided by the
proposition below.
Since its proof is long and very much involved, the proof details %of
%the proof
are collected in %Section~\ref{lemmaaffsec} of
Appendix~\ref{lemmaaffsec}.\vspace*{-2pt}
%The proof of \protect\ref{affbd} is similar to that of Lemma~6 in Cai
%and Zhou (2012).

%
%pr9 #&#
\begin{proposition}
\label{affbd} Let $\bU_{l}$ be independent with $\bU_{l}\sim N (
0,\bGamma(\theta)
+ ( a_{l}-1 ) \bI_{p} ) $, $l=1, \ldots, n$, with $\theta
\in\Theta$ and denote
the joint distribution by $\mathbb{P}_{\theta}$. For $a\in\{0,1\}$
and $
1\leq i\leq r$, define $\bar{\mathbb{P}}_{i,a}$ as in (\ref{avepi}). Then
there exists a constant $C_{1}>0$ such that
\[
\min_{1\leq i\leq r}\llVert\bar{\mathbb{P}}_{i,0}\wedge
\bar{\mathbb{P}}%
_{i,1}\rrVert\geq C_{1}
\]
uniformly over $\Theta$.
\end{proposition}

\textit{Step \textup{4:} Obtain the minimax lower bound.}
We obtain the minimax lower bound for estimating $\bGamma$ %volatility
%matrices
over $\mathcal{G}_{q}(\pi_{n}(p))$ by combining\vadjust{\goodbreak} together (\ref
{lowerbound*}) and
the bounds in Propositions~\ref{dffbd} and~\ref{affbd},
\begin{eqnarray*}
\inf_{\check{\bGamma}}\sup_{\mathcal{G}_{q}(\pi_{n}(p))}\mathbb
{E}\llVert
\check{ \bGamma}-\bGamma\rrVert_2^{2} &\geq& \inf
_{\check{\bGamma}} \max_{\bGamma(\theta)\in\mathcal{F}%
_{\ast}} \mathbb{E}_{\theta}
\bigl\llVert\check{\bGamma}-\bGamma( \theta) \bigr\rrVert_2^{2}
\geq\frac{ ( k\epsilon_{n,p} ) ^{2}}{p}\cdot\frac{r}{8}\cdot C_{1}
\\
&\geq&\frac{C_{1}}{16}(k\epsilon_{n,p})^{2} =
C_{2} \pi^2_{n}(p) \bigl( n^{-1/4}
\sqrt{\log p} \bigr)^{2-2q}
\end{eqnarray*}
for some constant $C_2>0$. %\hbox{\vrule width 4pt height 6pt depth 1.5pt}

%s6.3 #&#
\subsection{\texorpdfstring{Proof of (\protect\ref{lowerbd-d-norm}) for optimal
convergence rate under general matrix norm}
{Proof of (16) for optimal convergence rate under general matrix norm}}

%The upper bound holds under the general matrix $\ell_{d}$ norm for $1
The Riesz--Thorin interpolation theorem [\citet{Tho48}] implies for
$1\leq
d_{1}\leq d \leq d_{2}\leq\infty$,
%
%
%e43 #&#
\begin{equation}
\label{RT} \llVert\bA\rrVert_{d}\leq\max\bigl\{ \llVert\bA
\rrVert_{d_{1}}, \llVert\bA\rrVert_{d_{2}} \bigr\}. % \ for\ all\ 1
\end{equation}
Set $d_{1}=1$ and $d_{2}=\infty$, then (\ref{RT}
) yields $\llVert\bA\rrVert_{d}\leq\max\{ \llVert
\bA\rrVert_{1},\llVert\bA\rrVert_{\infty} \} $ for \mbox{$
1\leq d\leq\infty$}. When $\bA$ is symmetric, (\ref{ell-1infty-norm})
shows that $\llVert\bA\rrVert
_{1}=\llVert\bA\rrVert_{\infty}$. Then immediately we have $
\llVert\bA\rrVert_{d}\leq\llVert\bA\rrVert_{1}$,
which means that for a symmetric matrix estimator,
an upper bound under the matrix $\ell_{1}$ norm is also an upper bound under
the general matrix $\ell_{d}$ norm. Thus, as $\hat\bGamma$ is symmetric,
Theorem~\ref{Tthreshold} indicates that for $ 1 \leq d \leq\infty$,
\[
\sup_{\mathcal{P}_{q}(\pi_{n}(p))} \mathbb{E} \llVert\hat\bGamma
- \bGamma\rrVert
_d^2 \leq C^* \bigl[ \pi_n(p) \bigl(
n^{-1/4} \sqrt{\log p} \bigr) ^{1-q} \bigr]^2.
\]
Now consider the lower bound under %The lower bound can be extended to
the general matrix $\ell_{d}$ norm for $1\leq d\leq\infty$.
We will show
%
%
%e44 #&#
\begin{eqnarray}
\label{symmetric-lower-bound} \inf_{\check{\bGamma}_s}\sup
_{\mathcal{P}_{q}(\pi_{n}(p))} \mathbb{E}\llVert\check{\bGamma}_s-
\bGamma\rrVert_d^{2} &\geq&\inf_{\check{\bGamma}}
\sup_{\mathcal{P}_{q}(\pi_{n}(p))} \mathbb{E}\llVert\check
{\bGamma}-\bGamma\rrVert
_d^{2}
\nonumber
\\[-8pt]
\\[-8pt]
\nonumber
& \geq&\frac{1}{4} \inf_{\check{\bGamma}_s}
\sup_{\mathcal{P}_{q}(\pi_{n}(p))} \mathbb{E}\llVert\check
{\bGamma}_s-
\bGamma\rrVert_d^{2},
\end{eqnarray}
where $\check{\bGamma}$ denotes any matrix estimators of $\bGamma$, and
$\check{\bGamma}_s$ any
symmetric matrix estimators of $\bGamma$.
(\ref{symmetric-lower-bound}) indicates that it is enough to consider
estimators of symmetric matrices.

For symmetric $\bA$, (\ref{ell-12infty-norm}) shows that $\|\bA\|_2
\leq\|\bA\|_1 = \|\bA\|_\infty$.
For $d \in(1, \infty)$, $1/d + (d-1)/d=1$, by duality we have
$\llVert\bA\rrVert_{d}=\llVert
\bA\rrVert_{{d}/{(d-1)}}$. Also since $2$ is always between $d$
and $d/(d-1)$, applying (\ref{RT}) we
obtain that $\llVert\bA\rrVert_{2}\leq\max\{ \llVert
\bA\rrVert
_{d}, \llVert\bA\rrVert_{{d}/{(d-1)}} \} = \llVert
\bA\rrVert_{d}$.
This means that within the class of symmetric matrix estimators, a
lower bound under the matrix $\ell_{2}$ norm
is also a lower bound under the general matrix $\ell_{d}$ norm.
%Since $\bGamma$ are symmetric matrices in the minimax problem $
Thus (\ref{symmetric-lower-bound})
%together with $\left\Vert\bA\right\Vert_{2}\leq\left\Vert A\right
%for symmetric $A$ show that a lower bound under the spectral norm is
%also a lower bound under the
%general matrix $\ell_{d}$ norm for all $1\leq d\leq\infty$. Thus,
and Theorem~\ref{lowerbound} together imply that for
$ 1 \leq d \leq\infty$,
\[
\inf_{\check{\bGamma}}\sup_{\mathcal{P}_{q}(\pi_{n}(p))} \mathbb
{E}\llVert
\check{\bGamma}-\bGamma\rrVert_d^{2} \geq\frac{C_*}{4}
\bigl[ \pi_{n}(p) \bigl( n^{-1/4} \sqrt{\log p}
\bigr)^{1-q} \bigr]^2.
\]
To complete the proof we need to prove (\ref{symmetric-lower-bound}).
The first inequality of (\ref{symmetric-lower-bound}) is obvious.
%indicates that it is enough to consider estimators of symmetric
%matrices,
%imply a lower bound under the spectral norm is also a lower bound for
%the
%general matrix $\ell_{d}$ norm for all $1\leq d\leq\infty$. The
%argument for Equation (\ref{symmest}) is as follows.
For a given matrix estimator $\check{\bGamma}$ we project it onto the
parameter space of the minimax
problem $\mathcal{P}_{q}(\pi_{n}(p))$ by minimizing the matrix $\ell_d$
norm of
$\check{\bGamma} - \bGamma_*$ over all $\bGamma_*$ in the parameter
space. Denote its projection by $\check{\bGamma}_p$.
Since the parameter space consists of symmetric matrices, $\check
{\bGamma}_p$ is symmetric. Hence
%$\check{\bGamma}_p= \arg\min_{\bGamma} \left\Vert\check{\bGamma}-
%
\begin{eqnarray*}
&& \inf_{\check{\bGamma}_s}\sup_{\mathcal{P}_{q}(\pi_{n}(p))}
\mathbb{E}\llVert
\check{\bGamma}_s-\bGamma\rrVert_d^{2}
\\
&&\qquad\leq%
\sup_{\mathcal{P}_{q}(\pi_{n}(p))} \mathbb{E}\llVert
\check{\bGamma}_p-\bGamma\rrVert_d^{2}
\\
&&\qquad\leq 2 \sup_{\mathcal{P}_{q}(\pi_{n}(p))} \bigl[ \mathbb
{E}\llVert\check{
\bGamma}_p- \check{\bGamma} \rrVert_d^{2} +
\mathbb{E}\llVert\check{\bGamma}-\bGamma\rrVert_d^{2}
\bigr]
\\
&&\qquad \leq 2 \sup_{\mathcal{P}_{q}(\pi_{n}(p))} \bigl[ \mathbb
{E}\llVert\bGamma- \check{
\bGamma} \rrVert_d^{2} + \mathbb{E}\llVert\check{
\bGamma}-\bGamma\rrVert_d^{2} \bigr]
\\
&&\qquad \leq4 \sup_{\mathcal{P}_{q}(\pi_{n}(p))} \mathbb{E}\llVert
\check{\bGamma}-\bGamma
\rrVert_d^{2},
\end{eqnarray*}
where the second inequality is from the triangle inequality and the
third one follows from the definition of $\check{\bGamma}_p$. Since
the above inequality holds for every~$\check{\bGamma}$, we have
\begin{eqnarray*}
&&\inf_{\check{\bGamma}_s}\sup_{\mathcal{P}_{q}(\pi_{n}(p))}
\mathbb{E}\llVert
\check{\bGamma}_s-\bGamma\rrVert_d^{2}\\
&&\qquad \leq
4 \inf_{\check{\bGamma}} \sup_{\mathcal{P}_{q}(\pi_{n}(p))}
\mathbb{E}\llVert
\check{\bGamma}-\bGamma\rrVert_d^{2},
\end{eqnarray*}
which is equivalent to the second inequality of (\ref{symmetric-lower-bound}).

%
% \hangafter=1 }
%%\renewcommand\baselinestretch{0.9}
%%\tiny
%

%Yazhen Wang and Minjing Tao \\
%Department of Statistics \\
%University of Wisconsin-Madison\\
%Madison, WI 53706
%Harrison H. Zhou \\
%Department of Statistics \\
%Yale University \\
%New Haven, CT 06511
%
%Minjing Tao and Yazhen Wang & Harrison H. Zhou \\
%Department of Statistics & Department of Statistics \\
%University of Wisconsin-Madison \hspace*{1in} & Yale University\\
%Madison, WI 53706 & New Haven, CT 06511
%%

%
%sA #&#
\begin{appendix}\label{app}
\renewcommand{\thesection}{\Roman{section}}
\section{\texorpdfstring{Proofs of Propositions \lowercase{\protect\ref{p_v}}--\lowercase{\protect\ref{p_g1}}}
{Proofs of Propositions 5--7}}\label{appendixI}
\renewcommand{\thesubsection}{\Roman{section}.\arabic{subsection}}
%sA.1 #&#
\subsection{\texorpdfstring{Proof of Proposition \protect\ref{p_v}}{Proof of Proposition 5}} \label{sectionV}

%Under the assumptions of Theorem~\ref{UnivariateConv}, we have for
%$1\leq i,j \leq p$,
%P\left( \left|V_{ij}- (\Gamma_{ij}-\Gamma_{ij}^*)
% \right|\geq( d/4) \right)
%{ 128 \times9^2 C_t^2 C_\sigma^4}\right\}.

From the expression of $V_{ij}$ in terms of $V^{K_m}_{ij}$ given by
(\ref{VG}), we
have
%
%
%eA.1 #&#
%eA.2 #&#
%eA.3 #&#
\begin{eqnarray}
\label{V-V-Km}&& P \bigl( \vert V_{ij}- \Gamma_{ij}
\vert\geq d \bigr) \nonumber\\
&&\qquad\leq P \Biggl( \sum_{m=1}^N
|a_m| \bigl\vert V^{K_m}_{ij}-
\Gamma_{ij} \bigr\vert+ \zeta\bigl( \bigl\vert
V^{K_1}_{ij}- \Gamma_{ij} \bigr\vert+ \bigl
\vert V^{K_N}_{ij}- \Gamma_{ij} \bigr\vert
\bigr) \geq d \Biggr)\nonumber
\\
&& \qquad\leq P \Biggl( \sum_{m=1}^N
|a_m| \bigl\vert V^{K_m}_{ij}-
\Gamma_{ij} \bigr\vert\geq d/ 2 \Biggr)
\nonumber
\\[-8pt]
\\[-8pt]
\nonumber
&&\qquad\quad{}+ P \bigl( \zeta\bigl
\vert V^{K_1}_{ij}- \Gamma_{ij} \bigr\vert+
\zeta\bigl\vert V^{K_N}_{ij}- \Gamma_{ij} \bigr
\vert\geq d/ 2 \bigr)
\\
&&\qquad\leq \sum_{m=1}^N P \bigl( \bigl
\vert V^{K_m}_{ij}- \Gamma_{ij} \bigr\vert\geq
d/ (2 A) \bigr) + P \bigl( \zeta\bigl\vert V^{K_1}_{ij}-
\Gamma_{ij} \bigr\vert\geq d/4 \bigr) \nonumber\\
&&\qquad\quad{}+ P \bigl( \zeta\bigl\vert
V^{K_N}_{ij}- \Gamma_{ij} \bigr\vert\geq d/4
\bigr),\nonumber
\end{eqnarray}
where $A = \sum_{m=1}^N |a_m|= 9/2 + o(1)$.

The definition of $V_{ij}^{K_m}$ in (\ref{VG1}) shows
\begin{eqnarray*}
V_{ij}^{K_m} &=&\frac{1}{K_m} \sum
_{k_m=1}^{K_m} \sum_{r=2}^{|\btau^{k_m}|}
\bigl\{X_i \bigl(\tau_r^{k_m}
\bigr) - X_i \bigl(\tau_{r-1}^{k_m} \bigr) \bigr
\} \bigl\{X_j \bigl(\tau_r^{k_m} \bigr) -
X_j \bigl(\tau_{r-1}^{k_m} \bigr) \bigr\} \\
&\equiv&
\frac{1}{K_m} \sum_{k_m=1}^{K_m}
[X_i,X_j]^{(k_m)}
\end{eqnarray*}
and
\[
 V_{ij}^{K_m} - \Gamma_{ij} =
\frac1{K_m} \sum_{k_m=1}^{K_m}
\biggl[[X_i,X_j]^{(k_m)} - \int
_0^1 \gamma_{ij}(s) \,ds \biggr].
\]
%
%where
%( \Gamma^{*K_1}_{ij} - \Gamma^{*K_N}_{ij}), \nonumber\\
% |}^{k_m}}^1 \gamma_{ij}(s) ds \right). \nonumber
With the above expression for $V^{K_m}_{ij} - \Gamma_{ij}$ we obtain
that for $d_1 >0$ and
$1 \leq m \leq N$,
%
%
%eA.4 #&#
\begin{eqnarray}
\label{V-Km} P \bigl( \bigl\vert V^{K_m}_{ij}-
\Gamma_{ij} \bigr\vert\geq d_1 \bigr) &\leq& P \Biggl(
\frac1{K_m} \sum_{k_m=1}^{K_m}
\biggl\vert[X_i,X_j]^{(k_m)} - \int
_0^1 \gamma_{ij}(s) \,ds \biggr
\vert\geq d_1 \Biggr)
\nonumber\\
&\leq& \sum_{k_m=1}^{K_m} P \biggl( \biggl
\vert[X_i,X_j]^{(k_m)}-\int
_0^1 \gamma_{ij}(s) \,ds \biggr
\vert\geq d_1 \biggr)
\\
% \leq\ & C_1 \sqrt{n} \exp\left\{- \sqrt{n} [d - C_3( \tau_1^{k_m} + 1
%-
% \tau_{|\btau^{k_m}|}^{k_m}) ]^2/C_2 \right\} \\
&\leq& C_1 K_m
\exp\biggl\{- \frac{n}{K_m} \frac{d_1^2}{C_2} \biggr\} \leq
C_3 \sqrt{n} \exp\bigl\{- \sqrt{n} d_1^2/C_4
\bigr\},\nonumber
\end{eqnarray}
where the third inequality is from Lemma~\ref{V} below and the last
inequality is due to the fact
that $\sqrt{n} \leq K_m \leq2 \sqrt{n}$ and the maximum distance
between consecutive grids in $\tau^{k_m}$
is bounded by $K_m/n \leq2 /\sqrt{n}$.

Substituting (\ref{V-Km}) into (\ref{V-V-Km})
we immediately prove Proposition~\ref{p_v} as follows:
\begin{eqnarray*}
 P \bigl( \vert V_{ij}- \Gamma_{ij} \vert\geq d
\bigr) &\leq& C_3 N \sqrt{n} \exp\bigl\{- \sqrt{n} d^2/
\bigl(4 A^2 C_4 \bigr) \bigr\}\\
&&{} + 2 C_3
\sqrt{n} \exp\bigl\{- \sqrt{n} d^2/ \bigl(16 \zeta^2
C_4 \bigr) \bigr\}
\\
& \leq& C_5 n \exp\bigl\{- \sqrt{n} d^2/C_6
\bigr\}.
\end{eqnarray*}
%
%Finally, Condition A1 implies that $|\Gamma^{*K_m}_{ij}| \leq C
%n^{-1/2}$ and
%thus $|\Gamma_{ij}^*| \leq C n^{-1/2}$. Hence,
%&& P\left( \left| V_{ij}- \Gamma_{ij} \right|\geq d \right)
% \leq P\left( \left| V_{ij}- (\Gamma_{ij}-\Gamma_{ij}^* )\right|\geq d
% - C n^{-1/2} \right) \\
%&& \leq C_1 n \exp\left\{- \sqrt{n} (d - C n^{-1/2})^2/C_2 \right\}
% \leq C_1 n \exp\left\{- \sqrt{n} d^2/C_2 \right\}.

%
%le10 #&#
\begin{lemma}\label{V}
Under model (\ref{modeldiffusion}) and condition \textup{A2}, for any sequence
$0 =\nu_0 \leq\nu_1 < \nu_2 < \cdots< \nu_m \leq\nu_{m+1}=1$ satisfying
$\max_{1\leq r \leq m+1}|\nu_r - \nu_{r-1}| \leq C/m$, we have for
$1\leq i,j \leq p$ and small $d>0$,
\begin{eqnarray*}
&& P \Biggl( \Biggl\vert\sum_{r=2}^m
\bigl(X_{i}(\nu_r) - X_{i}(\nu
_{r-1}) \bigr) \bigl(X_{j}(\nu_r) -
X_{j}(\nu_{r-1}) \bigr) - \int_0^1
\gamma_{ij}(s)\,ds \Biggr\vert\geq d \Biggr)
\\
&&\qquad \leq C_1 \exp\bigl(- m d^2/C_2 \bigr).
\end{eqnarray*}
\end{lemma}
\begin{pf}
Let $X_i^*(t)= X_i(t) - \int_0^t \mu_{is} \,ds$ and $\bX^*(t)=(X_1^*(t),
\ldots, X^*_p(t))^T$. Then $\bX^*(t)$
is a stochastic integral with respect to $\bB_t$ and has the same
quadratic variation as $\bX(t)$.
Let $\bB_t=(B_{1}(t), \ldots, B_{p}(t))^T$. With $\bsigma_t=(\sigma
_{ij}(t))$ and $\bolds{\gamma}(t)=(\gamma_{ij}(t)) = \bsigma
_t^T \bsigma_t$ we have
\[
X_i^*(t) = \int_0^t \sum
_{\ell=1}^p \sigma_{\ell i}(s)
\,dB_{\ell}(s),\qquad i =1, \ldots, p,
\]
with %each $X_i^*$ has
quadratic variation
$\langle X^*_i, X^*_i \rangle_t = \int_0^t \gamma_{ii}(s) \,ds$. Also
$X_i^* \pm X_j^*$ have quadratic variations
\[
\bigl\langle X^*_i \pm X^*_j, X^*_i
\pm X^*_j \bigr\rangle_t = \int_0^1
\bigl[\gamma_{ii}(s) + \gamma_{jj}(s) \pm2
\gamma_{ij}(s) \bigr] \,ds.
\]
Define
\[
B^*_i(t) = \int_0^t
\gamma_{ii}^{-1/2}(s) \sum_{\ell=1}^p
\sigma_{\ell
i}(s) \,dB_{\ell}(s).
\]
Then
\[
X_i^*(t) = \int_0^t
\gamma_{ii}^{1/2}(s) \,dB^*_i(s),
\]
$B^*_i$ is a continuous-time martingale and has quadratic variation
\[
\bigl\langle B^*_i, B^*_i \bigr\rangle_t
= \int_0^t \gamma_{ii}^{-1}(s)
\sum_{\ell
=1}^p \sigma^2_{\ell i}(s)
\,ds = \int_0^t \gamma_{ii}^{-1}(s)
\gamma_{ii}(s) \,ds = t,
\]
and hence L\'{e}vy's martingale characterization of Brownian motion
shows that
$B^*_i$ is a one-dimensional Brownian motion; see Karatzas and Shreve
[(\citeyear{KarShr91}), Theorem~3.16]. We can apply Lemma~3 in \citet{FanLiYu12} to each
$X_i^*$ and obtain for $1 \leq i \leq p$,
%
%
%eA.5 #&#
\begin{eqnarray}
\label{RV-extra0} &&P \Biggl( \Biggl\vert\sum_{r=2}^m
\bigl[ X^*_{i}(\nu_r) - X^*_{i}(
\nu_{r-1}) \bigr]^2 - \int_{\nu_1}^{\nu_m}
\gamma_{ii}(s) \,ds \Biggr\vert\geq d \Biggr)
\nonumber
\\[-8pt]
\\[-8pt]
\nonumber
&&\qquad\leq4 \exp\bigl\{- m
d^2/C_0 \bigr\}.
\end{eqnarray}
Similarly for $X_i^* \pm X_j^*$, we define
\[
B^{\pm}_{ij}(s) = \int_0^t
\bigl[\gamma_{ii}(s) + \gamma_{jj}(s) \pm2 \gamma
_{ij}(s) \bigr]^{-1/2} \sum_{\ell=1}^p
\bigl[\sigma_{\ell i}(s) \pm\sigma_{\ell j}(s) \bigr]
\,dB_{\ell}(s).
\]
Then
\[
X_i^*(t) \pm X_j^*(t) = \int_0^t
\bigl[\gamma_{ii}(s) + \gamma_{jj}(s) \pm2
\gamma_{ij}(s) \bigr]^{1/2} \,dB^{\pm}_{ij}(s),
\]
$B^{\pm}_{ij}$ are continuous-time martingales with quadratic variations
\begin{eqnarray*}
 \bigl\langle B^{\pm}_{ij}, B^{\pm}_{ij}
\bigr\rangle_t &= &\int_0^t
\bigl[\gamma_{ii}(s) + \gamma_{jj}(s) \pm2
\gamma_{ij}(s) \bigr]^{-1}\\
&&\quad{}\times \sum_{\ell=1}^p
\bigl[ \sigma^2_{\ell i}(s) + \sigma^2_{\ell j}(s)
\pm2 \sigma_{\ell i}(s) \sigma_{\ell j}(s) \bigr] \,ds
\\
& =& \int_0^t \bigl[\gamma_{ii}(s)
+ \gamma_{jj}(s) \pm2 \gamma_{ij}(s) \bigr]^{-1}
\bigl[\gamma_{ii}(s) + \gamma_{jj}(s) \pm2
\gamma_{ij}(s) \bigr] \,ds = t,
\end{eqnarray*}
and hence L\'evy's martingale characterization of Brownian motion
implies that $B^{\pm}_{ij}$ are one-dimensional Brownian motions.
We can apply Lemma~3 in \citet{FanLiYu12} to each of $X_i^* + X_j^*$
and $X_i^* - X_j^*$ and obtain for $1 \leq i, j \leq p$,
%
%eA.6 #&#
%eA.7 #&#
\begin{eqnarray}\label{RV-extra1}
&&P \Biggl( \Biggl\vert\sum_{r=2}^m
\bigl( \bigl[X^*_{i}(\nu_r) - X^*_{i}(
\nu_{r-1}) \bigr] \pm\bigl[X^*_{j}(\nu_r) -
X^*_{j}(\nu_{r-1}) \bigr] \bigr)^2 \nonumber\\
&&\hspace*{84pt}{}- \int
_{\nu_1}^{\nu_m} \bigl[\gamma_{ii}(s) +
\gamma_{jj}(s) \pm2 \gamma_{ij}(s) \bigr] \,ds \Biggr
\vert\geq d \Biggr)
\\
&&\qquad \leq4 \exp\bigl\{- m d^2/C_0 \bigr
\}.\nonumber
\end{eqnarray}
Note that
\begin{eqnarray*}
&& 4 \gamma_{ij}(s) = \bigl[\gamma_{ii}(s) +
\gamma_{jj}(s) + 2 \gamma_{ij}(s) \bigr] - \bigl[
\gamma_{ii}(s) + \gamma_{jj}(s) - 2 \gamma_{ij}(s)
\bigr],
\\
&& 4 \sum_{r=2}^m \bigl(X^*_{i}(
\nu_r) - X^*_{i}(\nu_{r-1}) \bigr)
\bigl(X^*_{j}(\nu_r) - X^*_{j}(
\nu_{r-1}) \bigr)
\\
&&\qquad = \sum_{r=2}^m \bigl\{
\bigl[X^*_{i}(\nu_r) - X^*_{i}(
\nu_{r-1}) \bigr] + \bigl[X^*_{j}(\nu_r) -
X^*_{j}(\nu_{r-1}) \bigr] \bigr\}^2
\nonumber
\\
&&\qquad\quad{} - \sum_{r=2}^m \bigl\{
\bigl[X^*_{i}(\nu_r) - X^*_{i}(
\nu_{r-1}) \bigr] - \bigl[X^*_{j}(\nu_r) -
X^*_{j}(\nu_{r-1}) \bigr] \bigr\}^2,
\end{eqnarray*}
and thus
\begin{eqnarray*}
\label{RV-extra2}
&& 4 \Biggl\vert\sum_{r=2}^m
\bigl(X^*_{i}(\nu_r) - X^*_{i}(\nu
_{r-1}) \bigr) \bigl(X^*_{j}(\nu_r) -
X^*_{j}(\nu_{r-1}) \bigr) - \int_{\nu_1}^{\nu_m}
\gamma_{ij}(s)\,ds \Biggr\vert
\\[-2pt]
&&\qquad \leq\Biggl\vert\sum_{r=2}^m
\bigl\{ \bigl[X^*_{i}(\nu_r) - X^*_{i}(\nu
_{r-1}) \bigr] + \bigl[X^*_{j}(\nu_r) -
X^*_{j}(\nu_{r-1}) \bigr] \bigr\}^2\\[-2pt]
&&\hspace*{69pt}\qquad\quad{} - \int
_{\nu_1}^{\nu_m} \bigl[\gamma_{ii}(s) +
\gamma_{jj}(s) + 2 \gamma_{ij}(s) \bigr] \,ds \Biggr\vert
\\[-2pt]
&&\qquad\quad{} + \Biggl\vert\sum_{r=2}^m \bigl\{
\bigl[X^*_{i}(\nu_r) - X^*_{i}(
\nu_{r-1}) \bigr] - \bigl[X^*_{j}(\nu_r) -
X^*_{j}(\nu_{r-1}) \bigr] \bigr\}^2\\[-2pt]
&&\hspace*{81pt}\qquad\quad{} - \int
_{\nu_1}^{\nu_m} \bigl[\gamma_{ii}(s) +
\gamma_{jj}(s) - 2 \gamma_{ij}(s) \bigr] \,ds \Biggr
\vert.
\end{eqnarray*}
Combining (\ref{RV-extra1}) and above inequality we conclude
%
%
%eA.8 #&#
%eA.9 #&#
\begin{eqnarray}
\label{RV-extra3} && P \Biggl( \Biggl\vert\sum_{r=2}^m
\bigl(X^*_{i}(\nu_r) - X^*_{i}(\nu
_{r-1}) \bigr) \bigl(X^*_{j}(\nu_r) -
X^*_{j}(\nu_{r-1}) \bigr) - \int_{\nu_1}^{\nu_m}
\gamma_{ij}(s)\,ds \Biggr\vert\geq d \Biggr)
\nonumber\hspace*{-35pt}
\\[-9pt]
\\[-9pt]
\nonumber
&&\qquad \leq8 \exp\bigl\{- m (d/8)^2/C_0 \bigr\} = 8 \exp
\bigl\{- m d^2/(64 C_0) \bigr\}.\hspace*{-35pt}
\end{eqnarray}
On the other hand,
%
%
%eA.10 #&#
\begin{eqnarray}
\label{RV-extra4} && \sum_{r=2}^m
\bigl(X_{i}(\nu_r) - X_{i}(
\nu_{r-1}) \bigr) \bigl(X_{j}(\nu_r) -
X_{j}(\nu_{r-1}) \bigr)
\nonumber
\\[-2pt]
&&\qquad = \sum_{r=2}^m \biggl\{
\bigl[X^*_{i}(\nu_r) - X^*_{i}(
\nu_{r-1}) \bigr] + \int_{\nu_{r-1}}^{\nu_r}
\mu_{is} \,ds \biggr\}\nonumber\\[-2pt]
&&\hspace*{47pt}{}\times \biggl\{ \bigl[X^*_{j}(
\nu_r) - X^*_{j}(\nu_{r-1}) \bigr] + \int
_{\nu_{r-1}}^{\nu
_r} \mu_{js} \,ds \biggr\}
\nonumber
\\[-2pt]
&&\qquad = \sum_{r=2}^m \bigl(X^*_{i}(
\nu_r) - X^*_{i}(\nu_{r-1}) \bigr)
\bigl(X^*_{j}(\nu_r) - X^*_{j}(
\nu_{r-1}) \bigr) \\[-2pt]
&&\qquad\quad{}+ \sum_{r=2}^m
\int_{\nu_{r-1}}^{\nu_r} \mu_{is} \,ds \int
_{\nu_{r-1}}^{\nu_r} \mu_{js} \,ds
\nonumber\\[-2pt]
&&\qquad\quad{} + \sum_{r=2}^m \bigl[X^*_{i}(
\nu_r) - X^*_{i}(\nu_{r-1}) \bigr] \int
_{\nu
_{r-1}}^{\nu_r} \mu_{js} \,ds\nonumber\\[-2pt]
&&\qquad\quad{} + \sum
_{r=2}^m \bigl[X^*_{j}(
\nu_r) - X^*_{j}(\nu_{r-1}) \bigr] \int
_{\nu_{r-1}}^{\nu_r} \mu_{is} \,ds.\nonumber
\end{eqnarray}
From condition A2 we have that $\mu_i$ and $\mu_j$ are bounded by
$c_1$, and thus
%
%
%eA.11 #&#
\begin{equation}
 \Biggl\vert\sum_{r=2}^m \int
_{\nu_{r-1}}^{\nu_r} \mu_{is} \,ds \int
_{\nu
_{r-1}}^{\nu_r} \mu_{js} \,ds \Biggr\vert
\leq\frac{c^2_1}{m}. \label{RV-extra5}
\end{equation}
Applications of H\"older's inequality lead to
%
%
%eA.12 #&#
%eA.13 #&#
\begin{eqnarray}\label{RV-extra6}
 &&\Biggl\vert\sum_{r=2}^m
\bigl[X^*_{i}(\nu_r) - X^*_{i}(
\nu_{r-1}) \bigr] \int_{\nu
_{r-1}}^{\nu_r}
\mu_{js} \,ds \Biggr\vert^2 \nonumber\\
&&\qquad\leq\sum
_{r=2}^m \bigl[X^*_{i}(
\nu_r) - X^*_{i}(\nu_{r-1}) \bigr]^2
\sum_{r=2}^m \biggl\vert\int
_{\nu_{r-1}}^{\nu_r} \mu_{js} \,ds \biggr\vert
^2
\\
&&\qquad \leq\frac{c^2_1}{m} \sum_{r=2}^m
\bigl[X^*_{i}(\nu_r) - X^*_{i}(\nu
_{r-1}) \bigr]^2,
\nonumber\\
\label{RV-extra7}&& \Biggl\vert\sum_{r=2}^m
\bigl[X^*_{j}(\nu_r) - X^*_{j}(
\nu_{r-1}) \bigr] \int_{\nu
_{r-1}}^{\nu_r}
\mu_{is} \,ds \Biggr\vert^2
\nonumber
\\[-8pt]
\\[-8pt]
\nonumber
&&\qquad\leq\frac{c^2_1}{m} \sum
_{r=2}^m \bigl[X^*_{j}(
\nu_r) - X^*_{j}(\nu_{r-1})
\bigr]^2.
\end{eqnarray}
%
%where we use Condition A2 that $\mu_i$ and $\mu_j$ are bounded by
%$c_1$.
%Again from Condition A2 we have
%&& \sum_{r=2}^m [X^*_{i}(\nu_r) - X^*_{i}(\nu_{r-1})]^2 \rightarrow
%&& \sum_{r=2}^m [X^*_{j}(\nu_r) - X^*_{j}(\nu_{r-1})]^2 \rightarrow
From (\ref{RV-extra4}) we have
%that the last three terms on the right-hand side of (\ref{RV-extra4})
%are of order $m^{-1/2}$, thus
%there exists a constant $C_3$ such that with probability tending to
%one,
%
%
%eA.14 #&#
\begin{eqnarray}
\label{RV-extra10} && P \Biggl( \Biggl\vert\sum_{r=2}^m
\bigl(X_{i}(\nu_r) - X_{i}(\nu
_{r-1}) \bigr) \bigl(X_{j}(\nu_r) -
X_{j}(\nu_{r-1}) \bigr) - \int_{\nu_1}^{\nu_m}
\gamma_{ij}(s)\,ds \Biggr\vert\geq d \Biggr)
\nonumber
\\
&&\qquad \leq P \Biggl( \Biggl\vert\sum_{r=2}^m
\bigl(X^*_{i}(\nu_r) - X^*_{i}(\nu
_{r-1}) \bigr) \bigl(X^*_{j}(\nu_r) -
X^*_{j}(\nu_{r-1}) \bigr)\nonumber\\
&&\hspace*{158pt}\qquad{} - \int_{\nu_1}^{\nu_m}
\gamma_{ij}(s)\,ds \Biggr\vert\geq d/4 \Biggr)
\nonumber
\\
&&\qquad\quad{} + P \Biggl( \Biggl\vert\sum_{r=2}^m
\int_{\nu_{r-1}}^{\nu_r} \mu_{is} \,ds \int
_{\nu_{r-1}}^{\nu_r} \mu_{js} \,ds \Biggr\vert
\geq d/4 \Biggr) % \nonumber
\nonumber\\
&&\qquad\quad{}+ P \Biggl( \Biggl\vert\sum
_{r=2}^m \bigl[X^*_{i}(
\nu_r) - X^*_{i}(\nu_{r-1}) \bigr] \int
_{\nu_{r-1}}^{\nu_r} \mu_{js} \,ds \Biggr\vert
\geq d/4 \Biggr)
\\
% + P\left( \left| \frac{c^2_1}{m} \sum_{r=2}^m [X^*_{i}(\nu_r) -
%X^*_{i}(\nu_{r-1})]^2 \right|^{1/2} \geq d/4 \right) \nonumber\\
&&\qquad\quad{} + P \Biggl( \Biggl\vert\sum
_{r=2}^m \bigl[X^*_{j}(
\nu_r) - X^*_{j}(\nu_{r-1}) \bigr] \int
_{\nu_{r-1}}^{\nu_r} \mu_{is} \,ds \Biggr\vert
\geq d/4 \Biggr)
\nonumber\\
&&\qquad \leq8 \exp\bigl\{- m (d/4)^2/(64 C_0) \bigr\} + 1
\biggl( \frac
{c^2_1}{m} \geq d/4 \biggr)\nonumber\\
&&\qquad\quad{} + P \Biggl( \sum
_{r=2}^m \bigl[X^*_{i}(
\nu_r) - X^*_{i}(\nu_{r-1}) \bigr]^2
\geq m d^2/ \bigl(16 c_1^2 \bigr) \Biggr)
\nonumber\\
&&\qquad\quad{} + P \Biggl( \sum_{r=2}^m
\bigl[X^*_{i}(\nu_r) - X^*_{i}(
\nu_{r-1}) \bigr]^2 \geq m d^2/ \bigl(16
c_1^2 \bigr) \Biggr),\nonumber
\end{eqnarray}
where the last inequality is due to the bounds obtained from (\ref
{RV-extra3}) and (\ref{RV-extra5})--(\ref{RV-extra7}) for the four
respective probability terms.
We handle the last two terms on the right-hand side of (\ref
{RV-extra10}) as follows.
%First if $d > 4 c_1^2/m$, the probability in the second term
% P\left( \left| \sum_{r=2}^m \int_{\nu_{r-1}}^{\nu_r} \mu_{is} \,ds
% P\left( \frac{c^2_1}{m} \geq d/4 \right) = 0;
If $m d^2/(16 c_1^2) - c_2 > 0$ [or equivalently $d > 4 c_1 (c_2/m)^{1/2}$],
using condition A2 (which implies $\gamma_{ii} \leq c_2$ and $\gamma
_{jj} \leq c_2$) and~(\ref{RV-extra0}), we get
%
%
%eA.15 #&#
\begin{eqnarray}
\label{RV-extra8} && P \Biggl( \sum_{r=2}^m
\bigl[X^*_{i}(\nu_r) - X^*_{i}(
\nu_{r-1}) \bigr]^2 \geq m d^2/ \bigl(16
c_1^2 \bigr) \Biggr)\nonumber\\
&&\qquad\quad{} + P \Biggl( \sum
_{r=2}^m \bigl[X^*_{j}(
\nu_r) - X^*_{j}(\nu_{r-1}) \bigr]^2
\geq m d^2/ \bigl(16 c_1^2 \bigr) \Biggr)
\nonumber
\\
&&\qquad \leq P \Biggl( \sum_{r=2}^m
\bigl[X^*_{i}(\nu_r) - X^*_{i}(
\nu_{r-1}) \bigr]^2 - \int_{\nu_1}^{\nu_m}
\gamma_{ii}(s) \,ds \geq m d^2/ \bigl(16
c_1^2 \bigr) - c_2 \Biggr)
\\
&&\qquad\quad{} + P \Biggl( \sum_{r=2}^m
\bigl[X^*_{j}(\nu_r) - X^*_{j}(
\nu_{r-1}) \bigr]^2 - \int_{\nu_1}^{\nu_m}
\gamma_{jj}(s) \,ds \geq m d^2/ \bigl(16
c_1^2 \bigr) - c_2 \Biggr)
\nonumber\\
&&\qquad \leq P \Biggl( \Biggl\vert\sum_{r=2}^m
\bigl[X^*_{i}(\nu_r) - X^*_{i}(\nu
_{r-1}) \bigr]^2 - \int_{\nu_1}^{\nu_m}
\gamma_{ii}(s) \,ds \Biggr\vert\geq m d^2/ \bigl(16
c_1^2 \bigr) - c_2 \Biggr)
\nonumber
\\
&&\qquad\quad{} + P \Biggl( \Biggl\vert\sum_{r=2}^m
\bigl[X^*_{j}(\nu_r) - X^*_{j}(
\nu_{r-1}) \bigr]^2 - \int_{\nu_1}^{\nu_m}
\gamma_{jj}(s) \,ds \Biggr\vert\geq m d^2/ \bigl(16
c_1^2 \bigr) - c_2 \Biggr)
\nonumber
\\
&&\qquad \leq8 \exp\bigl\{- m \bigl[m d^2/ \bigl(16 c_1^2
\bigr) - c_2 \bigr]^2/C_0 \bigr\},\nonumber
\end{eqnarray}
which is bounded by $8 \exp\{- m d^2/C_0 \}$, if $m [m
d^2/(16 c_1^2) - c_2]^2 > m d^2$, which is true provided that
%
%
%eA.16 #&#
\begin{equation}
\label{RV-extra9} %&& 8 \exp\left\{- m d^2/C_0 \right\}, \nonumber\\
 %\mbox{ if } m [m d^2/(16 c_1^2) - c_2]^2 > m d^2 \mbox{ or
%equivalently}
d > \frac{8 c_1^2}{m}
+ \frac{4 c_1}{m} \bigl( 4 c_1^2 + m c_2
\bigr)^{1/2}.
\end{equation}
%
%third, the third and fourth terms on the right-hand side of (
%the above probability bound for the third term applies to the fourth
%term on right-hand side of (\ref{RV-extra10}).
Putting together (\ref{RV-extra10}) and the probability bound from
(\ref
{RV-extra8})--(\ref{RV-extra9}), %for the last two terms on the right
%hand side of (\ref{RV-extra10})
we conclude that if
%
%eA.17 #&#
\begin{eqnarray}
&&d > \max\biggl\{ \frac{4 c_1^2}{m}, \frac{4 c_1 c_2^{1/2}
}{m^{1/2}}, \frac{8 c_1^2}{m} +
\frac{4 c_1}{m} \bigl( 4 c_1^2 + m c_2
\bigr)^{1/2} \biggr\}\nonumber\\
&&\hphantom{d }= \frac{8 c_1^2}{m} + \frac{4 c_1}{m} \bigl( 4
c_1^2 + m c_2 \bigr)^{1/2},\nonumber
\\
\label{RV-1} && P \Biggl( \Biggl\vert\sum_{r=2}^m
\bigl(X_{i}(\nu_r) - X_{i}(\nu
_{r-1}) \bigr) \bigl(X_{j}(\nu_r) -
X_{j}(\nu_{r-1}) \bigr) \nonumber\\
&&\hspace*{145pt}{}- \int_{\nu_1}^{\nu_m}
\gamma_{ij}(s)\,ds \Biggr\vert\geq d \Biggr)
\nonumber
\\[-8pt]
\\[-8pt]
\nonumber
&& \qquad\leq8 \exp\bigl\{- m d^2/ (1024 C_0) \bigr\} + 8
\exp\bigl\{- m d^2/ C_0 \bigr\}
\\
&&\qquad \leq16 \exp\bigl\{- m d^2/ (1024 C_0) \bigr\}.\nonumber
\end{eqnarray}
From condition A2 we have $|\gamma_{ij}| \leq(\gamma_{ii} \gamma
_{jj})^{1/2} \leq c_2$ and
\begin{eqnarray*}
&&\biggl\vert\int_{\nu_1}^{\nu_m}
\gamma_{ij}(s)\,ds - \int_0^1
\gamma_{ij}(s)\,ds \biggr\vert\\
&&\qquad\leq c_2 (
\nu_1 + 1 - \nu_m) \leq2 c_2 /m.
\end{eqnarray*}
Then (\ref{RV-1}) and above inequality imply that if
%
%
%eA.18 #&#
\begin{equation}
\label{RV-extra12}  d > \max\biggl\{ \frac{4 c_2}{m}, \frac{8
c_1^2}{m} +
\frac{4 c_1}{m} \bigl( 4 c_1^2 + m c_2
\bigr)^{1/2} \biggr\},
\end{equation}
\begin{eqnarray*}
&& P \Biggl( \Biggl\vert\sum_{r=2}^m
\bigl(X_{i}(\nu_r) - X_{i}(\nu
_{r-1}) \bigr) \bigl(X_{j}(\nu_r) -
X_{j}(\nu_{r-1}) \bigr) - \int_0^1
\gamma_{ij}(s)\,ds \Biggr\vert\geq d \Biggr)
\\
&&\qquad \leq P \Biggl( \Biggl\vert\sum_{r=2}^m
\bigl(X_{i}(\nu_r) - X_{i}(\nu
_{r-1}) \bigr) \bigl(X_{j}(\nu_r) -
X_{j}(\nu_{r-1}) \bigr) \nonumber\\
&&\hspace*{155pt}\qquad{}- \int_{\nu_1}^{\nu_m}
\gamma_{ij}(s)\,ds \Biggr\vert\geq d/2 \Biggr)
\\
&&\qquad \leq16 \exp\bigl\{- m (d/2)^2/(1024 C_0) \bigr\}
= 16 \exp\bigl\{ - m d^2/ (4096 C_0) \bigr\}.
\end{eqnarray*}
This proves the lemma with $C_1=16$ and $C_2=4096 C_0$ for $d$
satisfies (\ref{RV-extra12}).

If (\ref{RV-extra12}) is not satisfied, we have
\begin{eqnarray*}
&& d \leq\max\biggl\{ \frac{4 c_2}{m}, \frac{8 c_1^2}{m} +
\frac{4
c_1}{m} \bigl( 4 c_1^2 + m c_2
\bigr)^{1/2} \biggr\} \leq\frac{8 c_1^2 + 4 c_2 + 4 c_1
c_2^{1/2}}{m^{1/2}} \equiv\frac{ C }{m^{1/2}}.
\end{eqnarray*}
%
%where $C = 8 c_1^2 + 4 c_2 + 4 c_1 c_2^{1/2}$. %$C$ is some constant
%depending on $c_1$ and $c_2$.
Then the tail probability bound in the lemma obeys
\[
C_1 \exp\bigl\{- m d^2/C_2 \bigr\} \geq
C_1 \exp\bigl\{- C^2/C_2 \bigr\},
\]
and we easily show the probability inequality in the lemma by choosing
$C_1=C_1^\prime$ and $C_2=C_2^\prime$, where
$C_1^\prime$ and $C_2^\prime$ satisfy $C_1^\prime\exp\{-
C^2/C_2^\prime\} \geq1$.

Finally taking $C_1 = \max(16, C_1^\prime)$ and $C_2=\max(4096 C_0,
C_2^\prime)$ we establish the tail probability,
regardless whether $d$ satisfies (\ref{RV-extra12}) or not, and
complete the proof.
%&& \sum_{r=2}^m \left([X_{i}(\nu_r) - X_{i}(\nu_{r-1})] \pm[X_{j}(
%&& \leq2 \sum_{r=2}^m ([X^*_{i}(\nu_r) - X^*_{i}(\nu_{r-1})] \pm
%[X^*_{j}(\nu_r) - X^*_{j}(\nu_{r-1})])^2 + 2 \sum_{r=2}^m \left| \int_{
%&& \int_{\nu_1}^{\nu_m} [\gamma_{ii}(s) + \gamma_{jj}(s) \pm2
%Except for the correction at the two end points $\nu_1$
%and $\nu_m$ the lemma is the same as Lemma~3 in Fan et al. (2011).
%a direct application of Lemma~3 in Fan et al. (2011) to $X_i + X_j$
% and $X_i - X_j$.
%The same arguments in the proof of Lemma~3 in Fan et al. (2011) lead to
%for $d > C /m$,
%&& P\left( \left| \sum_{r=2}^m (X_{i}(\nu_r) - X_{i}(\nu_{r-1}))(X_{j}(
%X_{j}(\nu_{r-1})) - \int_0^1 \gamma_{ij}(s)ds \right|\geq d \right) \\
%&& \leq P\left( \left| \sum_{r=2}^m (X_{i}(\nu_r) - X_{i}(
%X_{j}(\nu_{r-1})) - \int_{\nu_1}^{\nu_m} \gamma_{ij}(s)ds \right|\geq
%d - C/m \right) \\
%&& \leq C_1 \exp\left\{- m (d - C/m)^2/C_2 \right\} \leq C_1 \exp\left
%where the second inequality is from (\ref{RV-1}) and the last
%inequality is due to the fact that
%$- m (d - C/m)^2 = - m d^2 + 2 C - C/m \leq- m d^2 + 2 C$. This
%proves the lemma for $d > C /m$.
\end{pf}
%
%sA.2 #&#
\subsection{\texorpdfstring{Proof of Proposition \protect\ref{pG2}}{Proof of Proposition 6}}
\label{sectionG2}

%Under the assumptions of Theorem~\ref{UnivariateConv}, we have for
%$1\leq i,j \leq p$,
%P\left( \left|G_{ij}(2)\right|\geq d/4 \right) & \leq& C n^{1/2} C_
% \exp\left\{-\frac{\sqrt{n} d^2 }{\max\{ 648\eta_jC_\sigma^2,32 C_t^2
%C_\sigma^4 \}} \right\},\\
%P\left( \left|G_{ij}(3)\right|\geq d/4 \right) & \leq& C n^{1/2}C_
% \exp\left\{-\frac{\sqrt{n} d^2 }{\max\{ 648\eta_i C_\sigma^2,32C_t^2
%C_\sigma^4 \}} \right\}.

As the proofs for $G_{ij}(2)$ and $G_{ij}(3)$ are similar, we give arguments
only for $G_{ij}(2)$. Lemma~\ref{lemmaG2} below establishes the tail
probability for
$G^{K_m}_{ij}(2)$. Using the expression of $G_{ij}(2)$ in terms of
$G^{K_m}_{ij}(2)$
given by (\ref{VG}) and applying Lemma~\ref{lemmaG2}, we obtain
\begin{eqnarray*}
&& P \bigl( \bigl\vert G_{ij}(2) \bigr\vert\geq d \bigr) \\
&&\qquad\leq\sum
_{m=1}^N P \bigl( \bigl\vert
G^{K_m}_{ij}(2) \bigr\vert\geq d/(2 A) \bigr) + P \bigl(
\zeta\bigl\vert G^{K_1}_{ij}(2) \bigr\vert\geq d/4
\bigr) \\
&&\qquad\quad{}+ P \bigl( \zeta\bigl\vert G^{K_N}_{ij}(2) \bigr
\vert\geq d/4 \bigr)
\\
&&\qquad\leq C_1 N \sqrt{n} \exp\bigl\{-\sqrt{n} d^2 /
\bigl(4 A^2 C_2 \bigr) \bigr\} + 2 C_1
\sqrt{n} \exp\bigl\{-\sqrt{n} d^2 / \bigl(16 \zeta^2
C_2 \bigr) \bigr\}
\\
&&\qquad\leq C_3 n \exp\bigl\{-\sqrt{n} d^2
/C_4 \bigr\},
\end{eqnarray*}
where $A=\sum_{m=1}^N |a_m| = 9/2 + o(1)$.

%-------lemma-------------------------------
%
%
%le11 #&#
\begin{lemma}\label{lemmaG2}
Under the assumptions of Theorem~\ref{UnivariateConv}, we have for
$1\leq i,j \leq p$ and $ 1 \leq m \leq N$,
\[
P \bigl(\bigl|G^{K_m}_{ij}(2)\bigr|\geq d \bigr) \leq C_1
\sqrt{n} \exp\bigl\{-\sqrt{n} d^2/C_2 \bigr\}.
% C \frac{C_\sigma} {\sqrt{K_m}} \exp\left\{-\frac{K_m d^2}{648C_
% \eta_j}\right\} + C C_\sigma^2 \exp\left\{-\frac{\sqrt{n} d^2}{32
%C_t^2
% C_\sigma^4} \right\}
\]
\end{lemma}

\begin{pf}
Simple algebraic manipulations show
\begin{eqnarray*}
G_{ij}^{K_m}(2)& = & \frac1 {K_m} \sum
_{k_m=1}^{K_m} \sum_{r=2}^{|\btau
^{k_m}|}
\bigl[ X_i \bigl(\tau_{r}^{k_m} \bigr) -
X_i \bigl(\tau_{r-1}^{k_m} \bigr) \bigr] \bigl[
\varepsilon_j \bigl(\tau_{r}^{k_m} \bigr) -
\varepsilon_j \bigl(\tau_{r-1}^{k_m} \bigr)
\bigr]
\\
&= & \frac1 {K_m} \sum_{k_m=1}^{K_m}
\sum_{r=2}^{|\btau^{k_m}|} \bigl[ X_i
\bigl(\tau_{r}^{k_m} \bigr) - X_i \bigl(
\tau_{r-1}^{k_m} \bigr) \bigr] \varepsilon_j
\bigl( \tau_{r}^{k_m} \bigr)
\\
&&{} - \frac1{K_m} \sum_{k_m=1}^{K_m}
\sum_{r=2}^{|\btau^{k_m}|} \bigl[ X_i
\bigl(\tau_{r}^{k_m} \bigr) - X_i \bigl(
\tau_{r-1}^{k_m} \bigr) \bigr] \varepsilon_j
\bigl( \tau_{r-1}^{k_m} \bigr)
\\
&\equiv& R_5^{K_m} - R_6^{K_m}.
\end{eqnarray*}
The lemma is proved if we establish tail probabilities for both
$R_5^{K_m}$ and $R_6^{K_m}$.
Due to similarity, we give the arguments only for $R_5^{K_m}$. Since
$\bX_t$ and $\varepsilon_i(t_\ell)$
are independent, conditional on the\vadjust{\goodbreak} whole path of $\bX_t$, $R_5^{K_m}$
is the weighted sum of
$\varepsilon_j(\cdot)$.
%normal random variables with mean 0 and variance $\eta_j$.
Hence,
%
%
%eA.19 #&#
%eA.20 #&#
%eA.21 #&#
%eA.22 #&#
\begin{eqnarray}
\label{equ-R-5}  &&P \bigl(\bigl|R_5^{K_m}\bigr|\geq d \bigr)\nonumber\\
&&\qquad= E
\bigl[ P \bigl(\bigl|R_5^{K_m}\bigr|\geq d |\bX_t, t \in
[0,1] \bigr) \bigr]
\nonumber\hspace*{-15pt}\\
&&\qquad=  E \Biggl[ P \Biggl( \Biggl| \sum_{k_m=1}^{K_m}
\sum_{r=2}^{|\btau^{k_m}|} \bigl[ X_i
\bigl(\tau_{r}^{k_m} \bigr) - X_i \bigl(
\tau_{r-1}^{k_m} \bigr) \bigr] \varepsilon_j
\bigl( \tau_{r}^{k_m} \bigr) \Biggr| \geq d K_m \vert
\bX_t, t \in[0,1] \Biggr) \Biggr]
\nonumber\hspace*{-15pt}
\\[-8pt]
\\[-8pt]
\nonumber
&&\qquad\leq E \biggl[ c_0 \exp\biggl\{ - \frac{d^2 K_m}{2 \tau_0
V_{ii}^{K_m} \eta_j} \biggr\}
\biggr]\hspace*{-15pt}
\\
&&\qquad=  E \biggl[ c_0 \exp\biggl\{ - \frac{d^2 K_m}{2 \tau_0
V_{ii}^{K_m}\eta_j} \biggr\}1(
\Omega_0) \biggr] + E \biggl[ c_0 \exp\biggl\{ -
\frac{d^2 K_m }{2 \tau_0 V_{ii}^{K_m} \eta
_j} \biggr\}1 \bigl(\Omega^c_0 \bigr)
\biggr]\hspace*{-15pt} % \\ \equiv&\
\nonumber\\
&&\qquad\equiv R_{5,1}^{K_m} +
R_{5,2}^{K_m},\nonumber\hspace*{-15pt}
\end{eqnarray}
where the inequality is due to the subGaussianity of $\varepsilon
_j(\cdot)$ defined in (\ref{subGaussian-2}),
$\eta_j$ is the variance of $\varepsilon_j(\cdot)$, $V_{ii}^{K_m}$ is
given by (\ref{VG1}) with an expression
\begin{eqnarray*}
V_{ii}^{K_m} &=& \frac{1}{K_m} \sum
_{k_m=1}^{K_m}[X_i,X_i]^{k_m}
= \frac{1}{
K_m} \sum_{k_m=1}^{K_m} \sum
_{r=2}^{|\btau^{k_m}|} \bigl[ X_i \bigl(
\tau_{r}^{k_m} \bigr) - X_i \bigl(
\tau_{r-1}^{k_m} \bigr) \bigr] ^2
\end{eqnarray*}
and
\[
\Omega_0  =  \bigl\{ \bigl|V_{ii}^{K_m} -
\Gamma_{ii} \bigr| \geq d \bigr\}.
\]
%
%$V_{ii}^{K_m}$ is the same definition in subsection~\ref{sectionV}.
From the definition of $\Omega_0$ and conditions A1--A2, we have $\eta
_j \leq\kappa^2$, $\Gamma_{ii} \leq c_2$ and
$V_{ii}^{K_m} \leq\Gamma_{ii} + d \leq c_2 + d$ on $\Omega_0^c$. Thus
for small $d$ we have
%
%
%eA.23 #&#
\begin{eqnarray}
\label{equ-R-5-2} R_{5,2}^{K_m} &=& E \biggl[ c_0
\exp\biggl\{ - \frac{K_m d^2}{2 \tau_0 V_{ii}^{K_m}\eta_j} \biggr
\} 1 \bigl(\Omega^c_0
\bigr) \biggr]
\nonumber
\\[-8pt]
\\[-8pt]
\nonumber
&\leq& C_1 \exp\bigl\{- K_m
d^2/C_2 \bigr\} \leq C_1 \exp\bigl\{- \sqrt
{n} d^2/C_2 \bigr\}.
\end{eqnarray}
On the other hand, from (\ref{V-Km}) %Lemma~\ref{V} and
(in the proof of Proposition~\ref{p_v}) we have
\[
P( \Omega_0 ) \leq C_3 \sqrt{n} \exp\bigl\{-\sqrt{n}
d^2/C_4 \bigr\},
\]
and thus
%
%
%eA.24 #&#
\begin{eqnarray}
\label{equ-R-5-1} R_{5,1}^{K_m}& =& E \biggl[ c_0
\exp\biggl\{ - \frac{d^2 K_m}{2 \tau_0 V_{ii}^{K_m}\eta_j} \biggr
\} 1(\Omega_0) \biggr]
\leq c_0 P(\Omega_0)
\nonumber
\\[-8pt]
\\[-8pt]
\nonumber
& \leq& c_0
C_3 \sqrt{n} \exp\bigl\{- \sqrt{n} d^2/C_4
\bigr\}.
\end{eqnarray}
Finally substituting (\ref{equ-R-5-2}) and (\ref{equ-R-5-1}) into
(\ref
{equ-R-5}) we obtain
\begin{eqnarray*}
 P \bigl(\bigl|R_5^{K_m}\bigr|\geq d \bigr) &\leq& C_1
\exp\bigl\{-\sqrt{n} d^2/C_2 \bigr\} + c_0
C_3 \sqrt{n} \exp\bigl\{-\sqrt{n} d^2/C_4
\bigr\} \\
&\leq &C_5 \sqrt{n} \exp\bigl\{-\sqrt{n}
d^2/C_6 \bigr\}.
\end{eqnarray*}
\upqed\end{pf}
%
%sA.3 #&#
\subsection{\texorpdfstring{Proof of Proposition \protect\ref{p_g1}}{Proof of Proposition 7}}
\label{sectionG1}

%Under the assumptions of Theorem~\ref{UnivariateConv}, we have for
%$1\leq i,j \leq p$,
%P \left(\left| G_{ij}(1) \right|\geq d/4 \right) \leq C_{G1}
Denote by $\rho_{ij}(0)$ the correlation between $\varepsilon_i(t_1)$
and $\varepsilon_j(t_1)$.
From the expression of $G_{ij}(1)$ in terms of $G^{K_m}_{ij}(1)$ given
by (\ref{VG}) we obtain
%the definition of $\bG(1)=(G_{ij}(1))$ in (\ref{VG}),
that $P (\vert G_{ij}(1) \vert\geq d )$ is bounded by
\begin{eqnarray*}
%P \left(\left| G_{ij}(1) \right|\geq d \right) &\leq&
&& P \Biggl( \Biggl\vert\sum_{m=1}^N
a_m G_{ij}^{K_m}(1) + 2 \sqrt{
\eta_i \eta_j}\rho_{ij}(0) \Biggr\vert
\geq d/2 \Biggr)\\
&&\quad{} + P \bigl( \bigl\vert\zeta\bigl(G_{ij}^{K_1}(1)
- G_{ij}^{K_N}(1) \bigr) - 2 \sqrt{\eta_i
\eta_j}\rho_{ij}(0) \bigr\vert\geq d/2 \bigr)
\\
&&\qquad \leq C_1 \sqrt{n} \exp\bigl\{- \sqrt{n} d^2 /(4
C_2) \bigr\} + C_3 \exp\bigl\{- n d^2/( 4
C_4) \bigr\} \\
&&\qquad\leq C_5 \sqrt{n} \exp\bigl\{- \sqrt{n}
d^2/C_6 \bigr\},
\end{eqnarray*}
where the first inequality is from Lemmas~\ref{GPartI} and~\ref{GPartII}
below.
%the second inequality is from Lemmas~\ref{GPartI} and~\ref{GPartII}
%below.

%---------------- first lemma---------------------------------------
%
%
%le12 #&#
\begin{lemma} \label{GPartI}
Under the assumptions of Theorem~\ref{UnivariateConv}, we have for
\mbox{$1\leq i,j \leq p$},
\[
P \Biggl( \Biggl\vert\sum_{m=1}^N
a_m G_{ij}^{K_m}(1) + 2 \sqrt{
\eta_i \eta_j}\rho_{ij}(0) \Biggr\vert
\geq d \Biggr) \leq C_1 \sqrt{n} \exp\bigl\{- \sqrt{n}
d^2/C_2 \bigr\}.
\]
\end{lemma}

\begin{pf}
From the definition of $\bG^{K_m}=(G_{ij}^{K_m}(1))$ in (\ref{VG1}),
we have
\begin{eqnarray*}
\hspace*{-5pt} G_{ij}^{K_m}(1) &=& \frac{1}{K_m}\sum
_{k_m=1}^{K_m}\sum_{r=2}^{|\btau^{k_m}|}
\bigl[ \varepsilon_i \bigl( \tau_r^{k_m}
\bigr) - \varepsilon_i \bigl( \tau_{r-1}^{k_m}
\bigr) \bigr] \bigl[ \varepsilon_j \bigl( \tau_r^{k_m}
\bigr) - \varepsilon_j \bigl( \tau_{r-1}^{k_m}
\bigr) \bigr]
\\
\hspace*{-5pt}& =& \frac{1}{K_m} \sum_{k_m=1}^{K_m}\sum
_{r=2}^{|\btau^{k_m}|} \bigl[ \varepsilon_i
\bigl( \tau_r^{k_m} \bigr) \varepsilon_j
\bigl( \tau_r^{k_m} \bigr) - \varepsilon_i
\bigl( \tau_r^{k_m} \bigr) \varepsilon_j
\bigl( \tau_{r-1}^{k_m} \bigr) - \varepsilon_i
\bigl( \tau_{r-1}^{k_m} \bigr)\varepsilon_j \bigl(
\tau_r^{k_m} \bigr)\\
&&\hspace*{210pt}{} + \varepsilon_i \bigl(
\tau_{r-1}^{k_m} \bigr)\varepsilon_j \bigl(
\tau_{r-1}^{k_m} \bigr) \bigr]
\\
\hspace*{-5pt}& =& \frac{2}{K_m} \sum_{r=1}^n
\varepsilon_i( t_r) \varepsilon_j(
t_r ) - \frac{1}{K_m} \sum_{r=1}^{K_m}
\varepsilon_i( t_r) \varepsilon_j(
t_r ) - \frac{1}{K_m} \sum_{r=n-K_m+1}^n
\varepsilon_i( t_r) \varepsilon_j(
t_r )
\\
\hspace*{-5pt}&&{} - \frac{1}{K_m} \sum_{r=K_m+1}^n
\varepsilon_i( t_r) \varepsilon_j(
t_{r-K_m} ) - \frac{1}{K_m} \sum_{r=K_m+1}^n
\varepsilon_i( t_{r-K_m}) \varepsilon_j(
t_{r} )
\\
\hspace*{-5pt}& \equiv &I_0^{K_m} - I_1^{K_m} -
I_2^{K_m} - I_3^{K_m} -
I_4^{K_m}
\end{eqnarray*}
and
%
%
%eA.25 #&#
\begin{equation}
\label{termI-1234} \qquad\sum_{m=1}^N
a_m G_{ij}^{K_m}(1) = \sum
_{m=1}^N a_m I_0^{K_m}
- \sum_{i = 1}^4 \sum
_{m=1}^N a_m I_i^{K_m}
\equiv I_0 - I_1 - I_2 - I_3 -
I_4.
\end{equation}
Note that $\sum_{m=1}^N a_m/K_m = 0$, and
\[
\label{termI-0} I_0 = \sum_{m=1}^N
a_m I_0^{K_m} = \sum
_{m=1}^N \frac{a_m}{K_m} \sum
_{r=1}^n \varepsilon_i(
t_r) \varepsilon_j( t_r ) =0.
\]
Hence,
%
%
%eA.26 #&#
\begin{eqnarray}
\label{I1234} && P \Biggl( \Biggl\vert\sum_{m=1}^N
a_m G_{ij}^{K_m}(1) + 2 \sqrt{
\eta_i \eta_j}\rho_{ij}(0) \Biggr\vert
\geq d \Biggr)
\nonumber
\\[-8pt]
\\[-8pt]
\nonumber
&&\qquad\leq\sum_{i=1}^2 P \bigl(
\bigl|I_i - \sqrt{\eta_i \eta_j}
\rho_{ij}(0)\bigr| \geq d/4 \bigr)
 + \sum_{i=3}^4 P \bigl( \vert
I_i \vert\geq d/4 \bigr).
\end{eqnarray}
To prove the lemma we need to derive the four tail probabilities on the
right-hand side of (\ref{I1234}).
Below we will establish the tail probabilities for $I_1, I_2, I_3$ and
$I_4$ by using large deviation results for
the case of $m$-dependent random variables in \citet{SauSta91}.
Because of similarity, we give arguments only for the tail
probabilities of $I_1$ and $I_3$.

First for $I_1$, from the definition of $a_m$ in (\ref{coeff}) we have
%
%
%eA.27 #&#
\begin{eqnarray}
\label{I1-Km}  I_1 - \sqrt{\eta_i
\eta_j}\rho_{ij}(0)& =& \sum_{m=1}^N
a_m \bigl[I_1^{K_m} - \sqrt{
\eta_i \eta_j}\rho_{ij}(0) \bigr],
\nonumber
\\
\qquad P \bigl( \bigl|I_1 - \sqrt{\eta_i \eta_j}
\rho_{ij}(0)\bigr| \geq d/4 \bigr) &\leq&\sum_{m=1}^N
P \bigl( \bigl|I_1^{K_m} - \sqrt{\eta_i
\eta_j}\rho_{ij}(0)\bigr| \geq d/(4 A) \bigr), %&& \sum_{m=1}^N a_m
%G_{ij}^{K_m}(1) + 2 \sqrt{\eta_i \eta_j}
% \sum_{i=2}^3 \sum_{m=1}^N a_m I_i^{K_m} - \sum_{i=4}^5 \sum_{m=1}^N
%a_m
% [I_i^{K_m} - \sqrt{\eta_i \eta_j}\rho_{ij}(0)] \nonumber\\
%&& P \left(\left| \sum_{m=1}^N a_m G_{ij}^{K_m}(1) + 2 \sqrt{\eta_i
% \right|\geq d \right) \leq\sum_{i=2}^3 P\left( \left| \sum_{m=1}^N
% I_i^{K_m} \right| \geq d/4 \right) \nonumber\\
%&& + \sum_{i=4}^5 \sum_{m=1}^N P\left(
% |I_i^{K_m} - \sqrt{\eta_i \eta_j}\rho_{ij}(0)| \geq d/(4 A) \right)
\end{eqnarray}
where $A=\sum_{m=1}^N |a_m| = 9/2 + o(1)$. The $M$-dependence of
$(\varepsilon_1(t_\ell), \ldots, \varepsilon_p(t_\ell))$ in
condition A1 indicates that $\varepsilon_i (t_r) \varepsilon_j (t_r)$,
$r=1, \ldots, n$, are $M$-dependent,
$I_1^{K_m}$ is the average of $\varepsilon_i (t_r) \varepsilon_j
(t_r)$, $r=1, \ldots, K_m$, and
Lemma~\ref{ELmoments} below calculates $E(I_1^{K_m}) = \sqrt{\eta
_i\eta
_j}\rho_{ij}(0)$ and $\operatorname{Var} (I_1^{K_m}) \leq C n^{-1/2}$.
Also for any integer~$k$,
\begin{eqnarray*}
 E \bigl( \bigl\vert\varepsilon_i (t_r)
\varepsilon_j (t_r) \bigr\vert^k \bigr)
&\leq&\sqrt{E \bigl( \bigl\vert\varepsilon_i (t_r)
\bigr\vert^{2k} \bigr) E \bigl( \bigl\vert\varepsilon_j
(t_r) \bigr\vert^{2k} \bigr) }
\\
& \leq& c_0 (2k)! (2 \tau_0)^{2k} \leq
c_0 (k!)^2 \bigl(16 \tau_0^2
\bigr)^{k} \leq(k!)^2 \bigl[16 \tau_0^2
(c_0 \vee1) \bigr]^{k},
\end{eqnarray*}
%
%(2k-1)!! C_1^{-2k} \leq(k!)^2(4C_1^{-2})^k,$$
where the first inequality is from the Cauchy--Schwarz inequality, and
the second inequality is from
the subGaussian tails of $\varepsilon_i (t_r)$ and $\varepsilon_j
(t_r)$, which imply that their
$2k$-moments are bounded by $\int_0^\infty c_0 \exp[-x^{1/(2k)}/(2
\tau
_0)] \,dx = c_0 (2 k)! (2 \tau_0)^{2 k}$.
%Hence, $M$-dependent random variables $ \varepsilon_i (t_r)
%in Saulis and Statulevicius (1991, page 105), %with $H_1= 2 \tau_0$
%and some $\gamma_1>0$,
Applying Theorem~4.30 in \citet{SauSta91} to
$M$-dependent random variables
$ \varepsilon_i (t_r) \varepsilon_j (t_r) $ we obtain
%
%
%eA.28 #&#
\begin{equation}
\label{Mdep1} P \bigl( \bigl\vert I_1^{K_m} - \sqrt{
\eta_i\eta_j}\rho_{ij}(0) \bigr\vert\geq
d_1 \bigr) \leq C_1 \exp\biggl\{ - \frac{\sqrt{n}d_1^2}{C_2}
\biggr\}.
\end{equation}
Plugging (\ref{Mdep1}) with $d_1=d/(4 A)$ into (\ref{I1-Km}) we
establish the tail probability for $I_1$
%
%
%eA.29 #&#
\begin{eqnarray}
\label{I1-tail} P \bigl( \bigl|I_1 - \sqrt{\eta_i
\eta_j}\rho_{ij}(0)\bigr| \geq d/4 \bigr) &\leq& C_1
N \exp\biggl\{ - \frac{\sqrt{n} d^2}{16 A^2 C_2} \biggr\}
\nonumber
\\[-8pt]
\\[-8pt]
\nonumber
& \leq& C_3
\sqrt{n} \exp
\biggl\{ - \frac{\sqrt{n} d^2}{C_4} \biggr\}.
\end{eqnarray}

Second, consider $I_3$. We may express it as follows:
\[
I_3 = \sum_{m=1}^N \sum
_{r=K_m+1}^n \frac{a_m}{K_m}
\varepsilon_i (t_r)\varepsilon_j
(t_{r-K_m}) %%=\sum_{r=1}^{n-K_1} \sum_{m=1}^{(n+1-r-K_1)\wedge N}
= \sum
_{r=1}^{n-K_1} \sum_{m=1}^{(n - N -r)\wedge N}
\frac{a_m}{K_m} \varepsilon_j (t_r)
\varepsilon_i (t_{r+K_m}),
\]
%
%n+1-r -K_1
and Lemma~\ref{ELmoments} below derives $E(I_3) = 0$ and $\operatorname{Var}(I_3)
\leq C n^{-1/2}$.

As $(\varepsilon_1(t_{\ell}), \ldots, \varepsilon_p(t_{\ell}))$,
$\ell
=1, \ldots, n$, are serially $M$-dependent,
that is, for any integers $k$ and $k^\prime$, and
integer sets $\{\ell_1, \ldots, \ell_k\}$ and $\{\ell^\prime_1,
\ldots,
\ell^\prime_{k^\prime}\}$,
$\{\varepsilon_i(t_{\ell_1}),\break \ldots,  \varepsilon_i(t_{\ell_k}), i=1,
\ldots, p\}$ and
$\{\varepsilon_i(t_{\ell^\prime_1}), \ldots, \varepsilon_i(t_{\ell
^\prime_{k^\prime}}), i=1, \ldots, p\}$
are independent if every integer in $\{\ell_1, \ldots, \ell_k\}$
differs by more than $M$ from any integer in
$\{\ell^\prime_1, \ldots, \ell^\prime_{k^\prime}\}$.
Since $K_m >M$ for $n$ large enough, if integers $r$ and $r^\prime$
differ by more than $K_N+M$,
for two integer sets $\{r, r+K_m; m=1,\ldots, (n - N -r)\wedge N\}$ and
$\{r^\prime, r^\prime+K_m; m=1,\ldots, (n - N -r^\prime)\wedge N\}$,
every element in one integer set must be more than $M$ apart from any
element in the other integer set. Then
$\{\varepsilon_j (t_r), \varepsilon_i (t_{r+K_m}); m=1,\ldots, (n - N
-r)\wedge N\}$ and
$\{\varepsilon_j (t_{r^\prime}), \varepsilon_i (t_{r^\prime+K_m});
m=1,\ldots, (n - N -r^\prime)\wedge N\}$ are independent,
and thus $\varepsilon_j (t_r)\varepsilon_i (t_{r+K_m})$, $r=1, \ldots,
n-K_m$, are serially $(K_N+M)$-dependent. Also for any integer $k$,
\begin{eqnarray*}
 E \bigl(\bigl|\varepsilon_j (t_r)\varepsilon_i
(t_{r+K_m})\bigr|^{k} \bigr) &\leq&\sqrt{ E \bigl(\bigl|
\varepsilon_j (t_r)\bigr|^{2k} \bigr) E \bigl(\bigl|
\varepsilon_i (t_{r+K_m})\bigr|^{2k} \bigr) }
\\
& \leq& c_0 (2k)! (2 \tau_0)^{2k} \leq
c_0 (k!)^2 \bigl(16 \tau_0^2
\bigr)^{k} \leq(k!)^2 \bigl[16 \tau_0^2
(c_0 \vee1) \bigr]^{k},
\end{eqnarray*}
where the first inequality is from the Cauchy--Schwarz inequality, and
the second inequality is from
the subGaussian tails of $\varepsilon_j (t_r)$ and $\varepsilon_i
(t_{r+K_m})$.
Applying theorem 4.16 in \citet{SauSta91} we derive a
bound $(k!)^3 C_0^{k}$ on the
$k$th cumulant of $ n^{1/4} I_3$, and then using
Lemmas 2.3 and 2.4 in \citet{SauSta91} we establish the
tail probability for $I_3$ as follows:
%
%
%eA.30 #&#
\begin{equation}
\label{I3-tail} P\bigl ( |I_3| \geq d/4 \bigr) \leq C_1 \exp
\biggl\{ - \frac{\sqrt{n} (d/4)^2}{C_2} \biggr\} \leq C_3 \exp
\biggl\{ -
\frac{\sqrt{n} d^2}{C_4} \biggr\}.
\end{equation}

%%-----------------Calculate cumulants
%Theorem~4.16 (1): $k$-th cumulant of $n^{1/4} I_3$
%c_0^2 (k!)^2 (2 \tau_0)^{2k} } \]

%Theorem~4.30: the $k$-moment for each term of $n^{1/4} I_3$
%(t_{r-K_m})|^{k}) \leq(n^{-3/4})^k (k!)^2 (2 \tau_0)^{k}, \]
%$k$-th cumulant of $n^{1/4} I_3$

%Using the same conditional argument trick used in the proof of Lemma (
%$$P(|I_2 |\geq d) = E(P(|I_2 |\geq d | \varepsilon_i \text{'s in the
%terms})) \leq E \left(
%C_1 \exp\left\{ - \frac{d^2}{C_2 R} \right\} \right),$$
%where the inequality is due to the subGaussian tail of $\varepsilon_j
%$ and the fact that $\varepsilon_i (t_r)$ and $\varepsilon_j (t_r')$
%in $I_2 $ are independent since the lags are at least $K_1$ and $2M <
%K_1$. The term $R$ is defined as
%$$R = \sum_{r=1}^{n-K_1} \left(\sum_{m=1}^{(n+1-r-K_1)\wedge N}
%Set $\Omega= \{ |R-ER|\leq d \}$, where $ER = C/N + \text{ higher
%order}$. Due to the subGaussian tail of $\varepsilon_i $, we have
%$$P(\Omega^C) = P(|R-ER|\geq d) \leq C_1 \exp\{ -\frac{
%& P(|I_2 |\geq d) \leq E \left( C_1 \exp\left\{ - \frac{d^2}{C_2 R}
%& = E \left( C_1 \exp\left\{- \frac{d^2}{C_2 R} \right\} 1(\Omega)
%C_1 \exp\left\{- \frac{d^2}{C_2 R} \right\} 1(\Omega^C) \right) \\
%&\leq C_1 \exp\left\{ - \sqrt{n} d^2/C_2 \right\} + C_1 \exp\left\{ -
%d^2}{C_2 M}\right\}.

Since $I_2$ and $I_4$ have the same tail probabilities as $I_1$ and
$I_3$ given by (\ref{I1-tail}) and~(\ref{I3-tail}),
respectively, combining them with (\ref{I1234}) we conclude
\begin{eqnarray*}
&& P \Biggl( \Biggl\vert\sum_{m=1}^N
a_m G_{ij}^{K_m}(1) + 2 \sqrt{
\eta_i \eta_j}\rho_{ij}(0) \Biggr\vert
\geq d \Biggr)\\
&&\qquad \leq2 C_3 \sqrt{n} \exp\biggl\{ - \frac{\sqrt{n} d^2}{C_4}
\biggr\} + 2 C_3 \exp\biggl\{ - \frac{\sqrt{n} d^2}{C_4} \biggr\}
\\
&&\qquad\leq C_5 \sqrt{n} \exp\biggl\{ - \frac{\sqrt{n} d^2}{C_6} \biggr
\}.
\end{eqnarray*}
\upqed\end{pf}

%------- second lemma, the tail for

%
%le13 #&#
\begin{lemma} \label{GPartII}
Under the assumptions of Theorem~\ref{UnivariateConv}, we have for
\mbox{$1\leq i,j \leq p$},
%
%
%eA.31 #&#
\begin{equation}\quad
P \bigl( \bigl\vert\zeta\bigl(G_{ij}^{K_1}(1) -
G_{ij}^{K_N}(1) \bigr) - 2 \sqrt{\eta_i
\eta_j}\rho_{ij}(0) \bigr\vert\geq d \bigr) \leq
C_1 \exp\bigl\{-n d^2/C_2 \bigr\}.
\end{equation}
\end{lemma}

\begin{pf}
First consider $\zeta G_{ij}^{K_1}(1)$ term:
\begin{eqnarray*}
 \zeta G_{ij}^{K_1}(1)& = &\frac{K_N}{n(N-1)}\sum
_{k_1=1}^{K_1}\sum_{r=2}^{|\btau^{k_1}|}
\bigl( \varepsilon_i \bigl( \tau_r^{k_1}
\bigr) - \varepsilon_i \bigl( \tau_{r-1}^{k_1}
\bigr) \bigr) \bigl( \varepsilon_j \bigl( \tau_r^{k_1}
\bigr) - \varepsilon_j \bigl( \tau_{r-1}^{k_1}
\bigr) \bigr)
\\
& =& \frac{K_N}{n(N-1)} \sum_{r=K_1+1}^{n}
\bigl( \varepsilon_i( t_r) \varepsilon_j(
t_r ) + \varepsilon_i( t_{r-K_1})
\varepsilon_j(t_{r-K_1}) \\
&&\hspace*{80pt}{}- \varepsilon_i(t_r)
\varepsilon_j( t_{r-K_1}) - \varepsilon_i(
t_{r-K_1})\varepsilon_j(t_r) \bigr)
\\
& \equiv& R_1 + R_2 - R_3 -
R_4.
\end{eqnarray*}
Due to similarity, we show the tail probabilities only for $R_1$ and
$R_3$. Lemma~\ref{ELmoments} below calculates the mean
and variances of $R_1 $ and $R_3$. Since $R_1$ and $R_3$ have,
respectively, the same structures as $I_1^{K_m}$ and $I_3$ used in the
proof Lemma~\ref{GPartI}, the arguments for establishing the tail
probabilities for $I_1^{K_m}$ and $I_3$ can be used to derive
the tail probability bounds for $R_1$ and $R_3$. Consequently we obtain that
%
%
%eA.32 #&#
\begin{equation}
\label{G1-K1} \qquad P \biggl( \biggl\vert\zeta G_{ij}^{K_1}(1)
- \frac{2 K_N(n-K_1)}{n(N-1)} \sqrt{\eta_i\eta_j}
\rho_{ij}(0) \biggr\vert\geq d \biggr) \leq C_1 \exp
\bigl\{-n d^2/C_2 \bigr\}.
\end{equation}
As $G_{ij}^{K_N}(1)$ has the same structure as $\zeta G_{ij}^{K_1}(1)$,
similarly we can establish a tail probability for
$\zeta G_{ij}^{K_N}(1)$ as follows:
%
%
%eA.33 #&#
\begin{equation}
\label{G1-KN}  P \biggl( \biggl\vert\zeta G_{ij}^{K_N}(1)
- \frac{2 K_1(n-K_N)}{n(N-1)} \sqrt{\eta_i\eta_j}
\rho_{ij}(0) \biggr\vert\geq d \biggr) \leq C_1 \exp
\bigl\{-n d^2/C_2 \bigr\}.\hspace*{-25pt}
\end{equation}
Since
\[
\frac{K_N(n-K_1)}{n(N-1)} - \frac{K_1(n-K_N)}{n(N-1)} = 1,
\]
combining (\ref{G1-K1}) and (\ref{G1-KN})
%the tail probabilities for $\zeta G_{ij}^{K_1}(1)$ and $\zeta
%G_{ij}^{K_N}(1)$
we prove the lemma.
% P \left( \left| \zeta\left(G_{ij}^{K_1}(1) - G_{ij}^{K_N}(1)\right)
%- 2 \eta_i 1(i=j)
% \right|\geq d/8 \right) \leq C_1 \exp\{- n d^2/C_2 \}.
\end{pf}

%Saulis and Statulevi\v{c}ius (1991, page 46)] \label{Thm33}
%Under the definitions and assumptions of Theorem~3.3 in Saulis and
%Statulevi\v{c}ius (1991), if we allow $\{ a_{j,n} \}$ to be any array
%of
%numbers (not necessary nonnegative numbers), and re-define
%$\gamma_n = \max\{|a_{j,n}|, 1\leq j\leq n\}$, then all the results in
%Theorem~3.3 still valid.
%
%From Assumption $(B_\gamma)$ in Theorem~3.3 of Saulis and
%Statulevi\v{c}ius (1991), we immediately show that Inequality (3.20)
%in Saulis
%and Statulevi\v{c}ius (1991) still holds under the modified setting
%and thus
%prove the lemma.
%We need to show only that under the modified setting, inequality
%(3.20) in
%Saulis and Statulevi\v{c}ius (1991) still holds. With the notation
%below
%from Saulis and Statulevi\v{c}ius (1991) [not related to those in this
%paper], from assumption $(B_\gamma)$ in Theorem~3.3 of Saulis and
%Statulevi\v{c}ius (1991), we have
%& =\left| \sum_{j=1}^n a_{j,n}^k \Gamma_k(\xi_j) \right|
%=\left| \sum_{j=1}^n a_{j,n}^2 a_{j,n}^{k-2}\Gamma_k(\xi_j) \right| \\
%& \leq\sum_{j=1}^n a_{j,n}^2\left| a_{j,n}^{k-2} \right|
% \leq\sum_{j=1}^n a_{j,n}^2 \gamma_n^{k-2} (k )^{1+\gamma} (2\max\{K,
%& = (k )^{1+\gamma} (K_n \gamma_n)^{k-2} \tilde{B}_n^2,
%which is (3.20) in Saulis and Statulevi\v{c}ius (1991).

%-------------------------------------------- Lemma for moments
%--------------------------------------------
%
%
%le14 #&#
\begin{lemma} \label{ELmoments}
Under the assumptions of Theorem~\ref{UnivariateConv} and for large
enough $n$ so that $M<K_1$, we have
\begin{eqnarray*}
 E(I_3)& =& E(R_3) = 0,\qquad E \bigl(I_1^{K_m}
\bigr) = \sqrt{\eta_i \eta_j}\rho_{ij}(0),\\[2pt]
E(R_1) &= &\frac{K_N(n-K_1)}{n(N-1)} \sqrt{\eta_i
\eta_j}\rho_{ij}(0),
\\[2pt]
\operatorname{Var} \bigl(I_1^{K_m} \bigr) &\leq &C n^{-1/2},\qquad \operatorname{Var}
(I_3 ) \leq C n^{-1/2}, \qquad\operatorname{Var} (R_1 ) \leq C
n^{-1},\\[2pt]
 \operatorname{Var} (R_3 )& \leq& C n^{-1}. %= \frac{C}{\sqrt{n}}(1+o(1)), \operatorname{Var}
%&& \operatorname{Var} \left(R_1 \right) = \frac{C}{n}(1+o(1)), \operatorname{Var} \left(R_3 \right)
%= \frac{C}{n}(1+o(1)).
\end{eqnarray*}
\end{lemma}

\begin{pf} Because $K_m > M$, $\varepsilon_i (t_r)$ and $\varepsilon_j
(t_{r-K_m})$ are independent, so
\begin{eqnarray*}
E(I_3) %= E \left( \sum_{m=1}^N \frac{a_m}{K_m} \sum_{r = K_m +1}^n
&=& \sum
_{m=1}^N \frac{a_m}{K_m} \sum
_{K_m +1}^n E \bigl[\varepsilon_i
(t_r) \varepsilon_j (t_{r-K_m}) \bigr]\\[2pt]
 &=& \sum
_{m=1}^N \frac{a_m}{K_m} \sum
_{K_m +1}^n E \bigl[\varepsilon_i
(t_r) \bigr] E \bigl[\varepsilon_j (t_{r-K_m})
\bigr] =0,
\\[2pt]
%
%where the last equality is due to the fact that there is at least $K_1
%> 2M$ lags, and thus $\varepsilon_i (t_r)$ and $ \varepsilon_j
%(t_{r-K_m})$ are independent, and the fact that $\varepsilon_i(\cdot),
%i = 1,\ldots,p$ are mean zero.
%
E (R_3) &= &\frac{K_N}{n(N-1)} \sum_{r=K_1+1}^n
E \bigl[\varepsilon_i (t_r) \varepsilon_j
(t_{r-K_1}) \bigr] \\[2pt]
&=&\frac{K_N}{n(N-1)} \sum_{r=K_1+1}^n
E \bigl[\varepsilon_i (t_r) \bigr] E \bigl[
\varepsilon_j (t_{r-K_1}) \bigr] = 0.
\end{eqnarray*}
For $I_1^{K_m}$ and $R_1$, we have
\begin{eqnarray*}
E \bigl( I_1^{K_m} \bigr) & =& \frac{1}{K_m} \sum
_{r = 1}^{K_m} E \bigl[ \varepsilon_i
(t_r) \varepsilon_j (t_{r}) \bigr] =
\frac{1}{K_m} \sum_{r = 1}^{K_m} \sqrt{
\eta_i \eta_j}\rho_{ij}(0) = \sqrt{
\eta_i \eta_j}\rho_{ij}(0),
\\[2pt]
E R_1 & =& \frac{K_N}{n(N-1)} \sum_{r=K_1+1}^n
E \bigl[\varepsilon_i (t_r) \varepsilon_j
(t_{r}) \bigr] = \frac{K_N(n-K_1)}{n(N-1)} \sqrt{\eta_i
\eta_j}\rho_{ij}(0).
\end{eqnarray*}
With the $M$-dependence of $\varepsilon_i ( t_r) \varepsilon_j ( t_r
)$, we directly compute the variances of $I_1^{K_m}$ and $R_1$ as
follows:
\begin{eqnarray*}
 \operatorname{Var} \bigl(I_1^{K_m} \bigr) %= \frac{1}{K_m^2} \operatorname{Var}\left(
&=&
\frac{1}{K_m^2} \sum_{r=1}^{K_m} \operatorname{Var}
\bigl(\varepsilon_i ( t_r) \varepsilon_j
(t_r) \bigr) \\[2pt]
&&{}+ \frac{2}{K_m^2} \mathop{\sum\sum
}_{1 \leq r<r'\leq K_m} \operatorname{Cov} \bigl(\varepsilon_i (
t_r) \varepsilon_j (t_r),
\varepsilon_i ( t_{r'}) \varepsilon_j
(t_{r'}) \bigr)
\\[2pt]
&\leq&\frac{1}{K_m} \operatorname{Var} \bigl(\varepsilon_i (
t_1) \varepsilon_j (t_1) \bigr) \\[2pt]
&&{}+
\frac{2}{K_m} \sum_{\ell=2}^{M+1} \operatorname{Cov}
\bigl(\varepsilon_i ( t_1) \varepsilon_j
(t_1), \varepsilon_i ( t_{\ell})
\varepsilon_j (t_{\ell}) \bigr) \leq C n^{-1/2},
\\
 \operatorname{Var}(R_1 )& =& \biggl(\frac{K_N}{n(N-1)} \biggr)^2
\Biggl[ \sum_{r=K_1+1}^{n} \operatorname{Var} \bigl(
\varepsilon_i ( t_r) \varepsilon_j
(t_r) \bigr) \\
&&\hspace*{66pt}{}+ 2 \mathop{\sum\sum}_{K_1 + 1 \leq r<r'\leq n} \operatorname{Cov} \bigl(
\varepsilon_i ( t_r) \varepsilon_j
(t_r), \varepsilon_i ( t_{r'})
\varepsilon_j (t_{r'}) \bigr) \Biggr]
\\
&\leq&\biggl(\frac{K_N}{n(N-1)} \biggr)^2 \Biggl[ (n-
K_1) \operatorname{Var} \bigl(\varepsilon_i ( t_1)
\varepsilon_j (t_1) \bigr) \\
&&\hspace*{69pt}{}+ 2 (n-K_1) \sum
_{\ell=2}^{M+1} \operatorname{Cov} \bigl(\varepsilon_i
( t_1) \varepsilon_j (t_1),
\varepsilon_i ( t_{\ell}) \varepsilon_j
(t_{\ell}) \bigr) \Biggr]
\\
& \leq& C / n.
\end{eqnarray*}
We evaluate the variance of $I_3$ as follows:
\begin{eqnarray*}
 E \bigl(I_3^2 \bigr) &=& % E \left(\sum_{m=1}^N \frac{a_m}{K_m}
% = &
\sum_{m=1}^N \biggl(
\frac{a_m}{K_m} \biggr)^2 E \Biggl(\sum
_{r = K_m
+1}^n \varepsilon_i
(t_r) \varepsilon_j (t_{r-K_m})
\Biggr)^2
\\
&&{}+ 2 \mathop{\sum\sum}_{m<m'} \frac{a_m}{K_m} \frac{a_{m'}}{K_{m'}} E
\Biggl[ \Biggl(\sum_{r = K_m +1}^n
\varepsilon_i (t_r) \varepsilon_j
(t_{r-K_m}) \Biggr) \\
&&\hspace*{103pt}{}\times\Biggl(\sum_{r' = K_{m'} +1}^n
\varepsilon_i \bigl(t_r' \bigr)
\varepsilon_j (t_{r'-K_{m'}}) \Biggr) \Biggr]
\\
&= & \sum_{m=1}^N \biggl(
\frac{a_m}{K_m} \biggr)^2 \Biggl[\sum_{r=K_m+1}^n
E \bigl(\varepsilon^2_i (t_r)
\varepsilon^2_j (t_{r-K_m}) \bigr)\\
&&\hspace*{57pt}{} + 2 \mathop{\sum
\sum} _{r<r'} E \bigl(\varepsilon_i
(t_r)\varepsilon_j (t_{r-K_m})
\varepsilon_i (t_{r'})\varepsilon_j
(t_{r'-K_{m}}) \bigr) \Biggr]
\\
&&{} + 2\mathop {\sum\sum}_{m<m'} \frac{a_m}{K_m} \frac{a_{m'}}{K_{m'}}
\Biggl[ \sum_{r = K_{m'} +1 }^n E \bigl(
\varepsilon^2_i (t_r) \varepsilon_j
(t_{r-K_m}) \varepsilon_j (t_{r-K_{m'}}) \bigr)\\
&&\hspace*{92pt}{} + \sum
_{r=1}^{n-K_{m'}} E \bigl( \varepsilon_j^2(t_r)
\varepsilon_i(t_{r+K_m})\varepsilon_i(t_{r+K_{m'}})
\bigr)
\\
&&\hspace*{92pt}{} + 2 \mathop{\sum\sum}_{r<r'} E \bigl( \varepsilon_i
(t_r)\varepsilon_j (t_{r-K_m})
\varepsilon_i (t_{r'})\varepsilon_j
(t_{r'-K_{m'}}) \bigr) \Biggr]
\\
&= & \sum_{m=1}^N \biggl(
\frac{a_m}{K_m} \biggr)^2 \Biggl[ (n-K_m)\eta
_i\eta_j \\
&&\hspace*{60pt}{}+ 2 \sum_{\ell= 2}^{(n-K_m)\wedge(M+1)}
(n-K_m - \ell+1) E \bigl(\varepsilon_i(t_1)
\varepsilon_i(t_\ell) \bigr)\\
&&\hspace*{202pt}{}\times E \bigl(\varepsilon
_j(t_1)\varepsilon_j(t_\ell)
\bigr) \Biggr]
\\
&&{} + 2 \mathop{\sum\sum}_{m<m'<m+M+1} \frac{a_m}{K_m} \frac
{a_{m'}}{K_{m'}}
\Biggl[ (n-K_{m^\prime})\eta_i E \bigl( \varepsilon_j(t_1)
\varepsilon_j(t_{K_{m^\prime} - K_m+1}) \bigr)\\
&&\hspace*{126pt}{} + (n-K_m)
\eta_j  E \bigl(\varepsilon_i(t_n)
\varepsilon_i(t_{n-K_{m^\prime} +
K_m}) \bigr)
\\
& &{}+ 2 \sum_{\ell= 2}^{(n-K_{m^\prime})\wedge(M+1)}(n -
K_{m^\prime} - \ell+1) E \bigl(\varepsilon_i(t_1)
\varepsilon_i(t_\ell) \bigr)\\
&&\hspace*{179pt}{}\times E \bigl(
\varepsilon_j(t_1) \varepsilon_j(t_{\ell+K_{m^\prime}-K_m})
\bigr) \Biggr]
\\
%& \leq C_1 \eta_i\eta_j N n /N^4 + C_2 \eta_i \eta_j N \left(
%& + C_3 \eta_i \eta_j\left( \frac{1}{N^2}\right)^2 \sum_{\Delta_m =
%1}^{M} (N - \Delta_m)(n-K_m)(|\rho_i(\Delta_m)| + |\rho_j(\Delta_m)|) \
%& + C_4 \eta_i \eta_j\left( \frac{1}{N^2}\right)^2 N n \sum_{
%+ \Delta_m)| \\[-2pt]
& \leq& C_1
\eta_i \eta_j N \bigl(1/N^2
\bigr)^2 (n-K_1) + C_2 \eta_i
\eta_j \bigl(1/N^2 \bigr)^2
(n-K_1) \asymp C n^{-1/2},
\end{eqnarray*}
where the inequality is from the fact that the $M$-dependence of
$(\varepsilon_1(t_\ell), \ldots,\break \varepsilon_p(t_\ell))$ implies zero
expectations of $\varepsilon_i(\cdot) \varepsilon_j(\cdot)$ for lags
larger than $M$.

Similarly, we have
\begin{eqnarray*}
 E \bigl( R_3^2 \bigr) &=& \biggl(\frac{K_N}{n(N-1)}
\biggr)^2 E \Biggl(\sum_{r = K_1 +1}^n
\varepsilon_i (t_r) \varepsilon_j
(t_{r-K_1}) \Biggr)^2
\\
&= & \biggl(\frac{K_N}{n(N-1)} \biggr)^2 \Biggl[ \sum
_{r = K_1 +1}^n E \bigl(\varepsilon^2_i
(t_r) \varepsilon^2_j (t_{r-K_1})
\bigr)\\
&&\hspace*{68pt}{}+ 2 \mathop{\sum\sum}_{r<r'} E \bigl(\varepsilon_i
(t_r)\varepsilon_j (t_{r-K_1})
\varepsilon_i (t_{r'})\varepsilon_j
(t_{r'-K_{1}}) \bigr) \Biggr]
\\
&=& \biggl(\frac{K_N}{n(N-1)} \biggr)^2 \Biggl[ (n-K_1)
\eta_i \eta_j \\
&&\hspace*{68pt}{}+ 2 \sum_{\ell=2}^{(n-K_{1})\wedge(M+1)}
(n-K_1 - \ell+1) E \bigl(\varepsilon_i(t_1)
\varepsilon_i(t_\ell) \bigr)\\
&&\hspace*{208pt}{}\times E \bigl(
\varepsilon_j(t_1) \varepsilon_j(t_\ell
) \bigr) \Biggr]
\\
& \leq& C_1 \eta_i \eta_j
(1/n)^2 (n- K_1) \asymp C_2 /n, %& \leq C_1 \eta_i\eta_j /n + C_2
%K_1)| \leq C/n.
%& \leq C \eta_i\eta_j /n \sum_{\Delta_r} |\rho_i(\Delta_r)|
\end{eqnarray*}
where the inequality is from the fact that the $M$-dependence of
$(\varepsilon_1(t_\ell), \ldots,\break \varepsilon_p(t_\ell))$ implies zero
expectations of $\varepsilon_i(\cdot) \varepsilon_j(\cdot)$ for lags
larger than $M$.
\end{pf}

%
%sB #&#
\section{\texorpdfstring{Proof of Proposition \lowercase{\protect\ref{affbd}}}{Proof of Proposition 9}}
\label{lemmaaffsec}

We break the proof into a few major technical lemmas which are proved in
Sections~\ref{chisquaresec}--\ref{secRterms}. %\ref{secAlem}.
Without loss of generality we
consider only the case $i=1$ and prove that there exists a constant $C_{1}>0$
such that $\llVert\bar{\mathbb{P}}_{1,0}\wedge\bar{\mathbb{P}}%
_{1,1}\rrVert\geq C_{1}$.

The following lemma turns the problem of bounding the total variation
affinity into a chi-square distance calculation.
%For any $b\subseteq B$, let
%D_{\Lambda_{-i}}=\mathrm{Card}\left\{ \lambda:\lambda_{i}(\theta)=b
%which does not depend on values of $i$ and $b$ by the symmetry of the
%construction.
Denote the projection of $\theta\in\Theta$ to $\bGamma$ by $%
\gamma( \theta) = ( \gamma_{i} ( \theta)
) _{1\leq i\leq r}$ and to $\Lambda$ by $\lambda( \theta
) = ( \lambda_{i} ( \theta) ) _{1\leq i\leq r}$.
More generally, for a subset $A\subseteq\{ 1,2,\ldots,r \}
$, we
define a projection of $\theta$ to a subset of $\bGamma$ by $\gamma
_{A} ( \theta) = ( \gamma_{i} ( \theta)
)
_{i\in A}$. A particularly useful example of set $A$ is $ \{
1,\ldots,i-1,i+1,\ldots,r \}$
for which we use $\gamma_{-i} ( \theta) = ( \gamma
_{1} (
\theta),\ldots,\gamma_{i-1} ( \theta),\gamma
_{i+1} ( \theta),\gamma_{r} ( \theta) )
$. $%
\lambda_{A} ( \theta) $ and $\lambda_{-i} ( \theta
) $
are defined similarly. We define the set $\Lambda_{A}= \{ \lambda
_{A} ( \theta) \dvtx\theta\in\Theta\} $. For $a\in
\{
0,1 \} $, $b\in\{ 0,1 \} ^{r-1}$, and $c\in\Lambda
_{-i}\subseteq B^{r-1}$, let
\[
\Theta_{(i,a,b,c)}= \bigl\{ \theta\in\Theta\dvtx\gamma_{i}(
\theta)=a,\gamma_{-i}(\theta)=b\mbox{ and }\lambda_{-i}(
\theta)=c \bigr\},
\]
and $D_{(i,a,b,c)}=\mathrm{Card}(\Theta_{(i,a,b,c)})$ which depends
actually on the value of $c$, not values of $i$, $a$ and $b$ for the
parameter space $\Theta$ constructed in Section~\ref{lowerboundproofssc}. %Section~\ref{lowerbdsec}.
Define the mixture distribution
%
%
%eB.1 #&#
\begin{equation}
\bar{\mathbb{P}}_{ ( i,a,b,c ) }=\frac{1}{D_{(i,a,b,c)}}%
\sum
_{\theta\in\Theta_{(i,a,b,c)}}\mathbb{P}_{\theta}. \label{avepibd}
\end{equation}
In other words, $\bar{\mathbb{P}}_{(i,a,b,c)}$ is the mixture distribution
over all $\mathbb{P}_{\theta}$ with $\lambda_{i}(\theta)$ varying over
all possible values while all other components of $\theta$ remain fixed.
Define
\[
\Theta_{-1}= \bigl\{ ( b,c ) \dvtx\mbox{there exists a }\theta\in
\Theta\mbox{ such that }\gamma_{-1}(\theta)=b\mbox{ and }\lambda
_{-1}(\theta)=c \bigr\}.
\]

%
%le15 #&#
\begin{lemma}
\label{chisquareslem}If there is a constant $C_{2}<1$ such that
%
%
%eB.2 #&#
\begin{equation}
\mathop{\mathrm{Average}}_{ ( \gamma_{-1},\lambda_{-1} ) \in\Theta
_{-1}} \biggl\{ \int\biggl( \frac{d\bar{\mathbb{P}}_{ (
1,1,\gamma
_{-1},\lambda_{-1} ) }}{d\bar{\mathbb{P}}_{ ( 1,0,\gamma
_{-1},\lambda_{-1} ) }}
\biggr) ^{2}\,d\bar{\mathbb{P}}_{ (
1,0,\gamma_{-1},\lambda_{-1} ) }-1 \biggr\} \leq
C_{2}^{2}, \label{chi-squarebd}
\end{equation}
then $\llVert\bar{\mathbb{P}}_{1,0}\wedge\bar{\mathbb{P}}%
_{1,1}\rrVert\geq1-C_{2}>0$.
\end{lemma}
We can prove Lemma~\ref{chisquareslem} using the same arguments as the proof
of Lemma~8 in \citet{CaiZho}. To complete the proof of Proposition~\ref{affbd}
%Theorems~\ref{lowerbound} and Theorem~\ref{Operlowerbdthm}
we need to verify only equation (\protect\ref{chi-squarebd}).

\renewcommand{\thesubsection}{II.\arabic{subsection}}
\setcounter{subsection}{0}
%sB.1 #&#
\subsection{\texorpdfstring{Technical lemmas for proving equation (\protect\ref{chi-squarebd})}
{Technical lemmas for proving equation (71)}}
\label{chisquare-1sec}

From the definition of $\bar{\mathbb{P}}_{ ( 1,0,\gamma
_{-1},\lambda_{-1} ) }$ in equation (\ref{avepibd}) and
$\theta=(\gamma,\lambda)$ with $\gamma=(\gamma_1, \ldots, \gamma
_r)$ and
$\lambda=(\lambda_1, \ldots, \lambda_r)$,
%in Equation (\ref{theta}), ???
$\gamma_{1}=0$ implies $\bar{\mathbb{P}}_{ (
1,0,\gamma_{-1},\lambda_{-1} ) }$ is a product of $n$ multivariate
normal distributions each with a covariance matrix,
%
%
%eB.3 #&#
\begin{equation}
\qquad\Sigma_{l,0}= \pmatrix{ 1 &
\mathbf{0}_{1\times( p-1 ) }
\vspace*{2pt}\cr
\mathbf{0}_{ ( p-1 ) \times1} & \mathbf{S}_{ ( p-1 )
\times( p-1 ) }} + ( a_{l}-1 ) \bI_{p}\qquad\mbox{for }l=1,2,\ldots,n
\label{sigma0},
\end{equation}
where $\mathbf{S}_{ ( p-1 ) \times( p-1 ) }= (
s_{ij} ) _{2\leq i,j\leq p}$ is uniquely determined by $ (
\gamma
_{-1},\lambda_{-1} ) = ( (\gamma_{2},\ldots,\gamma
_{r}),(\lambda
_{2},\ldots,\lambda_{r}) ) $ with
\[
s_{ij}= \cases{ %
 1, & \quad$i=j$,
\vspace*{2pt}\cr
\epsilon_{n,p}, &\quad $\gamma_{i}=\lambda_{i}
( j ) =1$,
\vspace*{2pt}\cr
0, & \quad$\mbox{otherwise}.$}
\]
Let $n_{\lambda_{-1}}$ be the number of columns of $\lambda_{-1}$ with
column sum equal to $2k$ and $p_{\lambda_{-1}}= r %\left\lceil p/2
-n_{\lambda_{-1}}$. Since $n_{\lambda_{-1}}\cdot2k\leq r %\left
\cdot k$, the total number of $1$s in the upper triangular matrix, we have
$n_{\lambda_{-1}}\leq r /2$, which implies $p_{\lambda_{-1}}= r
-n_{\lambda_{-1}}\geq r /2\geq p/4-1$.
From equations (\ref{avepibd}) and
$\theta=(\gamma,\lambda)$ with $\gamma=(\gamma_1, \ldots, \gamma
_r)$ and
$\lambda=(\lambda_1, \ldots, \lambda_r)$, %and (\ref{theta}) ???
$\bar{\mathbb{P}}_{ ( 1,1,\gamma
_{-1},\lambda_{-1} ) }$ is an average of ${p_{\lambda
_{-1}}\choose k}$
number of products of multivariate normal distributions each with covariance
matrix of the following form:
%
%
%eB.4 #&#
\begin{equation}
\pmatrix{ 1 & \mathbf{r}_{1\times( p-1 ) }
\vspace*{2pt}\cr
\mathbf{r}_{ ( p-1 ) \times1} & \mathbf{S}_{ ( p-1 )
\times( p-1 ) }} + ( a_{l}-1 ) \bI_{p}\qquad\mbox{for }l=1,2,\ldots, n,
\label{matrixform}
\end{equation}
where $\llVert\mathbf{r}\rrVert_{0}=k$ with nonzero elements
of $r$
equal to $\epsilon_{n,p}$ and the submatrix $\mathbf{S}_{(p-1)\times(p-1)}$
is the same as the one for $\Sigma_{l,0}$ given in (\ref{sigma0}). Note
that the indices $\gamma_{i}$ and $\lambda_{i}$ are dropped from
$\mathbf{r%
}$ and $\mathbf{S}$ to simplify the notation.

With Lemma~\ref{chisquareslem} in place, it remains to establish equation
(\ref{chi-squarebd}) in order to prove Proposition~\ref{affbd}. The
following lemma
is useful for calculating the cross product terms in the chi-square distance
between Gaussian mixtures. The proof of the lemma is straightforward
and is thus omitted.

%
%le16 #&#
\begin{lemma}
\label{crossproduct} Let $g_{i}$ be the density function of $N (
0,\Sigma_{i} ) $ for $i=0,1$ and $2$, respectively. Then
\[
\int\frac{g_{1}g_{2}}{g_{0}}=\frac{1}{ [ \det( \bI-\Sigma
_{0}^{-2} ( \Sigma_{1}-\Sigma_{0} ) ( \Sigma_{2}-\Sigma
_{0} ) ) ] ^{1/2}}.
\]
\end{lemma}

Let $\Sigma_{l,i}$, $i=1$ or
$2$, be two covariance matrices of the form (\ref{matrixform}). Note
that $%
\Sigma_{l,i}$, $i=0,1$ or $2$, differs from each other only in the first
row/column. Then $\Sigma_{l,i}-\Sigma_{l,0}$, $i=1$ or $2$, has a very
simple structure. The nonzero elements only appear in the first row/column,
and in total there are $2k$ nonzero elements. This property immediately implies
the following lemma which makes the problem of studying the determinant
in Lemma~\ref{crossproduct} relatively easy.

%
%le17 #&#
\begin{lemma}
\label{sigmai} Let $\Sigma_{l,i}$, $i=1$ and $2$, be matrices of the
form (%
\ref{matrixform}). Define $J$ to be the number of overlapping
$\epsilon
_{n,p}$'s between $\Sigma_{l,1}$ and $\Sigma_{l,2}$ on the first row, and
\[
Q\stackrel{\bigtriangleup} {=} ( q_{ij} ) _{1\leq i,j\leq
p}= (
\Sigma_{l,1}-\Sigma_{l,0} ) ( \Sigma_{l,2}-\Sigma
_{l,0} ).
\]
There are index subsets $I_{r}$ and $I_{c}$ in $ \{ 1,2,\ldots,p \} $ with $\mathrm{Card} ( I_{r} ) =\mathrm{Card} (
I_{c} ) =k$ and $\mathrm{Card} ( I_{r}\cap I_{c} ) =J$ such
that
\[
q_{ij}= \cases{ %
J\epsilon_{n,p}^{2},
&\quad $i=j=1$,
\vspace*{2pt}\cr
\epsilon_{n,p}^{2}, &\quad $i\in I_{r}\mbox{ and }j\in
I_{c}$,
\vspace*{2pt}\cr
0, &\quad $\mbox{otherwise,}$}
\]
and the matrix $ ( \Sigma_{l,0}-\Sigma_{l,1} ) ( \Sigma
_{l,0}-\Sigma_{l,2} ) $ has rank $2$ with two identical nonzero
eigenvalues $J\epsilon_{n,p}^{2}$ when $J>0$.
\end{lemma}

Let
%
%
%eB.5 #&#
\begin{equation}
R_{l,\lambda_{1},\lambda_{1}^{
%TCIMACRO{\U{b4}}%
%BeginExpansion
\prime%
%EndExpansion
}}^{\gamma_{-1},\lambda_{-1}}=-\log\det\bigl( I-\Sigma_{l,0}^{-2}
( \Sigma_{l,0}-\Sigma_{l,1} ) ( \Sigma_{l,0}-\Sigma
_{l,2} ) \bigr), \label{R}
\end{equation}
where $\Sigma_{l,0}$ is defined in (\ref{sigma0}) and determined by
$ (
\gamma_{-1},\lambda_{-1} ) $, and $\Sigma_{l,1}$ and $\Sigma_{l,2}$
have the first row $\lambda_{1}$ and $\lambda_{1}^{%
%TCIMACRO{\U{b4}}%
%BeginExpansion
\prime%
%EndExpansion
}$, respectively. We drop the indices $\lambda_{1}$, $\lambda_{1}^{%
%TCIMACRO{\U{b4}}%
%BeginExpansion
\prime%
%EndExpansion
}$ and $ ( \gamma_{-1},\lambda_{-1} ) $ from $\Sigma_{i}$ to
simplify the notation. Define
\begin{eqnarray*}
\Theta_{-1} ( a_{1},a_{2} ) &=& \bigl\{ ( b,c )
\dvtx\mbox{ %%
there exist }\theta_{i}\in\Theta\mbox{,
}i=1, 2, \mbox{, such that }%
\lambda_{1}(
\theta_{i})=a_{i}\\
&&\hspace*{185pt}\mbox{and }\lambda_{-1}(
\theta_{i})=c \bigr\}.
\end{eqnarray*}
It is a subset of $\Theta_{-1}$ in which the element can pick both $a_{1}$
and $a_{2}$ as the first row to form parameters in $\Theta$. From
Lemma~\ref%
{crossproduct} the left-hand side of equation (\ref{chi-squarebd}) can be
written as
\begin{eqnarray} \label{chisquare}
&&\mathop{\mathrm{Average}}_{ ( \gamma_{-1},\lambda_{-1} ) \in
\Theta_{-1}} \Biggl\{ \mathop{\mathrm{Average}}_{\lambda_{1},\lambda_{1}^{
%%
%TCIMACRO{\U{b4}}%
%BeginExpansion
\prime%
%EndExpansion
}\in\Lambda_{1} ( \lambda_{-1} ) }
\Biggl[ \exp\Biggl( \frac{1}{2}\sum_{l=1}^{n}R_{l,\lambda
_{1},\lambda_{1}^{%
%TCIMACRO{\U{b4}}%
%BeginExpansion
\prime%
%EndExpansion
}}^{\gamma_{-1},\lambda_{-1}}
\Biggr)-1 \Biggr] \Biggr\}
\nonumber
\\[-8pt]
\\[-8pt]
\nonumber
&&\qquad=\mathop{\mathrm{Average}}_{\lambda_{1},\lambda_{1}^{%
%TCIMACRO{\U{b4}}%
%BeginExpansion
\prime%
%EndExpansion
}\in B} \Biggl\{ \mathop{\mathrm{Average}}_{ ( \gamma_{-1},\lambda
_{-1} ) \in\Theta_{-1} ( \lambda_{1},\lambda_{1}^{%
%TCIMACRO{\U{b4}}%
%BeginExpansion
\prime%
%EndExpansion
} ) }
\Biggl[ \exp\Biggl(\frac{1}{2}\sum_{l=1}^{n}R_{l,%
\lambda_{1},\lambda_{1}^{%
%TCIMACRO{\U{b4}}%
%BeginExpansion
\prime%
%EndExpansion
}}^{\gamma_{-1},\lambda_{-1}}
\Biggr)-1 \Biggr] \Biggr\},
\nonumber
\end{eqnarray}
where $B$ is defined in step 1.

Lemmas~\ref{sigmai} and~\ref{Rlem} below show that $R_{l,\lambda
_{1},\lambda_{1}^{%
%TCIMACRO{\U{b4}}%
%BeginExpansion
\prime%
%EndExpansion
}}^{\gamma_{-1},\lambda_{-1}}$ is approximately equal to
\[
-\log\det\bigl( I-a_{l}^{-2} ( \Sigma_{l,0}-
\Sigma_{l,1} ) ( \Sigma_{l,0}-\Sigma_{l,2} ) \bigr)
=-2\log\bigl( 1-a_{l}^{-2}J\epsilon_{n,p}^{2}
\bigr).
\]
Define
\begin{eqnarray*}
&&\Lambda_{1,J}= \bigl\{ \bigl( \lambda_{1},
\lambda_{1}^{\prime} \bigr) \in\Lambda_{1}\otimes
\Lambda_{1}\dvtx\mbox{the number of overlapping }\epsilon
_{n,p}\mbox{'s between }\lambda_{1}\\
&&\hspace*{282pt}\mbox{and
}\lambda_{1}^{\prime
}\mbox{ is }J \bigr\}.
\end{eqnarray*}

%
%le18 #&#
\begin{lemma}
\label{Rlem}For $R_{l,\lambda_{1},\lambda_{1}^{%
%TCIMACRO{\U{b4}}%
%BeginExpansion
\prime%
%EndExpansion
}}^{\gamma_{-1},\lambda_{-1}}$ defined in equation (\ref{R}), we have
%
%
%eB.6 #&#
\begin{equation}
R_{l,\lambda_{1},\lambda_{1}^{%
%TCIMACRO{\U{b4}}%
%BeginExpansion
\prime%
%EndExpansion
}}^{\gamma_{-1},\lambda_{-1}}=-2\log\bigl( 1-Ja_{l}^{-2}
\epsilon_{n,p}^{2} \bigr) +\delta_{l,\lambda_{1},\lambda_{1}^{%
%TCIMACRO{\U{b4}}%
%BeginExpansion
\prime%
%EndExpansion
}}^{\gamma_{-1},\lambda_{-1}},
\label{Rdecomp}
\end{equation}
where $\delta_{l,\lambda_{1},\lambda_{1}^{%
%TCIMACRO{\U{b4}}%
%BeginExpansion
\prime%
%EndExpansion
}}^{\gamma_{-1},\lambda_{-1}}$ satisfies
%
%
%eB.7 #&#
\begin{equation}
\mathop{\mathrm{Average}}_{ ( \lambda_{1},\lambda_{1}^{\prime} )
\in\Lambda
_{1,J}} \Biggl[ \mathop{\mathrm{Average}}_{ ( \gamma_{-1},\lambda
_{-1} ) \in\Theta_{-1} ( \lambda_{1},\lambda_{1}^{%
%TCIMACRO{\U{b4}}%
%BeginExpansion
\prime%
%EndExpansion
} ) }
\exp\Biggl( \frac{1}{2}\sum_{l=1}^{n}
\delta_{l,\lambda_{1},\lambda_{1}^{%
%TCIMACRO{\U{b4}}%
%BeginExpansion
\prime%
%EndExpansion
}}^{\gamma_{-1},\lambda_{-1}} \Biggr) \Biggr] \leq3/2, \label{R1eq}
\end{equation}
uniformly over all $J$ defined in Lemma~\ref{sigmai}.
\end{lemma}

We will prove Lemma~\ref{Rlem} in Section~\ref{secRterms}.

%sB.2 #&#
\subsection{\texorpdfstring{Proof of equation (\protect\ref{chi-squarebd})}{Proof of equation (71)}}
\label{chisquaresec}

We are now ready to establish equation (\ref{chi-squarebd}) using Lemma~\ref%
{Rlem}. %which is the key step in proving Lemma~\ref{affbd}.
It follows from equation (\ref{Rdecomp}) in Lemma~\ref{Rlem} that
\begin{eqnarray*}\hspace*{-5pt}
& &\mathop{\mathrm{Average}}_{\lambda_{1},\lambda_{1}^{%
%TCIMACRO{\U{b4}}%
%BeginExpansion
\prime%
%EndExpansion
}\in B} \Biggl\{\mathop{\mathrm{Average}} _{ ( \gamma_{-1},\lambda
_{-1} ) \in\Theta_{-1} ( \lambda_{1},\lambda_{1}^{%
%TCIMACRO{\U{b4}}%
%BeginExpansion
\prime%
%EndExpansion
} ) }
\Biggl[ \exp\Biggl( \frac{1}{2}%
\sum
_{l=1}^{n}R_{l,\lambda_{1},\lambda_{1}^{%
%TCIMACRO{\U{b4}}
%BeginExpansion
\prime
%EndExpansion
}}^{\gamma_{-1},\lambda_{-1}} \Biggr)
-1 \Biggr] \Biggr\} \\\hspace*{-5pt}
&&\qquad=\mathop{\mathrm{Average}_{J}} \Biggl\{ -\sum_{l=1}^{n}
\log\biggl( 1- \frac{J\epsilon_{n,p}^{2}}{a_{l}^{2}} \biggr) \\\hspace*{-5pt}
&&\hspace*{80pt}{}\times
\mathop{\mathrm{Average}}
_{ ( \lambda
_{1},\lambda_{1}^{\prime} ) \in\Lambda_{1,J}} \Biggl[\mathop{\mathrm
{Average}} _{ ( \gamma_{-1},\lambda_{-1} ) \in
\Theta
_{-1} ( \lambda_{1},\lambda_{1}^{
%TCIMACRO{\U{b4}}%
%BeginExpansion
\prime%
%EndExpansion
} ) }\exp\Biggl(
\frac{1}{2}\sum_{l=1}^{n}\delta
_{l,\lambda_{1},\lambda_{1}^{
%TCIMACRO{\U{b4}}%
%BeginExpansion
\prime%
%EndExpansion
}}^{\gamma_{-1},\lambda_{-1}} \Biggr) \Biggr] -1 \Biggr\}.
\end{eqnarray*}
Recall that $J$ is the number of overlapping $\epsilon_{n,p}$'s
between $%
\Sigma_{l,1}$ and $\Sigma_{l,2}$ on the first row. It is easy to see
that $%
J$ has the hypergeometric distribution with
%
%
%eB.8 #&#
\begin{eqnarray}\label{J}
&&\mathbb{P} \bigl( \mbox{number of overlapping }\epsilon_{n,p}
\mbox{'s}%
=J \bigr)\nonumber\\[-2pt] &&\qquad= \pmatrix{k
\cr
J}
\pmatrix{p_{\lambda_{-1}}-k
\cr
k-J} / \pmatrix{p_{\lambda
_{-1}}
\cr
k}
\\[-2pt]
&&\qquad\leq
\biggl( \frac{k^{2}}{p_{\lambda_{-1}}-k} \biggr) ^{J}.\nonumber
\end{eqnarray}
Equations (\ref{R1eq}) and (\ref{J}) imply
\begin{eqnarray*}
&&\mathop{\mathrm{Average}}_{ ( \gamma_{-1},\lambda_{-1} ) \in
\Theta_{-1}} \biggl\{ \int\biggl( \frac{d\bar{\mathbb{P}}_{ (
1,1,\gamma_{-1},\lambda_{-1} ) }}{d\bar{\mathbb{P}}_{ (
1,0,\gamma_{-1},\lambda_{-1} ) }}
\biggr) ^{2}\,d\bar{\mathbb{P}}%
_{ ( 1,0,\gamma_{-1},\lambda_{-1} ) }-1 \biggr\}
\\[-2pt]
&&\qquad\leq\sum_{J\geq0}\frac{{k\choose J}{p_{\lambda
_{-1}}-k\choose k-J}}{%
{p_{\lambda_{-1}}\choose k}} \Biggl\{ -\sum
_{l=1}^{n}\log\bigl( 1-J\epsilon
_{n,p}^{2}/a_{l}^{2} \bigr) \Biggr\}
\frac{3}{2}-1
\\[-2pt]
&&\qquad\leq C\sum_{J\geq1} \bigl( p^{{(\beta-1)}/{\beta}} \bigr)
^{-J}\exp\Biggl( 2J\sum_{l=1}^{n}a_{l}^{-2}
\cdot\frac{\upsilon^{2}\log p}{\sqrt{n}}%
\Biggr) +1/2
\\[-2pt]
&&\qquad\leq C\sum_{J\geq1} \bigl( p^{{(\beta-1)}/{\beta}} \bigr)
^{-J}\exp\biggl( 2 J c_{\kappa} \sqrt{n}\cdot
\frac{\upsilon^{2}\log p}{\sqrt{n}} \biggr) +1/2
\\[-2pt]
&&\qquad\leq C\sum_{J\geq1} \bigl( p^{{(\beta-1)}/{\beta}} \bigr)
^{-J}\exp\bigl( 2c_{\kappa}J\upsilon^{2}\log p
\bigr) +1/2\\[-2pt]
&&\qquad\leq C\sum_{J\geq
1} \bigl( p^{{(\beta-1)}/{(2\beta)}}
\bigr) ^{-J}+1/2<C_{2}^{2},
\end{eqnarray*}
where the third inequality is from (\ref{ck}), the fifth inequality is due
to (\ref{v1}) and the last inequality is obtained by setting
$C_{2}^{2}=3/4$. %\hbox{\vrule width
%4pt height 6pt depth 1.5pt}

%sB.3 #&#
\subsection{\texorpdfstring{Proof of Lemma \protect\ref{Rlem}}{Proof of Lemma 18}}
\label{secRterms}

Define
%
%
%eB.9 #&#
\begin{eqnarray}
\label{Adef} A_{l}&=& \bigl[ I-a_{l}^{-2} (
\Sigma_{l,0}-\Sigma_{l,1} ) ( \Sigma_{l,0}-
\Sigma_{l,2} ) \bigr] ^{-1}
\nonumber
\\[-8pt]
\\[-8pt]
\nonumber
&&{}\times \bigl( a_{l}^{2}
( \Sigma_{l,0} ) ^{-2}-I \bigr) a_{l}^{-2}
( \Sigma_{l,0}-\Sigma_{l,1} ) ( \Sigma_{l,0}-
\Sigma_{l,2} )
\end{eqnarray}
and
\[
\delta_{l,\lambda_{1},\lambda_{1}^{%
%TCIMACRO{\U{b4}}%
%BeginExpansion
\prime%
%EndExpansion
}}^{\gamma_{-1},\lambda_{-1}}=-\log\det( I-A_{l} ).\
\]
We rewrite $R_{l,\lambda_{1},\lambda_{1}^{\prime}}^{\gamma
_{-1},\lambda_{-1}}$ as follows:
%
%
%eB.10 #&#
%eB.11 #&#
\begin{eqnarray}
\label{minsR}  R_{l,\lambda_{1},\lambda_{1}^{\prime}}^{\gamma
_{-1},\lambda
_{-1}}& =& -\log\det\bigl[ I-a_{l}^{-2} ( \Sigma_{l,0}-
\Sigma_{l,1} ) ( \Sigma_{l,0}-\Sigma_{l,2} ) \nonumber\\
&&\hspace*{43pt}{}-
\bigl( a_{l}^{2}\Sigma_{l,0}^{-2}-I
\bigr) a_{l}^{-2} ( \Sigma_{l,0}-
\Sigma_{l,1} ) ( \Sigma_{l,0}-\Sigma_{l,2} ) \bigr]
\nonumber
\\
& = &-\log\det\bigl\{ [ I-A_{l} ] \cdot\bigl[ I-a_{l}^{-2}
( \Sigma_{l,0}-\Sigma_{l,1} ) ( \Sigma_{l,0}-\Sigma
_{l,2} ) \bigr] \bigr\}
\\
&=& -\log\det\bigl[ I-a_{l}^{-2} ( \Sigma_{l,0}-
\Sigma_{l,1} ) ( \Sigma_{l,0}-\Sigma_{l,2} )
\bigr] -\log\det( I-A_{l} )
\nonumber
\\
&=& -2\log\bigl( 1-J\epsilon_{n,p}^{2}/a_{l}^{2}
\bigr) +\delta_{l,\lambda
_{1},\lambda_{1}^{
%TCIMACRO{\U{b4}}%
%BeginExpansion
\prime%
%EndExpansion
}}^{\gamma_{-1},\lambda_{-1}},\nonumber
\end{eqnarray}
where the last equation follows from Lemma~\ref{sigmai}.

Now we are ready to establish equation (\ref{R1eq}). For simplicity we
will write matrix norm
$\| \cdot\|_2$ as $\| \cdot\|$ below. It is important to observe that
$\mathrm{rank} ( A_{l} ) \leq2$ due to the simple structure
of $%
( \Sigma_{l,0}-\Sigma_{l,1} ) ( \Sigma_{l,0}-\Sigma
_{l,2} ) $. Let $\varrho_{l}$ be an eigenvalue of $A_{l}$. It is easy
to see that
%
%
%eB.12 #&#
\begin{eqnarray}
\label{Anorm}  \vert\varrho_{l}\vert&\leq&\llVert
A_{l}\rrVert
\nonumber
\\
& \leq&\bigl\llVert a_{l}^{2}\Sigma_{l,0}^{-2}-I
\bigr\rrVert\cdot a_{l}^{-2}\llVert\Sigma_{l,0}-
\Sigma_{l,1}\rrVert\llVert\Sigma_{l,0}-
\Sigma_{l,2} \rrVert
\nonumber
\\[-8pt]
\\[-8pt]
\nonumber
&&{}/ \bigl( 1-a_{l}^{-2}\llVert
\Sigma_{l,0}-\Sigma_{l,1}\rrVert\llVert
\Sigma_{l,0}- \Sigma_{l,2}\rrVert\bigr)
\nonumber
\\
& \leq&\biggl( \biggl( \frac{3}{2} \biggr) ^{2}-1 \biggr)
\frac{1}{3}\cdot\frac{%
1}{3} \bigg/ \biggl( 1-\frac{1}{3}\cdot
\frac{1}{3} \biggr) =5/32<1/6,\nonumber
\end{eqnarray}
since $\llVert a_{l}^{-1} ( \Sigma_{l,0}-\Sigma_{l,1} )
\rrVert\leq\llVert a_{l}^{-1} ( \Sigma_{l,0}-\Sigma
_{l,1} ) \rrVert_{1}=2k\epsilon_{n,p}<1/3$ and\break $\lambda
_{\min
} ( a_{l}^{-1}\Sigma_{l,0} ) \geq1-\llVert
I-a_{l}^{-1}\Sigma
_{l,0}\rrVert\geq1-\llVert I-a_{l}^{-1}\Sigma_{l,0}\rrVert
_{1}>2/3$ from equation~(\ref{1normbound}).

Note that (\ref{Anorm}) and
\[
\bigl\vert\log( 1-x ) \bigr\vert\leq2\vert x\vert,\qquad \mbox{for }
\vert x\vert<1/6,
\]
imply
\[
\delta_{l,\lambda_{1},\lambda_{1}^{%
%TCIMACRO{\U{b4}}%
%BeginExpansion
\prime%
%EndExpansion
}}^{\gamma_{-1},\lambda_{-1}}\leq4\llVert A_{l}\rrVert,
\]
and then
%
%
%eB.13 #&#
\begin{equation}
\label{deltaA} \exp\Biggl( \frac{1}{2}\sum_{l=1}^{n}
\delta_{l,\lambda_{1},\lambda
_{1}^{%
%TCIMACRO{\U{b4}}%
%BeginExpansion
\prime%
%EndExpansion
}}^{\gamma_{-1},\lambda_{-1}} \Biggr) \leq\exp\Biggl( 2\sum
_{l=1}^{n}\llVert A_{l}\rrVert\Biggr).
\end{equation}

Since
%
%
%eB.14 #&#
\begin{equation}
\label{sigma0norm} \cases{ %
\bigl\llVert
I-a_{l}^{-1}\Sigma_{l,0} \bigr\rrVert\leq\bigl
\llVert I-a_{l}^{-1}\Sigma_{l,0} \bigr\rrVert
_{1}=2k\epsilon_{n,p}<1/3<1,
\vspace*{2pt}\cr
\bigl\llVert a_{l}^{-2} ( \Sigma_{l,0}-
\Sigma_{l,1} ) ( \Sigma_{l,0}-\Sigma_{l,2} ) \bigr
\rrVert\leq\frac{1}{3}\cdot\frac{1}{3}<1,}
\end{equation}
%
%I-a_{l}^{-1}\Sigma_{l,0}\right\Vert_{1}=2k\epsilon_{n,p}<1/3<1,
we write%
%
%
%eB.15 #&#
\begin{eqnarray} \label{Asigma0}
a_{l}^{2}\Sigma_{l,0}^{-2}-I &=&
\bigl( I- \bigl( I-a_{l}^{-1}\Sigma_{l,0} \bigr)
\bigr) ^{-2}-I\nonumber\\
&=& \biggl( I+\sum_{k=1} \bigl(
I-a_{l}^{-1}\Sigma_{l,0} \bigr) ^{k}
\biggr) ^{2}-I
\\
&=& \Biggl[ \sum_{m=0}^{\infty} ( m+2 )
\bigl( I-a_{l}^{-1}\Sigma_{l,0} \bigr)
^{m} \Biggr] \bigl( I-a_{l}^{-1}
\Sigma_{l,0} \bigr),
\nonumber
\end{eqnarray}
where
%
%
%eB.16 #&#
\begin{equation}
\label{Asigma0-series} \Biggl\llVert\sum_{m=0}^{\infty}
( m+2 ) \bigl( I-a_{l}^{-1}\Sigma_{l,0} \bigr)
^{m} \Biggr\rrVert\leq\sum_{m=0}^{\infty}
( m+2 ) \biggl( \frac{1}{3} \biggr) ^{m}<3.
\end{equation}
Define%
%
%
%eB.17 #&#
\begin{equation}
A_{l\ast}= \bigl( I-a_{l}^{-1}
\Sigma_{l,0} \bigr) \cdot a_{l}^{-2} (
\Sigma_{l,0}-\Sigma_{l,1} ) ( \Sigma_{l,0}-\Sigma
_{l,2} ). \label{Astardef}
\end{equation}
From equations (\ref{Adef}) and %(\ref{Asigma0})
(\ref{sigma0norm})--(\ref{Astardef}) we have
\begin{eqnarray*}
\llVert A_{l}\rrVert&\leq& \bigl\llVert\bigl[ I-a_{l}^{-2}
( \Sigma_{l,0}-\Sigma_{l,1} ) ( \Sigma_{l,0}-\Sigma
_{l,2} ) %
\bigr] ^{-1} \bigr\rrVert\\
&&{}\times\Biggl\llVert
\sum_{m=0}^{\infty} ( m+2 ) \bigl(
I-a_{l}^{-1}\Sigma_{l,0} \bigr) ^{m}
\Biggr\rrVert\llVert A_{l\ast
}\rrVert
\\
&<& \frac{1}{1-({1}/{3})\cdot({1}/{3})} \cdot3 \cdot\llVert
A_{l\ast
}\rrVert=
\frac{27}{8}\llVert A_{l\ast}\rrVert\leq\frac
{27}{8}%
\max\bigl\{ \llVert A_{l\ast}\rrVert_{1},\llVert
A_{l\ast
}\rrVert_{\infty} \bigr\}.
\end{eqnarray*}
The above result and (\ref{deltaA}) indicate that the proof of Lemma~\ref{Rlem} is complete if
we show
%It is then sufficient to show
%
%
%eB.18 #&#
\begin{eqnarray}
\label{A1infty}&& \mathop{\mathrm{Average}}_{ ( \lambda_{1},\lambda
_{1}^{\prime} )
\in\Lambda
_{1,J}} \Biggl[ {
\mathop\mathrm{Average}}_{ ( \gamma_{-1},\lambda
_{-1} ) \in\Theta_{-1} ( \lambda_{1},\lambda_{1}^{%
%TCIMACRO{\U{b4}}%
%BeginExpansion
\prime%
%EndExpansion
} ) }
\nonumber
\\[-8pt]
\\[-8pt]
\nonumber
&&\hspace*{51pt}{}\times\exp\Biggl( \frac{27}{2}\sum
_{l=1}^{n}\max\bigl\{ \llVert A_{l\ast}
\rrVert_{1},\llVert A_{l\ast
}\rrVert_{\infty} \bigr\}
\Biggr) \Biggr] \leq3/2,
\end{eqnarray}
where $\llVert A_{l\ast}\rrVert_1 $ and $\llVert A_{l\ast
}\rrVert_\infty$
depend on the values of $\lambda_{1},\lambda_{1}^{
%TCIMACRO{\U{b4}}%
%BeginExpansion
\prime
%EndExpansion
}$ and $ ( \gamma_{-1},\lambda_{-1} ) $. We dropped the
indices $%
\lambda_{1}$, $\lambda_{1}^{%
%TCIMACRO{\U{b4}}%
%BeginExpansion
\prime%
%EndExpansion
}$ and $ ( \gamma_{-1},\lambda_{-1} ) $ from $A_{l}$ to simplify
the notation.

Let $E_{m}= \{ 1,2,\ldots,r \} / \{ 1,m \} $. Let $%
n_{\lambda_{E_{m}}}$ be the number of columns of $\lambda_{E_{m}}$ with
column sum at least $2k-2$ for which two rows cannot freely take value $0$
or $1$ in this column. Then we have $p_{\lambda_{E_{m}}}=r-n_{\lambda
_{E_{m}}}$. Without loss of generality we assume that $k\geq3$. Since $
n_{\lambda_{E_{m}}}\cdot( 2k-2 ) \leq r\cdot k$, the total
number of $1$s in the upper triangular matrix by the construction of the
parameter set, we thus have $n_{\lambda_{E_{m}}}\leq r\cdot\frac{3}{4}$,
which immediately implies $p_{\lambda_{E_{m}}}=r-n_{\lambda
_{E_{m}}}\geq
\frac{r}{4}\geq p/8-1$. Thus we have for every nonnegative integer~$t$,
\begin{eqnarray*}
&&\mathbb{P} \bigl( %\left\Vert A_{l\ast}\right\Vert
\max\bigl\{ \llVert A_{l\ast}\rrVert
_{1},\llVert A_{l\ast
}\rrVert_{\infty} \bigr\} \geq2t
\cdot\epsilon_{n,p}\cdot k\epsilon_{n,p}^{2}\cdot
a_{l}^{-3} \bigr)
\\
&&\qquad\leq\mathbb{P} \bigl( \llVert A_{l\ast}\rrVert_{1}\geq2t
\cdot\epsilon_{n,p}\cdot k\epsilon_{n,p}^{2}\cdot
a_{l}^{-3} \bigr) +\mathbb{P}%
\bigl( \llVert
A_{l\ast}\rrVert_{\infty}\geq2t\cdot\epsilon_{n,p}
\cdot k\epsilon_{n,p}^{2}\cdot a_{l}^{-3}
\bigr)
\\
&&\qquad\leq 2\sum_{m}\mathop{\mathrm{Average}}_{\lambda_{E_{m}}}
\frac{{k\choose%
t}{p_{\lambda_{E_{m}}}\choose{k-t}}}{{p_{\lambda_{E_{m}}}\choose
k}}\leq2p \biggl( \frac{k^{2}}{p/8-1-k} \biggr) ^{t}
\end{eqnarray*}
from equation (\ref{J}), which immediately implies
%
%
%eB.19 #&#
\begin{eqnarray}\label{Astarbound}
&&\hspace*{-4pt}\mathop{\mathrm{Average}}_{ ( \lambda_{1},\lambda_{1}^{\prime} )
\in\Lambda
_{1,J}} \Biggl[\mathop{\mathrm{Average}} _{ ( \gamma_{-1},\lambda
_{-1} ) \in\Theta_{-1} ( \lambda_{1},\lambda_{1}^{%
%TCIMACRO{\U{b4}}%
%BeginExpansion
\prime%
%EndExpansion
} ) }
\exp\Biggl( \frac{27}{2} \sum_{l=1}^{n}
\max\bigl\{ \llVert A_{l\ast}\rrVert
_{1},\llVert A_{l\ast
}\rrVert_{\infty} \bigr\}
\Biggr) \Biggr]
\nonumber
\\
&&\hspace*{-4pt}\qquad\leq\exp\Biggl( \frac{27}{2}\sum_{l=1}^{n}
\frac{4\beta}{\beta
-1}\cdot\epsilon_{n,p}\cdot k\epsilon_{n,p}^{2}
\cdot a_{l}^{-3} \Biggr)
\nonumber
\\
&&\hspace*{-4pt} \quad\qquad{}+\int_{{2\beta}/{(\beta-1)}}^{\infty} \Biggl( 27k\epsilon
_{n,p}^{3}\sum_{l=1}^{n}a_{l}^{-3}
\Biggr)\nonumber\\
&&\hspace*{81pt}{}\times \exp\Biggl( \frac{27}{2}%
\sum
_{l=1}^{n}2t\cdot\epsilon_{n,p}\cdot k
\epsilon_{n,p}^{2}\cdot a_{l}^{-3}
\Biggr) 2p \biggl( \frac{k^{2}}{p/8-1-k} \biggr) ^{t-1}\,dt
\\
&&\hspace*{-4pt}\qquad\leq\exp\Biggl( 54 \cdot\Biggl( \sum_{l=1}^{n}a_{l}^{-3}
\Biggr) \cdot\frac{\beta}{\beta-1}\cdot k \epsilon_{n,p}^3
\Biggr)
\nonumber
\\
&&\hspace*{-4pt}\qquad\quad{}+ 2 p \int_{{2\beta}/{(\beta-1)}}^{\infty}\exp\Biggl[ %\log
( t+1 ) \cdot27 %\frac{27}{2}
\Biggl(\sum
_{l=1}^{n}a_{l}^{-3} \Biggr)
k \epsilon_{n,p}^3 %\pi_{n}
\nonumber\\
&&\hspace*{-4pt}\hspace*{90pt}\qquad\quad{}- ( t-1 )
\log\frac{p/8-1-k}{k^{2}} \Biggr] \,dt.\nonumber
\end{eqnarray}
Note that (\ref{ck}) implies
\[
\sum_{l=1}^{n}a_{l}^{-3}
\leq\sum_{l=1}^{n}a_{l}^{-2}
\leq c_{\kappa
}\sqrt{n},
\]
using (\ref{condc}) and (\ref{k-epsilon}) we have
\begin{eqnarray*}
2 \sqrt{n} k \epsilon_{n,p}^3 &\leq&\sqrt{n}
\pi_{n}(p)\epsilon_{n,p}^{3-q}\\
&\leq&\aleph
v^{3-q}n^{1/2}n^{(1-q)/4} (\log p)^{(q-3)/2}
n^{(q-3)/4} (\log p)^{(3-q)/2}\\
&=&\aleph v^{3-q},
\end{eqnarray*}
and thus we can bound the first term on the right-hand side of (\ref%
{Astarbound}),
\[
\exp\biggl( 54 \cdot c_{\kappa}\sqrt{n}\cdot\frac{\beta}{\beta
-1}\cdot k
\epsilon_{n,p}^{3} \biggr) \leq\exp\biggl(
\frac{
\beta}{\beta-1}\cdot27c_{\kappa}v^{2}\cdot\aleph
v^{1-q} \biggr) \leq\exp( 1/3 ) <3/2,
\]
where the second inequality is from (\ref{v}) and (\ref{v1}). We will show
that the second term on the right-hand side of (\ref{Astarbound}) is
negligible and hence establish~(\ref{A1infty}). Indeed, since we have just
shown that
\[
27 \Biggl( \sum_{l=1}^{n}a_{l}^{-3}
\Biggr) k \epsilon_{n,p}^{3}\leq\frac{\beta-1}{6\beta},
\]
the second term on the right-hand side of (\ref{Astarbound}) is bounded
by
\begin{eqnarray*}
&& 2p \int_{{2\beta}/{(\beta-1)}}^{\infty}\exp\biggl[ ( t+1 )
\frac{\beta-1}{6\beta}- ( t-1 ) \log\frac
{p/8-1-k}{%
k^{2}} \biggr] \,dt
\\
&&\qquad= 2 \biggl( \log\frac{p/8-1-k}{k^{2}}-\frac{\beta-1}{6\beta}
\biggr) ^{-1}
\\
&&\qquad\quad{}\times \exp\biggl[ \log p+ \biggl( \frac{2\beta}{\beta-1}+1 \biggr)
\frac
{\beta
-1}{6\beta}-
\biggl( \frac{2\beta}{\beta-1}-1 \biggr) \log\frac
{p/8-1-k}{%
k^{2}} \biggr]
\\
&&\qquad= O \bigl( p^{-1/\beta} [\log p]^{6 /(\beta-1)+2} \bigr) =o ( 1 ),
\end{eqnarray*}
where the second equality is from the fact that (\ref{condc}) and
(\ref
{k-epsilon}) together with $p\geq n^{\beta/2}$ indicate
\[
k^{2}\leq\pi_{n}(p)\epsilon_{n,p}^{-2q}/4
\leq\frac{\aleph
v^{-2q}\sqrt{n}}{%
4\log^{3}p}\leq\frac{\aleph v^{-2q}p^{1/\beta}}{4\log^{3}p},
\]
and then
\begin{eqnarray*}
&&\biggl( \frac{2\beta}{\beta-1}-1 \biggr) \log\frac{p/8-k}{k^{2}}
\\
&&\qquad =  \biggl(
\frac{2\beta}{\beta-1}-1 \biggr) \log\bigl( pk^{-2} \bigr) \bigl
[1 + o(1)
\bigr]
\\
&&\qquad\geq \biggl( \frac{2\beta}{\beta-1}-1 \biggr) \biggl[ \frac
{\beta
-1}{\beta
}\log p+3
\log\log p-\log\bigl(Mv^{-2q}/4 \bigr) \biggr]
\\
&&\qquad= \biggl( 1+\frac{1}{\beta} \biggr) \log p \bigl[1 + o(1) \bigr
]. %\hbox{\vrule width 4pt height 6pt depth
%1.5pt}
\end{eqnarray*}
%
%- k/p) \sim\log\left( p k^{-2} \right). \]
%that $\upsilon^{2}<\frac{\beta-1}{27 c_{\kappa}\beta}$ and $k^{2}=O
%}\right) $ as defined in Section~\ref{lowerboundproofssc} %
%Therefore we have
%%TCIMACRO{\U{b4}}%
%%BeginExpansion
%{\acute{}}%
%%EndExpansion
%}\right) }{\mathrm{Average}}\exp\left( \frac{27}{2}\sum_{l=1}^{n}\left
%A_{l\ast}\right\Vert\right) \leq3/2.
%$\ $\hbox{\vrule width 4pt height 6pt depth 1.5pt}
\end{appendix}

% zodis "Acknowledgments" paliekamas pagal autoriu

%suskaldyti doi

% imsref loaded by akundreckaite, 2013-07-30 16:20:42
%

\printaddresses

\end{document}